\documentclass[a4paper, 10pt]{article}

\usepackage[top=2.5cm, bottom=2.5cm, left=2.5cm, right=2.5cm]{geometry}
\usepackage[english]{babel}
\usepackage[utf8]{inputenc}
\usepackage[T1]{fontenc}
\usepackage{lmodern}
\usepackage{graphicx}
\usepackage{amsmath, amssymb, mathrsfs, mathtools, cancel}
\usepackage{xcolor}
\usepackage{environ, fancybox, framed}
\usepackage[thmmarks, amsmath]{ntheorem}
\usepackage[framemethod=TikZ]{mdframed}
\usepackage{tikz}
\usetikzlibrary{patterns.meta}
\usetikzlibrary{arrows.meta}
\usetikzlibrary{decorations.pathreplacing}
\usepackage[justification=centering]{caption}
\usepackage{enumitem, multirow}
\usepackage{multicol, contour}
\usepackage[ntheorem]{empheq}
\usepackage{biblatex}
\AtBeginBibliography{\small}
\usepackage[hidelinks, colorlinks=true, linkcolor=blue, urlcolor=cyan, citecolor=green!75!black]{hyperref}

\newcommand{\red}[1]{\textcolor{red}{#1}}
\setlist{itemsep=0pt, topsep=3pt}
\numberwithin{equation}{section}
\makeatletter
\newcommand{\setlabeleq}{\refstepcounter{equation}\textup{\tagform@{\theequation}}}
\makeatother

\makeatletter
\newlength{\reducBoite}
\newcommand{\fbeq}[1][1]{\setlength{\reducBoite}{-2.3ex}\vspace{#1\reducBoite}}
\newcommand{\fbi}[1][1]{\setlength{\reducBoite}{-1ex}\vspace{#1\reducBoite}}

\newcounter{paragraphe}[section]
\newcommand{\numpar}{\thesection.\theparagraphe{}}
\newcommand{\nvnumpar}{%
	\addtocounter{paragraphe}{1}%
	\phantomsection%
	\numpar%
	\xdef\@currentlabel{\numpar}}

\cornersize*{12pt}
\newcommand{\width@boite}{0.972\textwidth}

\newcommand{\avant@boite}{%
	\setlength{\fboxrule}{0.5pt}%
	\setlength{\fboxsep}{5.5pt}}

\newcommand{\dbt@boite}[1]{
	\vspace{1pt}%
	\setlength{\fboxsep}{3pt}%
	\textbf{#1}}
\newcommand{\titreligne}{\vspace{1pt}\\}

\newcommand{\fin@boite}{\vspace{1pt}}

\newcommand{\apres@boite}{%
	\vspace{14pt plus 2pt minus 1pt}}
\newcommand{\apres@par}{\vspace{14pt plus 2pt minus 1pt}}

\NewEnviron{boiteronde}[1]{%
	\avant@boite
	\noindent\ovalbox{\parbox{\width@boite}{\dbt@boite{#1}\BODY\fin@boite}}%
	\apres@boite} 

\mdfdefinestyle{styleboiteronde}{roundcorner=6pt, linewidth=0.5pt, innerleftmargin=5.5pt, innerrightmargin=5.5pt, innertopmargin=5.5pt, innerbottommargin=5.5pt}
\newenvironment{boiterondesecable}[1]{%
	\avant@boite
	\begin{mdframed}[style=styleboiteronde]\dbt@boite{#1}}
	{\fin@boite\end{mdframed}\apres@boite}

\NewEnviron{boitecarree}[1]{%
	\avant@boite
	\noindent\fbox{\parbox{\width@boite}{\dbt@boite{#1}\BODY\fin@boite}}%
	\apres@boite}
	
\mdfdefinestyle{styleboitecarre}{linewidth=0.5pt, innerleftmargin=5.5pt, innerrightmargin=5.5pt, innertopmargin=5.5pt, innerbottommargin=5.5pt}
\newenvironment{boitecarreesecable}[1]{%
	\avant@boite
	\begin{mdframed}[style=styleboitecarre]\dbt@boite{#1}}
	{\end{mdframed}\apres@boite}

\NewEnviron{grosseboitecarree}[1]{%
	\avant@boite
	\setlength{\fboxrule}{0.5pt}%
	\noindent\doublebox{\parbox{0.95\textwidth}{\dbt@boite{#1}\BODY\fin@boite}}%
	\apres@boite}
	
\newenvironment{paragraphesansboite}[1]{\textbf{#1}}{\apres@par}
\newenvironment{textesansboite}{}{\apres@par}

\NewEnviron{definition}[1]{%
	\avant@boite
	\noindent\ovalbox{\parbox{\width@boite}{\dbt@boite{Definition \nvnumpar{}: #1}\BODY\fin@boite}}
	\apres@boite} 

\NewEnviron{proposition}[1]{%
	\avant@boite
	\noindent\fbox{\parbox{\width@boite}{\dbt@boite{Proposition \nvnumpar{}: #1}\BODY\fin@boite}}%
	\apres@boite}

\NewEnviron{lemme}[1]{%
	\avant@boite
	\noindent\fbox{\parbox{\width@boite}{\dbt@boite{Lemma \nvnumpar{}: #1}\BODY\fin@boite}}%
	\apres@boite}

\NewEnviron{corollaire}[1]{%
	\avant@boite
	\noindent\fbox{\parbox{\width@boite}{\dbt@boite{Corollary \nvnumpar{}: #1}\BODY\fin@boite}}%
	\apres@boite}

\NewEnviron{theoreme}[1]{%
	\avant@boite
	\noindent\fbox{\parbox{\width@boite}{\dbt@boite{Theorem \nvnumpar{}: #1}\BODY\fin@boite}}%
	\apres@boite}

\definecolor{grisfonce}{gray}{0.3}
\theoremstyle{empty}
\theorembodyfont{\normalfont}
\theoremsymbol{\ensuremath{\scriptstyle\Box}}
\newtheorem{demonstration}{}
\newenvironment{demo}[1][Proof]{%
	\begingroup\color{grisfonce}%
	\begin{demonstration}\underline{#1:} %
	}{\end{demonstration}\endgroup\vspace{-3pt}\apres@par}

\newenvironment{remarque}{\textbf{Remark:}}{\apres@par}
\newenvironment{remarques}{\textbf{Remarks:}}{\apres@par}
\newenvironment{exemple}{\textbf{Example:}}{\apres@par}

\newlength{\sousdemo@length}
\newenvironment{sousdemo}{\vspace{-1.5ex}%
    \setlength{\sousdemo@length}{\textwidth-\linewidth}
	\def\FrameCommand{\hspace{4mm}\hspace{\sousdemo@length} {\vrule width 0.6pt}\hspace{-\sousdemo@length}{\vrule width 0pt} \fboxsep=\FrameSep}%
	\MakeFramed {\advance \hsize-\width\FrameRestore}}
	{\endMakeFramed \vspace{-4ex}\hspace{4mm} \rule{0.6pt}{1mm}\rule{1.5mm}{0.6pt}\par}
\makeatother

\let\Smallrightarrow\Rightarrow
\let\Smallleftarrow\Leftarrow

\renewcommand{\Rightarrow}{\mathrel{\mbox{$=\hspace{-6.5pt}\Smallrightarrow$}}}
\renewcommand{\Leftarrow}{\mathrel{\mbox{$\Smallleftarrow\hspace{-6.5pt}=$}}}
\renewcommand{\Leftrightarrow}{\mathrel{\mbox{$\Smallleftarrow\hspace{-6pt}\Smallrightarrow$}}}

\newcommand{\couche}{\Lambda}
\newcommand{\bordgauche}{\Sigma_\Omega}
\newcommand{\borddroit}{\Sigma_\Lambda}
\newcommand{\Rcoin}{R_{\mathrm{c}}}
\newcommand{\fSource}{f_{\mathrm{s}}}

\makeatletter
\newcommand{\TZ}{\@ifstar{\frac\pi\Theta \mathbb{Z}^*}{\frac\pi\Theta \mathbb{Z}}}
\makeatother
\newcommand{\TN}{\frac\pi\Theta \mathbb{N}^*}
\newcommand{\NN}{\mathbb{P}}

\newcommand{\Rbasic}[2]{R_{#1}^{\raisebox{1pt}{$\scriptscriptstyle #2$}}}
\newcommand{\ROmegaPi}{R_{\Delta}}
\newcommand{\RGammaPi}{R_{\mathrm{N}}}
\newcommand{\RcouchePi}{R_{\partial_Y^2}}
\newcommand{\ROmegaOmega}{\Rbasic{\Delta}{\Omega}}
\newcommand{\RGammaOmega}{\Rbasic{\mathrm{D}}{\Omega}}
\newcommand{\Aperp}{\mathcal{A}_{\perp}}
\newcommand{\phiOmega}{\phi^{\raisebox{1pt}{$\scriptscriptstyle\Omega$}}}
	
\newcommand{\bfu}{\mathbf{u}}
\newcommand{\bfv}{\mathbf{v}}
\newcommand{\bfr}{\mathbf{r}}
\newcommand{\bfue}{\mathbf{u}^\varepsilon}
\newcommand{\Abarre}{\overline{\mathcal{A}}}
\newcommand{\Amacro}{\mathcal{A}^+\hspace{-0.15ex}}
\newcommand{\Amicro}{\mathcal{A}^-\hspace{-0.15ex}}
\newcommand{\Amaicro}{\mathcal{A}^\pm\hspace{-0.15ex}}
\newcommand{\Aemacro}{\mathcal{A}^+_\varepsilon(\Pi)}
\newcommand{\Aemicro}{\mathcal{A}^-_\varepsilon(\Pi)}
\newcommand{\Aemaicro}{\mathcal{A}^\pm_\varepsilon(\Pi)}
\newcommand{\RIma}[1]{L_{#1}}

\contourlength{0.01em}
\DeclareMathOperator*{\fsum}{\mathchoice
	{\contour*{.}{$\displaystyle\sum$}}%
	{\contour*{.}{$\textstyle\sum$}}%
	{\contour*{.}{$\scriptstyle\sum$}}%
	{\contour*{.}{$\scriptscriptstyle\sum$}}}

\newcommand{\bfuDA}{\bfu^0}
\newcommand{\SDA}{S^\infty}
\newcommand{\uDA}{u^0}
\newcommand{\sumdmacro}{\!\!\!\fsum_{d\in\TZ*\cap(\NN-p)}\!\!\!\!}
\newcommand{\sumdmicro}{\!\!\!\fsum_{d\in\TZ*\cap(p-\NN)}\!\!\!\!}
\newcommand{\Rmacro}{\mathcal{R}^+_\varepsilon}
\newcommand{\Rmicro}{\mathcal{R}^-_\varepsilon}
\newcommand{\Rmaicro}{\mathcal{R}^\pm_\varepsilon}

\newcommand{\Rmacron}{(\Rmacro)^n}
\newcommand{\Rmicron}{(\Rmicro)^n}
\newcommand{\Eps}{\mathcal{H}_\varepsilon}

\newcommand{\progphimacro}{\pi^+_\sigma}
\newcommand{\progphimicro}{\pi^-_\sigma}
\newcommand{\progphimaicro}{\pi^\pm_\sigma}
\newcommand{\progphimiacro}{\pi^\mp_\sigma}
\newcommand{\cmacro}{c^{\bfu\leftarrow S}}
\newcommand{\cmicro}{c^{S\leftarrow\bfu}}
\newcommand{\cfinal}{c^{\bfu\leftarrow\bfu}}
\newcommand{\cfinalbis}{\tilde{c}^{\bfu\leftarrow\bfu}}
\newcommand{\intvl}[1]{[\![#1]\!]_{\NN\!,\mathbb{N}}}

\newcommand{\ptronc}{P}
\newcommand{\Vmacro}{\mathtt{V}^{\mathtt{u}}}
\newcommand{\Vmicro}{\mathtt{V}^{\mathtt{S}}}
\newcommand{\Pmacro}{\mathtt{P}}
\newcommand{\Pmicro}{\mathtt{Q}}
\newcommand{\Pmacroo}[2]{\Pmacro^{#1}_{#2}}
\newcommand{\Pmicroo}[2]{\Pmicro^{#1}_{#2}}
\newcommand{\smacro}{\mathtt{u}}
\newcommand{\smicro}{\mathtt{S}}
\newcommand{\sMacro}{\sigma(\bfu)}
\newcommand{\sMicro}{\sigma(S)}
\newcommand{\Vcoin}{V}
\newcommand{\crown}{C_\eta}
\newcommand{\eps}{\epsilon}

\newcommand{\mindeg}{\operatorname{deg_{min}}}
\newcommand{\maxdeg}{\operatorname{deg_{max}}}
\newcommand{\khiz}{\chi_0}
\newcommand{\khii}{\chi_\infty}
\newcommand{\khif}{\chi_f}

\newcommand{\prodH}[1]{\mathcal{H}^{#1}_\times}
\newcommand{\Res}{\operatorname{Res}}

\makeatletter
\newcommand{\sbd@textstyle}[1]{
	\def\@sbdHmarge{3mu} \def\@sbdVmarge{0.2ex} \def\@sbdVthickness{-0.05ex}
	\mkern\@sbdHmarge{\overbracket[0.2pt][0pt]{
	\raisebox{0pt}[\height+\@sbdVthickness]{$\overbracket[0.2pt][0pt]{\mkern-\@sbdHmarge
	\raisebox{0pt}[\height+\@sbdVmarge]{$\textstyle #1$}
	\mkern-\@sbdHmarge}$}}}\mkern\@sbdHmarge}
\newcommand{\sbd@scriptstyle}[1]{
	\def\@sbdHmarge{2mu} \def\@sbdVmarge{-0.05ex} \def\@sbdVthickness{-0.05ex}
	\mkern\@sbdHmarge{\overbracket[0.2pt][0pt]{
	\raisebox{0pt}[\height+\@sbdVthickness]{$\overbracket[0.2pt][0pt]{\mkern-\@sbdHmarge
	\raisebox{0pt}[\height+\@sbdVmarge]{$\scriptstyle #1$}
	\mkern-\@sbdHmarge}$}}}\mkern\@sbdHmarge}
\newcommand{\sbd}[1]{\mathchoice{\sbd@textstyle{#1}}{\sbd@textstyle{#1}}{\sbd@scriptstyle{#1}}{\sbd@scriptstyle{#1}}}
\makeatother

\begin{document}
\title{Asymptotic analysis at any order of Helmholtz's problem in a corner with a thin layer: an algebraic approach}
\author{Cédric Baudet*}
\date{August 24, 2025}
\maketitle

\begin{paragraphesansboite}{Abstract:}
We consider the Helmholtz equation in an angular sector partially covered by a homogeneous layer of small thickness, denoted $\varepsilon $. We propose in this work an asymptotic expansion of the solution with respect to $\varepsilon $ at any order. This is done using matched asymptotic expansion, which consists here in introducing different asymptotic expansions of the solution in three subdomains: the vicinity of the corner, the layer and the rest of the domain. These expansions are linked through matching conditions. The presence of the corner makes these matching conditions delicate to derive because the fields have singular behaviors. Our approach is to reformulate these matching conditions purely algebraically by writing all asymptotic expansions as formal series. By using algebraic calculus we reduce the matching conditions to scalar relations linking the singular behaviors of the fields. These relations have a convolutive structure and involve some coefficients that can be computed analytically. Our asymptotic expansion is justified rigorously with error estimates.
\end{paragraphesansboite}

\begin{paragraphesansboite}{Keywords:}
asymptotic analysis, Helmholtz's equation, matched asymptotic expansions, corner singularities, algebraic formal series.
\end{paragraphesansboite}

\begin{textesansboite}
* POEMS, CNRS, Inria, ENSTA, Institut Polytechnique de Paris, 91120 Palaiseau, France.\\
Email: cedric.baudet@ensta.fr
\end{textesansboite}

\section*{Introduction}

\begin{textesansboite}
Problems that involve thin layers appear in many areas, from composite materials engineering \cite{FakKha14} to biology \cite{CauHadLiNgu16}, including elasticity \cite{BelGeyKra16, FriJamMul02}, fluid mechanics \cite{RicRuyVil16, AchPirVal98, JagMik01} and electrochemistry \cite{VerMaiNov10}. Applications are especially numerous in electromagnetism, let us mention the studies of thin dielectric layers \cite{HadJia15, Mak08, SenVol95}, ferromagnetic films \cite{AmmHalHam00, HadJol01} and the skin effect \cite{ZutKno05}. All these situations are numerically challenging because they require finely meshing the thin structure, which is very costly when its thickness is very small compared to the wavelength and the size of the objects. In this work we propose to overcome this difficulty by using an asymptotic expansion of the solution, such that each term of the expansion is cheaper to compute than the solution itself.\\

Infinite planar layers and smooth curved layers were studied during the 90s in \cite{BenLem96, EngNed93}. Their method is to stretch the layer in its transverse direction into a standard layer of thickness 1, and look for a Taylor-type asymptotic expansion as a sum of integer powers of the original thickness, denoted $\varepsilon $. The terms of this expansion can be computed by induction. Those results were later extended to heterogeneous and periodic layers in \cite{AmmHe97, AmmHe98, AmmLat98} and more recently in \cite{DelHadJol12, Bou23}.\\

Here we want to handle more realistic situations where the coating has angles or covers only partially the obstacle. We consider a two-dimensional model where the domain is the union of an infinite angular sector and the coating, potentially with a perturbation at the corner of size proportional to $\varepsilon $. This was studied for Poisson's problems in \cite{CalCosDauVia06, AuvVia18, AuvVia19}, providing an asymptotic expansion at any order and approximate models. These works show the presence of non-integer powers of $\varepsilon $ and integer powers of $\ln \varepsilon $ in the asymptotic expansion, that are linked to the corner singularities of the solution. See similar studies for eigenvalue problems in \cite{Naz11, GomLobNazPer06a, GomLobNazPer06b}. That asymptotic expansion at any order was generalized to periodic layer in \cite{DelSchSem16}, still for Poisson's problems. As for the Helmholtz equation, \cite{JolTor06} proved a similar asymptotic expansion at any order in the related case of a half plane with a thin slot. In comparison, Helmholtz's problems involving finite layers not only present the same difficulties, but they also lead to much more complex singularities, which prompted us to introduce new and more efficient algebraic calculus tools in order to obtain an expansion at any order. Let us mention also that \cite{BenMakTor11, Mak08}, resp. \cite{SemDelSch18}, propose asymptotic expansions of Helmholtz's problems up to order 2 in presence of homogeneous, resp. periodic, layers.\\ 

We can identify in these works two methods of analysis: multiscale asymptotic expansions and matched asymptotic expansions (see \cite{VisLyu57, Naz81a, Naz81b} and \cite{Van64, Ili92, Eck94} respectively for a general presentation). They both involve two types of fields: ``far fields'' depending on the macroscopic scale described by $(x,y)$ and ``near fields'' depending on the microscopic scale described by $( \frac{x}{\varepsilon } , \frac{y}{\varepsilon } )$. In multiscale expansions, far and near fields are defined in the whole domain and the near fields tend to 0 towards infinity so that they describe a boundary layer effect in the ``near zone'' (the vicinity of the corner or the layer, depending on the situation). In contrast, matched asymptotic expansions involve near fields only in the near zone and far fields only in the ``far zone'' (the rest of the domain), and the near and far fields have to coincide in an intermediate zone.\\

In this paper, we chose the method of matched expansions. In addition, we propose a new algebraic approach to derive the matching conditions, as they seem too intricate in our problem to be reasonably obtained at any order by classical means, especially if one wants to treat any corner angle. It lies on three main ideas : a general explicit expression of the singularities, operators to generate and manipulate efficiently families of singularities, and formal series to perform rigorous calculations on all orders at once. It avoids specific cumbersome calculations, replacing them with abstract generalizable ones. We believe that this approach gives a better understanding of the structure of the asymptotic expansion at any order. It reveals a convolutive structure and it provides explicit expressions to compute exactly and very cheaply the constants that appear in the obtained matching formulas.\\

We consider the Helmholtz equation with absorption because it brings obvious well-posedness and stability of the problem uniformly in $\varepsilon $, which allows us to focus on asymptotic expansion techniques. The case without absorption requires to design a specific radiation condition, that will be the object of a future paper. We do not put any restriction on the corner angle. Moreover, we apply a Dirichlet condition on the boundary. The extension to Neumann is not obvious and will be presented in a forthcoming article.\\

This paper is organized as follows. In Section~\ref{sec: def du pb}, we define the problem, state the main result and introduce the method based on matched asymptotic expansion. The matching condition around the corner are derived using an algebraic approach in Section~\ref{sec: raccord de coin}. It is the most original part of the article. In Section~\ref{sec: construction des champs}, we introduce appropriate frameworks which allow to define uniquely the terms of the asymptotic expansion. Error estimates are performed in Section~\ref{sec: estimations d'erreur}, proving the main result of the paper.
\end{textesansboite}

\begin{paragraphesansboite}{Acknowledgment:}
I would like to thank Sonia Fliss and Patrick Joly for the helpful discussions we had about the writing of this paper.
\end{paragraphesansboite}

\tableofcontents

\section{Setting of the problem and the method}
\label{sec: def du pb}

\subsection{Definition of the problem and main result}
\label{subsec: definition du probleme etudie}

\begin{figure}[h]
\begin{center}
\includegraphics{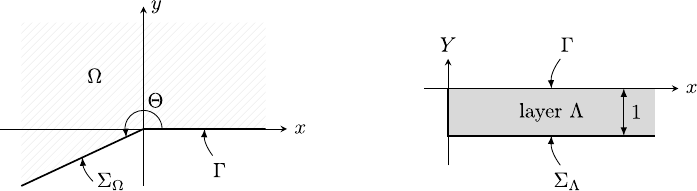}
\caption{The domains $\Omega$ (on the left) and $\couche$ (on the right)}
\label{fig:Omega et couche}
\end{center}
\end{figure}

\begin{figure}[h]
\begin{center}
\hspace*{10mm}\includegraphics{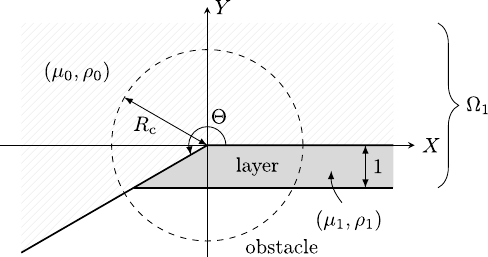}
\caption{The domain $\Omega_1$ with a configuration example of $\Omega_1 \cap B(0,\Rcoin)$}
\label{fig:Omega1}
\end{center}
\end{figure}

\begin{figure}[h]
\begin{center}
\includegraphics{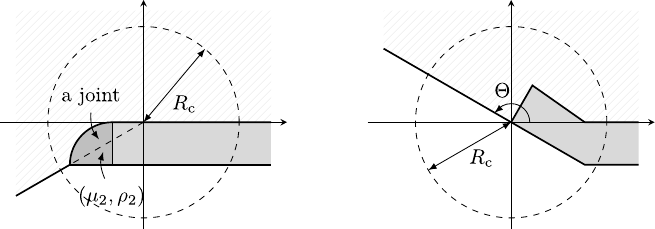}
\caption{Other configuration examples of $\Omega_1 \cap B(0,\Rcoin)$ for different values of $\Theta$}
\label{fig:autres exemples de Omega1}
\end{center}
\end{figure}

\begin{textesansboite}
To describe the domain, let us introduce $\Theta  \in (0,2\pi )$, $\Omega  := \{(r \cos \theta ,r \sin \theta ) \mid r\in \mathbb{R}_+^*,\theta \in (0,\Theta )\}$, $\couche := \mathbb{R}_+^* \times (-1,0)$, $\Gamma  = \mathbb{R}_+^* \times \{0\}$, $\bordgauche=\{(r \cos \Theta ,r \sin \Theta ) \mid r\in \mathbb{R}_+^*\}$ and $\borddroit := \mathbb{R}_+^* \times \{-1\}$. All these sets are shown in Figure~\ref{fig:Omega et couche}. Then let $\Omega _{1} \subset \mathbb{R}^2 $ be an open set that coincides with $\Omega \cup \Gamma \cup \couche$ outside of the disc $B(0,\Rcoin)$ for some $\Rcoin \in \mathbb{R}_+^*$. In addition, let $\mu ,\rho  \in L^{\infty } (\Omega _{1} )$ be two functions greater than a positive constant (ellipticity assumption), and equal to $\mu _{0}$ and $\rho _{0}$ in $\Omega  \setminus B(0,\Rcoin)$ and to $\mu _{1}$ and $\rho _{1}$ in $\couche \setminus B(0,\Rcoin)$. See Figures~\ref{fig:Omega1} and \ref{fig:autres exemples de Omega1} for different configurations.
\end{textesansboite}

\begin{textesansboite}
Let $\varepsilon >0$. The physical domain is given by $\Omega _{\varepsilon } := \{(x,y)\in \mathbb{R}^2  \mid ( \frac{x}{\varepsilon } , \frac{y}{\varepsilon } )\in \Omega _{1} \}$. We introduce the scaled coefficients $\mu _{\varepsilon } :(x,y)\in \Omega _{\varepsilon } \mapsto  \mu ( \frac{x}{\varepsilon } , \frac{y}{\varepsilon } )$ and $\rho _{\varepsilon } :(x,y)\in \Omega _{\varepsilon } \mapsto  \rho ( \frac{x}{\varepsilon } , \frac{y}{\varepsilon } )$, and the scaled variables $X := \frac{x}{\varepsilon }$ and $Y:= \frac{y}{\varepsilon }$. Let $\omega \in \mathbb{C} \setminus \mathbb{R}$ and $\fSource\in H^{-1} (\Omega )$ a source term s.t. $\operatorname{dist}(\operatorname{supp}(\fSource),\Gamma )>0$. We denote $u_{\varepsilon }$ the unique solution in $H^{1} _{0} (\Omega _{\varepsilon } )$ of 
\begin{equation}{\label{eq: pb etudie}}
\operatorname{div}(\mu _{\varepsilon } \nabla  u_{\varepsilon } )+\omega ^2 \rho _{\varepsilon } u_{\varepsilon } =\fSource \qquad \text{ in }\Omega _{\varepsilon } 
\end{equation}
$\Im(\omega ) \neq 0$ is a technical assumption that makes this problem well-posed (it suffices to use the Lax-Milgram theorem) with a stability constant independent of $\varepsilon $ :
\begin{equation}{\label{eq: uniforme coercivite}}
\exists C>0,\;\forall \fSource \in H^{-1} (\Omega ), \;\forall \varepsilon >0, \quad  \|u_{\varepsilon } \|_{H^{1} (\Omega _{\varepsilon } )} \leqslant  C \|\fSource\|_{H^{-1} (\Omega )}
\end{equation}
The case $\Im(\omega )=0$ is an open question and will be the object of a future work.\\

The main result of this paper is given in the following theorem, proven in Section~\ref{sec: estimations d'erreur}, page~\pageref{demo du resultat principal}.
\end{textesansboite}

\begin{theoreme}{asymptotic expansion of $u_{\varepsilon }$\titreligne}
\label{resultat principal}%
Let $\NN := \mathbb{N}+ \frac{\pi }{\Theta } \mathbb{N}$. There exist $(n_{p} )\in \mathbb{N}^{\NN}$ and a family $(u_{p,\ell } )_{p\in \NN,\ell \in [\![0,n_{p} ]\!]}$ of elements of $H^{1} _{\mathrm{loc}} (\Omega )$ that can be build recursively w.r.t. $p$ (see Theorem~\ref{construction directe des u p,l} for the construction) such that
\[\forall P\in \mathbb{R}_+,\forall \delta >0, \qquad  \bigg\|u_{\varepsilon } -\sum _{p\in \NN\cap [0,\ptronc]} \sum _{\ell =0} ^{n_{p}} \varepsilon ^{p} \ln ^{\ell } \!\varepsilon \, u_{p,\ell } \bigg\|_{H^{1} (\Omega  \setminus B(0,\delta ))} = o(\varepsilon ^{\ptronc} ) \qquad  \text{when}  \ \varepsilon \rightarrow 0.\]\fbeq
\end{theoreme}

\begin{textesansboite}
The presence of integer powers of $\varepsilon $ is entirely classical in asymptotic analysis. Integer powers of $\varepsilon ^{\pi /\Theta }$ and $\ln \varepsilon $ can be found in other asymptotic expansions involving corners, see \cite{CalCosDauVia06, AuvVia18, AuvVia19}. Theorem~\ref{resultat principal} can be extended to the case where $\bordgauche$ is covered by another layer (see Remark~\ref{rq: espaces A pour 2 couches} for a useful point).
\end{textesansboite}

\begin{paragraphesansboite}{Notations:}
We denote $(x,y)$ the cartesian coordinates, $(r,\theta )$ the polar coordinates with $\theta  \in [0,2\pi )$, $B(0,r)$ the disc of $\mathbb{R}^2 $ of radius $r$ centered at $(0,0)$, $k_{i} := \omega  \sqrt{\rho _{i} /\mu _{i} }$ for any $i\in \{0,1\}$ and $\alpha  :=e^{-\mathrm{i}\Theta }$.
\end{paragraphesansboite}

\subsection{The matched asymptotic expansion method}
\label{subsec: presentation du dvlpmt raccorde}

\begin{figure}[h]
\begin{center}
\hspace*{12mm}\includegraphics{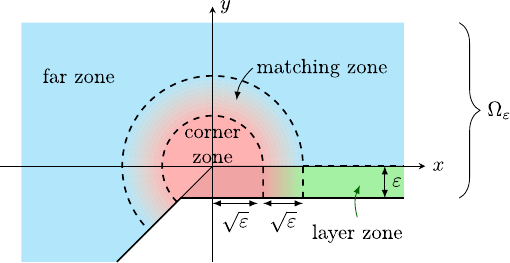}
\caption{Zones of the matched asymptotic expansion}
\label{fig: zones du dvlpmt raccorde}
\end{center}
\end{figure}

\begin{textesansboite}
To take into account the different behaviors of the solution in the layer, near the corner and far from the corner and the layer, we divide $\Omega _{\varepsilon }$ in three zones, illustrated in Figure~\ref{fig: zones du dvlpmt raccorde}. In each zone we postulate an asymptotic expansion in powers of $\varepsilon $ and $\ln \varepsilon $, called ``ansatz''.
\end{textesansboite}

\begin{textesansboite}
Let us assume that for any $(p,\ell )\in \NN\times \mathbb{N}$ there exist three functions independent of $\varepsilon $ – namely $u_{p,\ell }$ defined on $\Omega $ called ``far field'', $U_{p,\ell }$ defined on $\couche$ called ``layer field'' and $S_{p,\ell }$ defined on $\Omega _{1}$ called ``corner field'' – such that $u_{\varepsilon }$ is formally written as:
\begin{itemize}
	\item $u_{\varepsilon } (x,y) =\displaystyle \sum _{p\in \NN} \sum _{\ell \in \mathbb{N}} \varepsilon ^{p} \ln ^{\ell } \!\varepsilon \;u_{p,\ell } (x,y)$ when $r:=\|(x,y)\|\geqslant \sqrt{\varepsilon }$ and $(x,y)\not\in \mathbb{R}_+\times (-\varepsilon ,0)$ (far zone),
	\item $u_{\varepsilon } (x,y) =\displaystyle \sum _{p\in \NN} \sum _{\ell \in \mathbb{N}} \varepsilon ^{p} \ln ^{\ell } \!\varepsilon \;U_{p,\ell } \Big(x, \frac{y}{\varepsilon } \Big)$ when $x\geqslant \sqrt{\varepsilon }$ and $y\in (-\varepsilon ,0)$ (layer zone),
	\item $u_{\varepsilon } (x,y) =\displaystyle \sum _{p\in \NN} \sum _{\ell \in \mathbb{N}} \varepsilon ^{p} \ln ^{\ell } \!\varepsilon \;S_{p,\ell } \Big( \frac{x}{\varepsilon } , \frac{y}{\varepsilon } \Big)$ when $r\leqslant 2\sqrt{\varepsilon }$ (corner zone).
\end{itemize}\fbi
\end{textesansboite}

\begin{remarque}
We will see in Proposition~\ref{borne sur l} that: $\forall p\in \NN,\exists n_{p} \in \mathbb{N},\forall \ell >n_{p} ,\; u_{p,\ell } = 0 \text{ and } U_{p,\ell } = 0 \text{ and } S_{p,\ell } = 0$.
\end{remarque}

\begin{textesansboite}
Injecting the above sums in the Helmholtz equation, using that $\partial _{x} ^2  \big[\varphi ( \frac{x}{\varepsilon } )\big]=\varepsilon ^{-2} [\partial _{X} ^2 \varphi ](X = \frac{x}{\varepsilon } )$ and $\partial _{y} ^2  \big[\varphi ( \frac{y}{\varepsilon } )\big]=\varepsilon ^{-2} [\partial _{Y} ^2 \varphi ](Y = \frac{y}{\varepsilon } )$ for any function $\varphi $, and formally identifying the powers of $\varepsilon $ and $\ln \varepsilon $, one can easily derive the following volume equations and edge conditions for the various fields.

\begin{tabular}{@{}c@{}c@{}}
\begin{minipage}{0.48\linewidth}
\begin{equation}{\label{eq: eq vol champs lointains}}
\left\{\begin{array}{r@{\;=\;}l@{\;}l}
\mu _{0} \Delta  u_{p,\ell } +\omega ^2 \rho _{0} u_{p,\ell } & f\,\delta _{p,0} \,\delta _{\ell ,0} & \text{ in }  \Omega  \\[0.5mm]
u_{p,\ell } & U_{p,\ell } & \text{ on }  \Gamma \\[0.5mm]
u_{p,\ell } & 0 & \text{ on }  \bordgauche
\end{array}\right. \end{equation}
\end{minipage}
&
\begin{minipage}{0.52\linewidth}
\begin{equation}{\label{eq: eq vol champs de couche}}
\left\{\begin{array}{r@{\;=\;}l@{\;}l}
\mu _{1} \partial _{Y} ^2  U_{p,\ell } & -(\mu _{1} \partial _{x} ^2 +\omega ^2 \rho _{1} )U_{p-2,\ell } & \text{ in }  \couche \\[0.5mm]
\mu _{1} \partial _{Y} U_{p,\ell } & \mu _{0} \,\partial _{y} u_{p-1,\ell } & \text{ on }  \Gamma  \\[0.5mm]
U_{p,\ell } & 0 & \text{ on }  \borddroit
\end{array}\right. \end{equation}
\end{minipage}
\end{tabular}
\vspace{1mm}

\begin{equation}{\label{eq: eq vol champs de coin}}
\left\{\begin{array}{r@{\;=\;}ll}
\operatorname{div}(\mu  \nabla  S_{p,\ell } ) & -\omega ^2 \rho  \,S_{p-2,\ell } & \text{ in }  \Omega _{1} \\[0.5mm]
S_{p,\ell } & 0 & \text{ on }  \partial  \Omega _{1}
\end{array}\right. \vspace{1.5mm} \end{equation}

where we denote by convention $u_{p,\ell } =0$, $U_{p,\ell } =0$ and $S_{p,\ell } =0$ for any $(p,\ell )\in \mathbb{R} \setminus \NN \times \mathbb{N}$, and $\delta _{i,j} :=1$ if $i=j$ and 0 if not.
\end{textesansboite}

\begin{remarques}
\begin{itemize}
	\item The condition $\mu _{1} \partial _{Y} U_{p,\ell } = \mu _{0} \partial _{y} u_{p-1,\ell }$ is included in the problem satisfied by $U_{p,\ell }$ whereas $u_{p,\ell } =U_{p,\ell }$ is included in the problem satisfied by $u_{p,\ell }$ so that the construction is inductive: $u_{p-1,\ell }$ allows to build $U_{p,\ell }$, which allows to build $u_{p,\ell }$.
	\item The problem satisfied by $U_{p,\ell }$ depends only on $Y$, the variable $x$ playing the role of a parameter.\\
\end{itemize}\fbi
\end{remarques}

\begin{textesansboite}
\eqref{eq: eq vol champs lointains}--\eqref{eq: eq vol champs de coin} would be sufficient to uniquely define the fields, if they were in their natural variational spaces (e.g. $H^{1} (\Omega )$ for $u_{p,\ell }$). But we need to take into account a matching condition: the far and corner fields must coincide in the intersection of the far and corner zones, and similarly for the layer and corner fields. These intersections form the matching zone (see Figure~\ref{fig: zones du dvlpmt raccorde}). Given that $\varepsilon \rightarrow 0$ and $\frac{\sqrt{\varepsilon } }{\varepsilon } \rightarrow \infty $, this zone tends to $(0,0)$ w.r.t. the far and layer fields, but it tends to infinity w.r.t. the corner fields. Thus, the matching condition links the asymptotic behavior of far and layer fields at the corner to the one of corner fields at infinity:
\begin{equation}{\label{eq: hypothese de raccord}}
\left\{\begin{array}{cl}
	\displaystyle \sum  \varepsilon ^{p} \ln ^{\ell } \!\varepsilon \;u_{p,\ell } (x,y) \approx  \sum  \varepsilon ^{p} \ln ^{\ell } \!\varepsilon \;S_{p,\ell } ({\textstyle \frac{x}{\varepsilon } , \frac{y}{\varepsilon } }) &\text{ in }\Omega  \text{ when }  r\rightarrow 0 \text{ and } \frac{r}{\varepsilon } \rightarrow \infty \\[1mm]
	\displaystyle \sum  \varepsilon ^{p} \ln ^{\ell } \!\varepsilon \;U_{p,\ell } (x,{\textstyle \frac{y}{\varepsilon } }) \approx  \sum  \varepsilon ^{p} \ln ^{\ell } \!\varepsilon \;S_{p,\ell } ({\textstyle \frac{x}{\varepsilon } , \frac{y}{\varepsilon } }) &\text{ in }\couche \text{ when }  x\rightarrow 0 \text{ and } \frac{x}{\varepsilon } \rightarrow \infty 
\end{array}\right.
\end{equation}
We will see that the far fields $u_{p,\ell }$ have an asymptotic expansion at the corner which is roughly a sum of powers of $r$, some of which are positive (like in a Taylor expansion e.g.). The matching conditions imply that these positive powers of $r$ have to appear in the asymptotic expansions of the corner fields at infinity. We call them singularities for the corner fields. Conversely, the asymptotic expansion of the corner fields at infinity contain negative powers of $r$ corresponding to the decay of the variational part, and these powers must be found in the far fields, which corresponds to singularities at the corner. Thus the fields cannot be searched in their natural variational spaces. To derive problems that define uniquely the fields, we must specify their singular part (see Theorems~\ref{cadre fonctionnel des champs lointains} and \ref{cadre fonctionnel pour les champs de coin}). These parts are fixed by the matching conditions as explained in the next section.
\end{textesansboite}

\section{Matching conditions}
\label{sec: raccord de coin}

\begin{figure}[h]
\begin{center}
\includegraphics{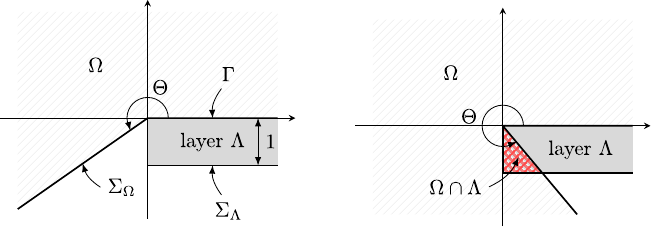}
\caption{The domain $\Pi $ is equal to $\Omega \sqcup \Gamma \sqcup \couche$ when $\Theta \leqslant  \frac{3\pi }{2}$ (open subset of $\mathbb{R}^2 $) and it is a flat Riemannian manifold when $\Theta > \frac{3\pi }{2}$ ($\Omega \cap \couche \neq  \varnothing $).}
\label{fig:Pi}
\end{center}
\end{figure}

\begin{textesansboite}
This section establishes the matching condition linking corner fields to far and layer fields. This is by far the most difficult relation to derive, while all the others have been easily stated in \eqref{eq: eq vol champs lointains}--\eqref{eq: eq vol champs de coin}. In this section we assume that the various fields exist and that they satisfy \eqref{eq: eq vol champs lointains}--\eqref{eq: eq vol champs de coin} and we give a necessary and sufficient condition for the matching assumptions \eqref{eq: hypothese de raccord} to be satisfied. Our approach is based on an algebraic formulation of the problem, that reveals the structure of the matching relations by a rigorous algebraic calculus.\\

\phantomsection\label{page: definition de Pi}
To perform the matching of the corner fields with the far and the layer fields at the same time, we merge the latter two into a single field denoted $\bfu_{p,\ell }$ and called ``far-and-layer field''. It is defined on a new domain $\Pi $, defined as follows:
\begin{itemize}
	\item If $\Theta \leqslant  \frac{3\pi }{2}$, then $\Pi  := \Omega \sqcup \Gamma \sqcup \couche$ (disjoint union) and it is an open of $\mathbb{R}^2 $.
	\item If $\Theta > \frac{3\pi }{2}$, then $\Omega $ and $\couche$ intersect as subsets of $\mathbb{R}^2 $, so the previous definition is not valid anymore (see Figure~\ref{fig:Pi}). Thus, we define $\Pi $ as the disjoint gluing of $\Omega $ and $\couche$ on $\Gamma $ (which is a flat Riemannian manifold).
\end{itemize}\fbi
\end{textesansboite}

\begin{textesansboite}
For all $(p,\ell )\in \NN\times \mathbb{N}$, we define \;$\bfu_{p,\ell } := \left\{\begin{array}{@{\;}l@{\,}l}
	u_{p,\ell } &\text{ in }\Omega  \\[0.5mm]
	U_{p,\ell } &\text{ in }\couche
\end{array}\right.$ and the generalized radial variable $\bfr := \left\{\begin{array}{@{\;}l@{\,}l}
	r &\text{ in }\Omega  \\[0.5mm]
	x &\text{ in }\couche
\end{array}\right.$\\

A straightforward reformulation of \eqref{eq: eq vol champs lointains} and \eqref{eq: eq vol champs de couche} gives that for any $(p,\ell )\in \NN\times \mathbb{N}$: \begin{equation}{\label{eq: eq vol de bfu p,l}}
\left\{\begin{array}{r@{\;=\;}ll}
	\mu _{0} \Delta  \bfu_{p,\ell } +\omega ^2 \rho _{0} \bfu_{p,\ell } & f\,\delta _{p,0} \,\delta _{\ell ,0} &\text{ in }\Omega  \\[0.5mm]
	\mu _{1} \partial _{Y} ^2  \bfu_{p,\ell } & -(\mu _{1} \partial _{x} ^2 +\omega ^2 \rho _{1} )\bfu_{p-2,\ell } &\text{ in }\couche \\[0.5mm]
	\mu _{1} \partial _{Y} \bfu_{p,\ell |Y=0^{-}} & \mu _{0} \,\partial _{y} \bfu_{p-1,\ell |y=0^{+}}  &\text{ on }\Gamma  \\[0.5mm]
	\bfu_{p,\ell |y=0^{+}} - \bfu_{p,\ell |Y=0^{-}} & 0 &\text{ on }\Gamma \\[0.5mm]
	\bfu_{p,\ell } & 0 &\text{ on }\bordgauche \cup \borddroit
\end{array}\right.\end{equation}
where by convention $\bfu_{p,\ell } =0$ when $(p,\ell )\in \mathbb{R} \setminus \NN\times \mathbb{N}$.\\

Let us give some starting point ideas to dive into this section. The matching assumption under study links the asymptotic behaviors of $\sum _{p,\ell } \varepsilon ^{p} \ln ^{\ell } \!\varepsilon \, \bfu_{p,\ell }$ when $\bfr\rightarrow 0$ and $\sum _{p,\ell } \varepsilon ^{p} \ln ^{\ell } \!\varepsilon \, S_{p,\ell }$ when $\bfr\rightarrow \infty $. So we can begin with a look at the asymptotic of $\bfu_{p,\ell }$ when $\bfr\rightarrow 0$, especially on $\Omega $ because it is the most interesting part. First, by \eqref{eq: eq vol de bfu p,l}, $u_{0,0}$ satisfies $\mu _{0} \Delta u_{0,0} +\omega ^2  \rho _{0} u_{0,0} =0$ in the vicinity of the corner in $\Omega $, with homogeneous Dirichlet condition on $\partial \Omega $. So using separation of variables, it is easy to show that:
\begin{equation}{\label{eq: decomposition de u 0,0 en fcts de Bessel}}
u_{0,0} (r,\theta ) \underset{r\ll 1}{=} \sum _{d\in \TN} \sigma _{d} (\bfu_{0,0} ) \,J_{d} (k_{0} r) \,\sin (d\theta ) = \sum _{d\in \TN} \sigma _{d} (\bfu_{0,0} ) \,\Bigg(\sum_{n=0}^\infty  a_{d,n} r^{d+2n} \Bigg) \,\sin (d\theta )\end{equation}
with $\sigma _{d} (\bfu_{0,0} )$ and $a_{d,n}$ some coefficients in $\mathbb{C}$, $J_{d}$ Bessel functions of the first kind and $k_{0} := \omega  \sqrt{\rho _{0} /\mu _{0} }$. Then, one can show that \eqref{eq: eq vol de bfu p,l} implies that $u_{1,0}$ satisfies the Helmholtz equation with condition $\mu _{1} u_{1,0} = \mu _{0} \,\partial _{y} u_{0,0|y=0^{+}}$ on $\Gamma $. Using \eqref{eq: decomposition de u 0,0 en fcts de Bessel}, one can show that there exist some functions $f_{q,i}$ and coefficients $\sigma _{d} (\bfu_{1,0} )$ s.t.:
\[u_{1,0} \underset{r\ll 1}{=} \sum _{d\in \NN} r^{d-1} (f_{d,0} (\theta )+\ln (r)\,f_{d,1} (\theta ))+ \sum _{d\in \TN} \sigma _{d} (\bfu_{1,0} ) \,J_{d} (k_{0} r) \,\sin (d\theta )\]
where the first sum is a particular solution of Helmholtz's equation that has trace $\frac{\mu _{0} }{\mu _{1}} \,\partial _{y} u_{0|y=0^{+}}$ on $\Gamma $, and the second one is a homogeneous solution. More generally, the behavior of $\bfu_{p,\ell }$ when $\bfr\rightarrow 0$ (resp. $S_{p,\ell }$ when $\bfr\rightarrow \infty $) has the following form:
\[\sum_{n=0}^\infty  r^{d_{n}} \ln ^{k_{n}} \!r \cdot  g_{n,1} (\theta ) \quad \text{ in }\Omega  \qquad  \text{ and } \qquad  \sum_{n=0}^\infty  r^{d_{n}} \ln ^{k_{n}} \!r \cdot  g_{n,2} (Y) \quad \text{ in }\couche\]
where $d_{n}$ is an increasing (resp. decreasing) sequence in $\mathbb{Z}+\TZ$ and $k_{n} \in \mathbb{N}$. In order to know whether the behaviors of the far-and-layer fields and the corner fields match, we need to describe the functions $g_{n,i}$. In the literature, they are usually built inductively, by solving 1D problems \cite{CalCosDauVia06, DelSchSem16, SemDelSch18}. However the several inductive source terms of \eqref{eq: eq vol de bfu p,l} imply that the $g_{n,i}$ do not depend on each other only through linear sequences, but rather through a tree of dependencies that grows exponentially as one builds more terms of the asymptotic expansion (see Remark~\ref{rq: utilite des operateurs R}).\\

To overcome the difficulty of matching such complex singularities, 
the first original idea of this section is to decompose them on a basis of functions. Thus, only the coefficients of the decompositions need to be matched. This somehow mimics the simplicity of infinite periodic layers, for which a basis consists of the integer powers of the distance to the layer, since the fields behaviors in the matching zone are polynomial. We provide explicit expressions of those basis functions. Therefore all singularities have analytical expressions, which allow to compute them both exactly and very quickly. This last point is an original contribution compared to the existing literature.\\

The second idea of this section is to define operators in order to describe how singularities are linked to each other. These operators are defined from the solution of Poisson-like problems in $\Pi $, given by \eqref{eq: pb resolu dans A}. They can be composed, added and factorized in all possible ways, which, in contrast to linear sequences of functions, offers sufficient flexibility to describe all singularities. An important result (Theorem~\ref{champs totaux en fct de sigma}) is that all singularities are generated, through these operators, by the fundamental singularities $\phi_{d} :=r^{d} \sin (d \theta )$, $d\in \TZ*$ (which are part of the basis). It follows that it is sufficient to match only the coefficients of the fields behaviors on the $\phi_{d}$ (denoted $\sigma _{d} (\bfu_{p,\ell } )$ for the far fields and $\sigma _{d} (S_{p,\ell } )$ for the near fields). More precisely only the coefficients involved in the non variational behaviors of the fields need to be imposed, i.e. the $\sigma _{d} (\bfu_{p,\ell } )$ with $d<0$ and the $\sigma _{d} (S_{p,\ell } )$ with $d>0$.\\

The last main idea is to use formal series to treat all orders at once and to focus on the calculative nature of the matching, without having to deal with sum truncations and asymptotic remainders like $o(r^{d} )$. Formal series are usually used in asymptotic analysis in an intuitive way like we did in Section~\ref{sec: def du pb}. Here we provided rigorous foundations to formal series in order to use them as a proof tool and to justify non trivial operations on them, such as applying operators, and even formal series of operators. Let us mention that \cite{Fao02} also used an approach with formal series to study a problem with a thin layer, but with more elementary algebraic tools.\\

All these ideas form a system of powerful calculus tools that are suitable to face the sheer complexity of the matching.\\

In Section~\ref{subsec: preliminaires algebriques} we introduce the algebraic formal series. In Section~\ref{subsec: def de A} we define the spaces $\mathcal{A}({\dots })$ of explicit functions that contain the singularities. In Section~\ref{subsec: operateurs de resolution dans A} we build the mentioned operators on singularities. In Section~\ref{subsec: DA des champs lointains et de coin}, we write the asymptotic expansions w.r.t. $\bfr$ of the various fields. Finally, in Section~\ref{subsec: raccordement algebrique}, we re-express the matching conditions \eqref{eq: hypothese de raccord} with equations that can be used to build the fields (Theorem~\ref{condition de raccord}).
\end{textesansboite}

\subsection{Algebraic preliminaries}
\label{subsec: preliminaires algebriques}

\newcommand{\ev}{E}
\newcommand{\evbis}{F}
\begin{textesansboite}
To handle infinite series that may not converge, e.g. ``$\sum _{p,\ell } \varepsilon ^{p} \ln ^{\ell } \!\varepsilon \, u_{p,\ell }$'', we use the algebraic notion of \emph{formal series} introduced in this section. Let $\ev$ be a vector space and $(\ev_{i} )_{i\in I}$ be a family of vector subspaces of $\ev$. To begin, let us remind that $\sum _{i\in I} \ev_{i}$ designates the vector subspace of $\ev$ made of \emph{finite} sums of elements of the $\ev_{i}$. If this sum is direct, we denote it $\bigoplus_{i\in I} \ev_{i}$. From now on, we assume that the sum is direct. In order to deal with \emph{infinite} sums we introduce the following definition.
\end{textesansboite}

\begin{boiteronde}{Notation \nvnumpar{}: }
\label{def: serie formelle}%
Let us denote 
\[\forall (\varphi _{i} )\in \prod _{i\in I} \ev_{i,} \quad  \fsum_{i\in I} \varphi _{i} :=(\varphi _{i} ) \qquad  \text{ and } \qquad  \fsum_{i\in I} E_{i} := \prod _{i\in I} E_{i.}\]
Note the boldness of the symbol $\fsum$. $\fsum_{i\in I} \varphi _{i}$ is not a real sum that can be computed, but just a notation called ``formal series''. Its support is defined as $\{i\in I \mid \varphi _{i} \neq 0\}$. In additional for any $J\subset I$ that contains this support, we also denote $\fsum_{j\in J} \varphi _{j} :=(\varphi _{i} )_{i\in I}$.
\end{boiteronde}

\begin{figure}[h]
\begin{center}
\def\widthFig{4}
\def\heightFig{0.9}
\def\degfFig{1}
\def\icourbe{(exp(\i/4)-1)*4}
\def\diFig{(\icourbe/3.2-1)*\widthFig}
\begin{tikzpicture}
	\draw [>=stealth,->] (-\widthFig,\heightFig) -- (\widthFig,\heightFig) node[above]{$\;\ d$} node[right]{$\;\ \ \Big\}\bigoplus\limits_{d\in\mathbb{R}}^{\rule{0ex}{1ex}} \ev_d$};
	\draw [>=stealth,->] (-\widthFig,0) -- (\widthFig,0) node[above]{$\;\ d$} node[right]{$\;\ \ \Big\}\bigoplus\limits_{d\in\mathbb{R}}^{\rule{0ex}{1ex}} \evbis_d$};
	\draw  (-\widthFig-1,0) node{\ }; 
	\foreach \i in {1,2,3} {
		\draw ({\diFig},\heightFig) node{\footnotesize$\bullet$} node[above]{$\ev_{d_\i}$};
		\draw ({\diFig},0) node{\footnotesize$\bullet$} node[below]{$\evbis_{d_\i}$};
		\draw ({\diFig+\degfFig},0) node{\footnotesize$\bullet$} node[below]{$\evbis_{d_\i +d_0}$};
		\draw[arrows={->[length=6pt,width=6pt]}, dashed] ({\diFig},\heightFig) -- node[midway,right]{\rule[-6pt]{0pt}{1ex}$f$} ({\diFig+\degfFig},0);}
\end{tikzpicture}\par
\caption{Schematic illustration of an operator that has degree $d_{0}$ (here $d_{0} >0$)}
\label{fig: degre d'un operateur}
\end{center}
\end{figure}
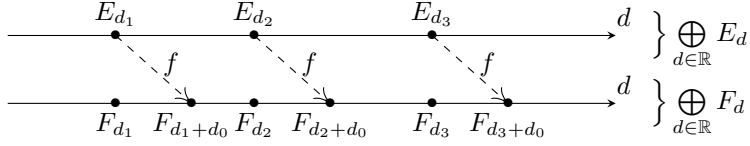

\begin{textesansboite}
There is a canonical injection $\bigoplus \ev_{i} \rightarrow  \fsum \ev_{i}$, that maps any sum $\sum _{i\in I} \varphi _{i}$ with finite support (and $\forall i,\, \varphi _{i} \in \ev_{i}$) to the formal series $\fsum_{i\in I} \varphi _{i}$. So we can consider in practice that $\bigoplus \ev_{i}$ is included in $\fsum \ev_{i}$.\\

We will use Notation~\ref{def: serie formelle} with $I=\mathbb{R}$ and $\ev_{d} =\mathcal{A}_{d} ({\dots })$ a space of functions that behave like $\bfr^{d}$ defined in Section~\ref{subsec: def de A}. In Section~\ref{subsec: operateurs de resolution dans A} we build some operators in the spaces $\mathcal{A}$ that have a translation action on the index $d$. We say that they have a ``degree'' (cf. Definition~\ref{def: degre d'un operateur} and Figure~\ref{fig: degre d'un operateur}). That allows us to naturally extend them to the formal series of the spaces $\mathcal{A}$ via the construction below.
\end{textesansboite}

\begin{definition}{operators with a degree\titreligne}
\label{def: degre d'un operateur}%
Let $\evbis$ be another vector space and $(\evbis_{d} )_{d\in \mathbb{R}}$ be a family of subspaces of $\evbis$ s.t. the sum $\sum \evbis_{d}$ is direct. Let $f:\bigoplus \ev_{d} \rightarrow  \bigoplus \evbis_{d}$ be a linear map and let $d_{0} \in \mathbb{R}$.\\
We say that $f$ has degree $d_{0}$ iff: $\forall d\in \mathbb{R},\, \forall \varphi \in \ev_{d} ,\; f(\varphi )\in \evbis_{d+d_{0}}$. In this case we denote $\deg f:=d_{0}$ and we extend $f$ from $\fsum \ev_{d}$ to $\fsum \evbis_{d}$ by setting $f\big(\fsum_{d\in \mathbb{R}} \varphi _{d} \big):= \fsum_{d\in \mathbb{R}} f(\varphi _{d} )$ for any $(\varphi _{d} )$.
\end{definition}

\begin{figure}[h]
\begin{center}
\def\widthFig{4.3}
\def\dminFig{-\widthFig+1}
\def\pasVFig{0.55}
\def\pasHFig{0.25}
\def\NFig{3}
\def\NmoinsUnFig{2}
\def\degfFig{2}
\def\ij{(\i*\degfFig+\j)}
\def\Nj{(\NFig*\degfFig+\j)}
\def\basFig{-(\NFig+1.55)*\pasVFig}
\begin{tikzpicture}
	\foreach \i in {0,...,\NFig} {
		\draw[>=stealth,->] (-\widthFig,-\i*\pasVFig) -- (\widthFig,-\i*\pasVFig) node[above]{$\;\ d$};}
	\draw[>=stealth,->] (-\widthFig,{\basFig}) -- (\widthFig,{\basFig}) node[above]{$\;\ d$};
	\foreach \i in {1,...,\NFig} {
		\draw (-\widthFig-0.2,{-(\i-0.5)*\pasVFig}) node[left]{$+$};}
	\draw (-\widthFig-0.2,{-(\NFig+0.5)*\pasVFig}) node[left]{\scalebox{0.75}{$\vdots$}\;};
	\draw (-\widthFig-0.2,{\basFig}) node[left]{$=$};
	\draw  (-\widthFig-1.2,0) node{\ };
	\draw (\widthFig,0) node[right]{$\;\ \ \big\}\; \fsum\varphi_d $};
	\draw (\widthFig,-\pasVFig) node[right]{$\;\ \ \big\}\; f\big(\fsum\varphi_d \big)$};
	\foreach \i in {2,...,\NFig} {
		\draw  (\widthFig,-\i*\pasVFig) node[right]{$\;\ \ \big\}\; f^\i\big(\fsum\varphi_d \big)$};}
	\draw (\widthFig,{(-\NFig*\pasVFig+\basFig)/2+0.05}) node[right]{$\qquad \quad$\scalebox{0.75}{$\vdots$}};
	\draw (\widthFig,{\basFig}) node[right]{$\;\ \ \big\}\; \sum f^n \big(\fsum\varphi_d \big)$};
	\begin{scope}
		\clip (-\widthFig,0.1) rectangle (\widthFig-0.15,{\basFig-0.1});
		\foreach \j in {0,4,9,12,18,20,27} {
			\foreach \i in {0,...,\NFig} {
				\draw[dotted, gray] ({\dminFig+\ij*\pasHFig},{\basFig}) -- ({\dminFig+\ij*\pasHFig},{-\i*\pasVFig});
				\draw ({\dminFig+\ij*\pasHFig},{-\i*\pasVFig}) node{\footnotesize$\bullet$};}
			\foreach \i in {0,...,40} {
				\draw ({\dminFig+\ij*\pasHFig},{\basFig}) node{\footnotesize$\bullet$};}
			\foreach \i in {0,...,\NmoinsUnFig} {
				\draw[arrows={->[length=5pt,width=4pt]}, densely dashed] ({\dminFig+\ij*\pasHFig},{-\i*\pasVFig}) -- ({\dminFig+(\ij+\degfFig)*\pasHFig},{-(\i+1)*\pasVFig});}
			\draw[arrows={->[length=5pt,width=4pt]}, densely dashed] ({\dminFig+\Nj*\pasHFig},{-\NFig*\pasVFig}) -- ({\dminFig+(\Nj+\degfFig*0.6)*\pasHFig},{-(\NFig+0.6)*\pasVFig});}
	\end{scope}
	\draw[<-,>=latex] ({\dminFig},0.05) -- ++(0,0.25) node[above]{$d_{\mathrm{inf}}$};
	\foreach \i in {0} {
		\foreach \j in {4} {
			\draw[<-,>=latex] ({\dminFig+\ij*\pasHFig+0.03},{-\i*\pasVFig+0.05}) to[bend left=25] ++(0.5,0.4) node[right]{$\varphi_d \neq0$};}
		\foreach \j in {15} {
			\draw[<-,>=latex] ({\dminFig+\ij*\pasHFig+0.03},{-\i*\pasVFig+0.05}) to[bend left=25] ++(0.5,0.4) node[right]{$\varphi_d =0$};}}
	\draw[<-,>=latex] ({\dminFig+4*\pasHFig+0.03},{\basFig-0.05}) to[bend right=25] ++(0.5,-0.35) node[right]{finite sum};
\end{tikzpicture}\par
\caption{Schematic illustration of $\sum_{n=0}^\infty  f^{n}$ for $f$ a linear map that has a positive degree}
\label{fig: serie geometrique d'operateur}
\end{center}
\end{figure}
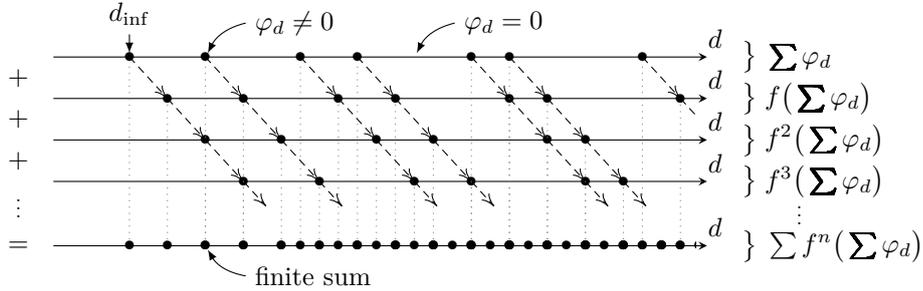

\begin{definition}{}
\label{def: serie geometrique d'operateurs 1D}%
Let $\fsum_{d\in \mathbb{R}} \varphi _{d}$ be an element of $\fsum E_{d}$ with a support bounded from below and $d_{\mathrm{\inf }} := \inf \operatorname{supp}(\varphi _{d} )$. For any linear map $f:\bigoplus E_{d} \rightarrow \bigoplus E_{d}$ that has a positive degree, we define
\begin{equation}{\label{eq: serie geometrique d'un operateur a degre 1D}}
\sum_{n=0}^\infty  f^{n} \bigg(\fsum_{d\in \mathbb{R}} \varphi _{d} \bigg) := \fsum_{d\in \mathbb{R}} \sum _{n\in \mathbb{N}} f^{n} (\varphi _{d-n\deg f} ) \end{equation}
(where $f^{n}$ is the $n$-th iterated composition of $f$). See Figure~\ref{fig: serie geometrique d'operateur}. More generally, for any finite set $\mathcal{F}$ of linear maps $\bigoplus E_{d} \rightarrow \bigoplus E_{d}$ that have positive degrees, we set
\begin{equation}{\label{eq: serie geometrique de plusieurs operateurs a degre 1D}}
\sum_{n=0}^\infty  \sum _{(f_{1} ,\dots , f_{n} )\in \mathcal{F}^{n} } f_{1} \circ \cdots \circ f_{n} \bigg(\fsum_{d\in \mathbb{R}} \varphi _{d} \bigg)
:= \fsum_{d\in \mathbb{R}} \sum _{n\in \mathbb{N}} \sum _{(f_{1} ,\dots , f_{n} )\in \mathcal{F}^{n} } f_{1} \circ \cdots \circ f_{n} (\varphi _{d-\Sigma _{1\leqslant i\leqslant n} \deg  f_{i} }). \end{equation}
We also denote it as $\langle \mathcal{F}\rangle \big(\fsum \varphi _{d} \big)$, or $\langle  \tilde{f}_{1} ,\dots , \tilde{f}_{k} \rangle \big(\fsum \varphi _{d} \big)$ if $\mathcal{F}=\{ \tilde{f}_{1} ,\dots , \tilde{f}_{k} \}$. This is well-defined because, for any $d$, the sums over $n$ in the right-hand sides of \eqref{eq: serie geometrique d'un operateur a degre 1D}--\eqref{eq: serie geometrique de plusieurs operateurs a degre 1D} have a finite number of non-zero terms and belong to $E_{d}$.\\
If $\operatorname{supp}(\varphi _{d} )$ is bounded from above and the elements of $\mathcal{F}$ have a negative degree, we can do the same definition.
\end{definition}

\newcommand{\evtilde}{\tilde{\ev}}
\newcommand{\vecAEps}{v}
\newcommand{\evEps}{v}
\newcommand{\setP}{P}
\begin{textesansboite}
In Section~\ref{subsec: raccordement algebrique}, we introduce formal series with powers of $\varepsilon $ in order to express the matching condition. They are defined similarly to Notation~\ref{def: serie formelle}: for any set $\setP \subset \mathbb{R}$ and any family $(\evtilde_{p} )_{p\in \setP}$ of vector spaces, we denote
\begin{equation}{\label{eq: serie formelle en puissance de eps}}
\forall (\varphi _{p} ) \in \prod _{p\in \setP} \evtilde_{p} ,\quad  \fsum_{p\in \setP} \varepsilon ^{p} \varphi _{p} := (\varphi _{p} ) \qquad  \text{ and } \qquad  \fsum_{p\in \setP} \varepsilon ^{p} \evtilde_{p} :=\prod _{p\in \setP} \evtilde_{p} .\end{equation}
Again this is only a notation and here $\varepsilon $ is not a real number but an algebraic indeterminate. This is similar to Notation~\ref{def: serie formelle} and we will later choose to use either notation depending on the physical meaning of the formal series.\\

For instance, the case $\evtilde=\mathbb{C}$ and $\setP=\mathbb{N}$ gives the classical set of formal power series, usually denoted $\mathbb{C}[[\varepsilon ]]$ (see \cite{Niv69, Sam23}). The Taylor approximations at 0 of any smooth function $f$ can be represented by $\sum _{p=0} ^{\infty } \varepsilon ^{p} \frac{f^{(p)} (0)}{p!} \in \mathbb{C}[[\varepsilon ]]$. Truncations of this series give approximations at a given order. We will use similar representations for the asymptotic expansion of $u_{\varepsilon }$.\\

Let us take $\setP:=\mathbb{R}$, $\evtilde := \fsum_{d\in \mathbb{R}} E_{d}$ and $\evtilde_{p} :=\evtilde$ for any $p$. Let $p_{0} \in \mathbb{R}$ and $f:\bigoplus E_{d} \rightarrow \bigoplus E_{d}$ a linear map that has a degree $d_{0}$. For any $(\varphi _{p} )\in \evtilde^{\mathbb{R}}$, we define:
\begin{equation}{\label{eq: operateur a degre selon p et d}}
(\varepsilon ^{p_{0}} f)\bigg(\fsum_{p\in \mathbb{R}} \varepsilon ^{p} \varphi _{p} \bigg) := \fsum_{p\in \mathbb{R}} \varepsilon ^{p} f(\varphi _{p-p_{0}} ) \end{equation}
These kind of linear maps $\fsum_{p} \varepsilon ^{p} \evtilde \rightarrow  \fsum_{p} \varepsilon ^{p} \evtilde$ are the one said to have a degree. We denote $\deg (\varepsilon ^{p_{0}} f):=(p_{0} ,d_{0} )$.
\end{textesansboite}

\begin{definition}{}
\label{def: serie geometrique d'operateurs 2D}%
Let $\fsum_{(p,d)\in \mathbb{R}^2 } \varepsilon ^{p} \varphi _{p,d} \in \fsum_{p\in \mathbb{R}} \varepsilon ^{p} \evtilde$. Let $\mathcal{G}$ be a finite set of linear maps $\fsum \varepsilon ^{p} \evtilde \rightarrow \fsum \varepsilon ^{p} \evtilde$ that have degrees. There is a finite set $\mathcal{F}$ of linear maps $\bigoplus E_{d} \rightarrow \bigoplus E_{d}$ and $(p_{f} )_{f} \in \mathbb{R}^{\mathcal{F}}$ s.t. $\mathcal{G} = \{\varepsilon ^{p_{f}} f \mid f\in \mathcal{F}\}$. We assume that there is $\vecAEps \in \mathbb{R}^2 $ s.t. $\{\langle (p,d),v\rangle  \mid (p,d)\in \mathbb{R}^2 , \varphi _{p,d} \neq 0\}$ is bounded from below and $\forall g\in \mathcal{G},\; \langle \deg  g,v\rangle >0$. We denote:
\[\sum_{n=0}^\infty  \sum _{(g_{1} ,\dots , g_{n} )\in \mathcal{G}^{n} } \! g_{1} \circ \cdots \circ g_{n} \bigg(\fsum_{(p,d)\in \mathbb{R}^2 } \! \varepsilon ^{p} \varphi _{p,d} \bigg) := \!\!\fsum_{(p,d)\in \mathbb{R}^2 } \!\! \varepsilon ^{p} \sum _{n\in \mathbb{N}} \sum _{(f_{1} ,\dots , f_{n} )\in \mathcal{F}^{n} } \! f_{1} \circ \cdots \circ f_{n} (\varphi _{p-\Sigma _{i} p_{f_{i}} ,d-\Sigma _{i} \deg  f_{i} }).\]
which is well-defined in $\fsum \varepsilon ^{p} \evtilde$. We also denote it as $\langle \mathcal{G}\rangle \big(\fsum_{p,d} \varepsilon ^{p} \varphi _{p,d} \big)$.
\end{definition}

\begin{definition}{}
Let ``$\ln \varepsilon $'' be here an algebraic indeterminate independent from the indeterminate $\varepsilon $. We denote $\ev[\ln \varepsilon ]$ the set of polynomials with coefficients in $\ev$. More precisely it is the set of elements of $\ev^{\mathbb{N}}$ with finite support and, for any $(\varphi _{\ell } )\in \ev[\ln \varepsilon ]$, we denote $\fsum_{\ell \in \mathbb{N}} \ln ^{\ell } \!\varepsilon \,\varphi _{\ell } := (\varphi _{\ell } )$.
\end{definition}

\subsection{Definition of the spaces $\mathcal{A}$}
\label{subsec: def de A}

\begin{textesansboite}
In \cite[p.10]{CosDau96}, Costabel and Dauge build a similar asymptotic expansion for the Poisson equation in the half plane with mixed boundary condition: Neumann in a part of the boundary and Robin $u+\varepsilon \partial _{n} u=0$ in another. They quickly mention that their singularities can be written as $\Re[(-z)^{q} P(\ln (-z))]$ with $z=x+\mathrm{i}y$, $q\in \mathbb{R}$ and $P$ a real polynomial. To define the spaces $\mathcal{A}$, we adapted this idea to take into account the layer, the angle $\Theta $ and the Helmholtz equation. These simple expressions give both powerful algebraic tools for the theory and fast precise algorithms for the numerical resolution (see Section~\ref{subsec: operateurs de resolution dans A}).
\end{textesansboite}

\begin{definition}{the spaces $\mathcal{A}$\titreligne}
\label{def: espaces A}%
Let $\alpha  := e^{-\mathrm{i}\Theta }$. We define in $\Omega $ the complex variable $z:=x+\mathrm{i}y=re^{\mathrm{i}\theta }$. For all $q\in \mathbb{R}$, we take the following conventions: $(\alpha z)^{q} := r^{q} e^{q\mathrm{i}(\theta -\Theta )}$, $\overline{\alpha z}^{q} := r^{q} e^{-q\mathrm{i}(\theta -\Theta )}$ and $\log (\alpha z) := \ln r+\mathrm{i}(\theta -\Theta )$. Let $d\in \mathbb{R}$. We denote:
\begin{itemize}
	\item $\mathcal{A}_{d} (\Omega )$ the vector subspace of $\mathcal{C}^{0} (\Omega ,\mathbb{C})$ generated by the functions $z\mapsto \Im[(\alpha z)^{q} \,\overline{\alpha z}^{k} P(\log (\alpha z))]$ with $q\in \mathbb{R}$, $k\in \mathbb{N}$, $q+k=d$ and $P\in \mathbb{R}[T]$,
	\item $\mathcal{A}_{d} (\couche) := \{(x,Y)\mapsto  x^{d} \,Q(\ln x,Y) \mid Q\in \mathbb{C}[T,Y] \text{ and } Q(T,-1)=0\}$,
	\item $\mathcal{A}_{d} (\Gamma ) := \{x\mapsto  x^{d} \,Q(\ln x) \mid Q\in \mathbb{C}[T]\}$,
	\item $\mathcal{A}_{d} (\Pi ) := \{\varphi \in \mathcal{C}^{0} (\Pi ,\mathbb{C}) \mid \varphi _{|\Omega } \in \mathcal{A}_{d} (\Omega ) \text{ and } \varphi _{|\couche} \in \mathcal{A}_{d} (\couche)\}$,
	\item and for any $D\in \{\Pi ,\Omega ,\couche,\Gamma \}$, $\mathcal{A}(D) :=\displaystyle \sum _{d\in \mathbb{R}} \mathcal{A}_{d} (D)$ (cf. the introduction of Section~\ref{subsec: preliminaires algebriques}).
\end{itemize}\fbi
\end{definition}

\begin{textesansboite}
Note that elements of $\mathcal{A}(\Omega )$ and $\mathcal{A}(\Pi )$ vanish on $\bordgauche$, and elements of $\mathcal{A}(\couche)$ and $\mathcal{A}(\Pi )$ vanish on $\borddroit$. In addition, elements of $\mathcal{A}(\Omega )$ are naturally functions depending on the polar coordinates. For instance:
\begin{itemize}
	\item $\Im[(\alpha z)^{q} \,\overline{\alpha z}^{k} ] = r^{q+k} \sin \!\big((q-k)(\theta -\Theta )\big)$
	\item $\Im[(\alpha z)^{q} \,\overline{\alpha z}^{k} \log (\alpha z)] = r^{q+k} \left[ \ln (r)\,\cos \!\big((q-k)(\theta -\Theta )\big)-(\theta -\Theta )\sin \!\big((q-k)(\theta -\Theta )\big) \right]$
\end{itemize}
Note also that in this definition we used the variables $x,y$, which are relevant for far fields, but all the tools developed in this section can also be used for corner fields, replacing $(x,y)$ by $(X,Y)$.
\end{textesansboite}

\begin{paragraphesansboite}{Remark \nvnumpar{}:}
\label{rq: espaces A pour 2 couches}%
Definition~\ref{def: espaces A} can be extended to the case where $\bordgauche$ is covered by another layer, by defining $\mathcal{A}_{d} (\Omega )$ as the vector space generated by the functions $z\mapsto \Im[(\alpha z)^{q} \,\overline{\alpha z}^{k} P(\log (\alpha z))]$ and $z\mapsto \Im[z^{q} \bar{z}^{k} P(\log  z)]$ with $q\in \mathbb{R}$, $k\in \mathbb{N}$, $q+k=d$ and $P\in \mathbb{R}[T]$.
\end{paragraphesansboite}

\begin{textesansboite}
In order to build particular solutions of PDEs in $\mathcal{A}$, we will need the three following lemmas. The proof of the first one can be found in Appendix~\ref{sec: annexe A}.
\end{textesansboite}

\begin{lemme}{}
\label{decomposition de A en sev}%
For any $D\in \{\Pi ,\Omega ,\Gamma ,\couche\}$, we have the following decomposition: $\mathcal{A}(D) =\smash{\displaystyle \bigoplus_{d\in \mathbb{R}} \mathcal{A}_{d} (D)}$.\\[0.5mm]
Furthermore, for any $d\in \mathbb{R}$, $\mathcal{A}_{d} (\Omega )$ can itself be decomposed as follows: 
\begin{equation}{\label{eq: decomposition de A en sev 2}}
\mathcal{A}_{d} (\Omega ) = \bigoplus_{\substack{(q,k)\in \mathbb{R}\times \mathbb{N}\\ q+k=d}} \mathrm{\operatorname{Span}}_{\mathbb{C}} \big(\big\{z\mapsto  \Im[(\alpha z)^{q} \,\overline{\alpha z}^{k} P(\log (\alpha z))] \ \big|\ P\in \mathbb{R}[T] \text{ and } \mathcal{P}(q,k,P)\big\}\big)\vspace{-1mm} \end{equation}
where $\mathcal{P}$ is the property defined by $\mathcal{P}(q,k,P) := (q\not\in \mathbb{N} \text{ or } q>k \text{ or } P(0)=0)$.
\end{lemme}

\begin{remarque}
The condition $\mathcal{P}$ is a way to exclude the functions $z\mapsto  \Im[(\alpha z)^{q} \overline{\alpha z}^{k} ]$ with $q\in \mathbb{N}$ and $q<k$, which are already present in the direct sum as they are equal to $z\mapsto  -\Im[(\alpha z)^{k} \overline{\alpha z}^{q} ]$.
\end{remarque}

\begin{textesansboite}
Let $\varphi $ be a function of $\mathcal{A}(\Omega )$ of the form $\Im[(\alpha z)^{q} \,\overline{\alpha z}^{k} P(\log (\alpha z))]$. Note that on $\Gamma $, $\varphi $ is equal to $x^{q+k} \Im[\alpha ^{q -k } P(\ln x-\mathrm{i}\Theta )]$. Let us define $\Im[\alpha ^{q-k} P(T-\mathrm{i}\Theta )] :=\sum _{i=1} ^{\deg P} \Im(a_{i} ) T^{i} $ in $\mathbb{C}[T]$, where $\sum _{i=1} ^{\deg P} a_{i} T^{i} := \alpha ^{q-k} P(T-\mathrm{i}\Theta )$. Then $\varphi _{|\Gamma } (x) = x^{q+k} Q(\ln x)$ for some $Q\in \mathbb{R}[X]$, which implies $\varphi _{|\Gamma } \in \mathcal{A}(\Gamma )$. Conversely, for future constructions, it will be important to solve the equation:
\begin{equation}{\label{eq: resolution de Im[alphad P(T-iTheta)]=Q}}
\text{given $Q\in \mathbb{R}[T]$, \;find $P\in \mathbb{R}[T]$ \;s.t. \ $\Im[\alpha ^{d} P(T-\mathrm{i}\Theta )]=Q(T)$.} \end{equation}\fbeq
\end{textesansboite}

\begin{lemme}{}
\label{resolution de Im alphad P de T moins iTheta egal Q}%
Let $\alpha :=e^{-\mathrm{i}\Theta }$, $d\in \mathbb{R}$ and $Q\in \mathbb{R}[T]$.
\begin{enumerate}
	\item If $d\in \mathbb{R} \setminus \TZ$, then there is a unique solution $P\in \mathbb{R}[T]$ of \eqref{eq: resolution de Im[alphad P(T-iTheta)]=Q}. Moreover $\deg P=\deg Q$. We denote the solution $\RIma{d}(Q)$.
	\item If $d\in \TZ$, then the set of solutions of \eqref{eq: resolution de Im[alphad P(T-iTheta)]=Q} is of the form $\{P_{0} +c \mid c\in \mathbb{R}\}$ with $P_{0} \in \mathbb{R}[T]$ and $\deg (P_{0} )=\deg (Q)+1$. We denote $\RIma{d}(Q)$ the unique solution that vanishes at 0.
\end{enumerate}
In both cases, $\RIma{d}$ is a linear map from $\mathbb{R}[T]$ into itself.
\end{lemme}

\begin{demo}
There are two cases whether the coefficient of degre $m$ of $\Im[\alpha ^{d} (T-\mathrm{i}\Theta )^{m} ]$ vanishes or not.
\begin{enumerate}
	\item If $d\in \mathbb{R} \setminus \TZ$, then $\alpha ^{d} \in \mathbb{C} \setminus \mathbb{R}$, so: $\forall m\in \mathbb{N},\; \deg \Im[\alpha ^{d} (T-\mathrm{i}\Theta )^{m} ] =m$. Therefore, $(\Im[\alpha ^{d} (T-\mathrm{i}\Theta )^{m} ])_{m\in \mathbb{N}}$ is a basis of $\mathbb{R}[T]$. So writing $Q$ in this basis gives a unique solution of \eqref{eq: resolution de Im[alphad P(T-iTheta)]=Q}.
	\item If $d\in \TZ$, since $\alpha ^{d} \in \mathbb{R}$, we have: $\forall m\in \mathbb{N},\; \deg \Im[\alpha ^{d} (T-\mathrm{i}\Theta )^{m} ] =m-1$. So in this case $(\Im[\alpha ^{d} (T-\mathrm{i}\Theta )^{m} ])_{m\in \mathbb{N}^*}$ is a basis of $\mathbb{R}[T]$. Thus, $F:P\mapsto \Im[\alpha ^{d} P(T-\mathrm{i}\Theta )]$ is surjective and its kernel is the set of constant polynomials. Its restriction to $E:=\{P\in \mathbb{R}[T] \mid P(0)=0\}$ is therefore an isomorphism and we set $\RIma{d}:= (F_{|E} )^{-1}$. Finally, for any $Q\in \mathbb{R}[T]$ we have $F^{-1} (\{Q\})= \{\RIma{d}(Q)\}+\operatorname{Ker} F$.
\end{enumerate}\fbi
\end{demo}

\begin{textesansboite}
For any $q\in \mathbb{R}$, the maps $Q\in \mathbb{C}[T] \mapsto  \Im[x^{q} Q(\ln x)]\in \mathcal{A}(\Gamma )$ and $Q\in \mathbb{C}[T,Y] \mapsto  \Im[x^{q} Q(\ln x,Y)]\in \mathcal{A}(\couche)$ are clearly injective. In the following lemma we investigate the injectivity of $P\in \mathbb{R}[T] \mapsto  \Im[(\alpha z)^{q} \,\overline{\alpha z}^{k} P(\log (\alpha z))]\in \mathcal{A}(\Omega )$.
\end{textesansboite}

\begin{lemme}{}
\label{unicite de P et Q}%
Let $(q,k)\in \mathbb{R}\times \mathbb{N}$. The map $P \mapsto  \Im[(\alpha z)^{q} \,\overline{\alpha z}^{k} P(\log (\alpha z))]$ is injective from the set of real polynomials $P$ for which $\mathcal{P}(q,k,P)$ is true into $\mathcal{A}(\Omega )$.
\end{lemme}

\begin{demo}
For any $\theta \in (0,\Theta )$ and $r\in \mathbb{R}_+^*$, we have $\varphi (re^{\mathrm{i}\theta } )=r^{q+k} \Im[e^{\mathrm{i}(q-k)(\theta -\Theta )} P(\ln r+\mathrm{i}(\theta -\Theta ))] = r^{q+k} \Im[\alpha '^{q-k} P(\ln r-\mathrm{i}\Theta ')]$ with $\Theta ':=\Theta -\theta $ and $\alpha ' :=e^{-\mathrm{i}\Theta '}$.
\begin{itemize}
	\item If $q\neq k$, we can choose $\theta $ so that $q-k \in \mathbb{R} \setminus \frac{\pi }{\Theta '} \mathbb{Z}$. So Lemma~\ref{resolution de Im alphad P de T moins iTheta egal Q} applied to $(\Theta ',\alpha ')$ instead of $(\Theta ,\alpha )$ implies that $P$ is unique.
	\item Otherwise, $q-k=0 \in  \frac{\pi }{\Theta '} \mathbb{Z}$ for any $\theta $. So according to \ref{resolution de Im alphad P de T moins iTheta egal Q}, $P$ is unique up to a constant a priori. But the property $\mathcal{P}$ implies that $P(0)=0$, so this constant is fixed.
\end{itemize}\fbi
\end{demo}

\subsection{Tools for solving the Poisson and Helmholtz equations in the spaces $\mathcal{A}$}
\label{subsec: operateurs de resolution dans A}

\newcommand{\inc}{\varphi}
\newcommand{\sm}{\psi}
\begin{textesansboite}
In this section, we show how to solve canonical problems set in $\Pi $ in the spaces $\mathcal{A}$. More precisely, let $(\sm_{\Omega } ,\sm_{\couche} ,\sm_{\Gamma } ) \in \mathcal{A}(\Omega )\times \mathcal{A}(\couche) \times \in \mathcal{A}(\Gamma )$, we look for the solutions $\inc \in \mathcal{A}(\Pi )$ of:
problems of the form
\begin{equation}{\label{eq: pb resolu dans A}}
\left\{\begin{array}{r@{\;=\;}ll}
	\Delta  \inc & \sm_{\Omega } &\text{ in }\Omega  \\[0.5mm]
	\partial _{Y} ^2 \inc & \sm_{\couche} &\text{ in }\couche \\[0.5mm]
	\partial _{Y} \inc_{|Y=0^{-}} & \sm_{\Gamma } &\text{ on }\Gamma  
\end{array}\right.\end{equation}
Note that by definition of $\mathcal{A}(\Pi )$, $\inc$ also satisfies $\inc_{|y=0^{+}} - \inc_{|Y=0^{-}}$ on $\Gamma $ and $\inc_{|\bordgauche \cup \borddroit} =0$. Solving this system will enable us to build in Section~\ref{subsec: DA des champs lointains et de coin} the asymptotic expansion of $\bfu_{p,\ell }$ and $S_{p,\ell }$. Indeed, note for instance that this system is identical to \eqref{eq: eq vol de bfu p,l} except for the first line. We first describe the homogeneous solutions of \eqref{eq: pb resolu dans A}, then build explicitly some particular solutions. Since functions of $\mathcal{A}$ are uniquely determined by some polynomials (see Lemma~\ref{unicite de P et Q}), we are able to code an \emph{exact, very fast and memory-thrifty} solver of \eqref{eq: pb resolu dans A}. This is one of the key advantages of the $\mathcal{A}$ framework.
\end{textesansboite}

\begin{definition}{}
\label{def: phi d}%
For any $d\in \TZ*$, we define on $\Omega $ the function $\phiOmega_{d} := (-1)^{d\Theta /\pi } r^{d} \sin (d\theta ) = r^{d} \sin \!\big(d(\theta -\Theta )\big) =\Im[(\alpha z)^{d} ] \in \mathcal{A}_{d} (\Omega )$ and $\phi_{d} \in \mathcal{A}_{d} (\Pi )$ its extension by 0 in $\couche$.
\end{definition}

\begin{textesansboite}
These functions play an important role in the sequel because they solve the homogeneous Laplace equation in $\Omega $, resp. $\Pi $.
\end{textesansboite}

\begin{proposition}{the Laplace problem in $\mathcal{A}$\titreligne}
\label{pb de Laplace homogene dans A}\vspace{-5mm}%
\begin{enumerate}
	\item $\operatorname{Span}(\{\phiOmega_{d} \mid d\in \TZ*\})$ is the set of solutions in $\mathcal{A}(\Omega )$ of 
	\begin{equation}{\label{eq: pb de Laplace homogene dans A(Omega)}}
	\left\{\begin{array}{r@{\;=\;}ll}
		\Delta  \varphi  & 0 &\text{ in }\Omega \\[0.5mm]
		\varphi  & 0 &\text{ on }\Gamma \cup \bordgauche
	\end{array}\right.\end{equation}
	\item $\operatorname{Span}(\{\phi_{d} \mid d\in \TZ*\})$ is the set of solutions in $\mathcal{A}(\Pi )$ of 
	\begin{equation}{\label{eq: pb de Laplace homogene dans A(Pi)}}
	\left\{\begin{array}{r@{\;=\;}ll}
		\Delta  \varphi  & 0 &\text{ in }\Omega \\[0.5mm]
		\partial _{Y} ^2 \varphi  & 0 &\text{ in }\couche \\[0.5mm]
		\partial _{Y} \varphi _{|Y=0^{-}} & 0 &\text{ on }\Gamma 
	\end{array}\right.\end{equation}
\end{enumerate}\fbeq
\end{proposition}

\begin{demo}
\begin{enumerate}
	\item $\phiOmega_{d}$ is clearly solution of \eqref{eq: pb de Laplace homogene dans A(Omega)} for any $d\in \TZ*$. Conversely, let $\varphi $ be a solution.	Let us denote 
	\[\forall d\in {\textstyle \TN},\;\forall r\in \mathbb{R}_+^*, \quad  c_{d} (r) := \frac{2}{\Theta } \int _{0} ^{\Theta } \varphi (r,\theta )\,\sin (d\theta )\,\mathrm{d}\theta .\]
	Using separation of variables and $\varphi _{|\partial \Omega } =0$, it is easy to show that  $\varphi =\sum _{d\in \TN} c_{d} (r) \sin (d\theta )$ with convergence in $H^2 (\Omega \cap \{r_{1} <r<r_{2} \})$ for any $0<r_{1} <r_{2} <\infty $. Since $\Delta \varphi =0$, we have $(r \frac{\mathrm{d}}{\mathrm{d}r} )^2 c_{d} = d^2 c_{d}$ for any $d$. Hence we get: $\forall d\in \TN, \exists a_{d} ,a_{-d} \in \mathbb{C}, \forall r\in \mathbb{R}_+^*,\; c_{d} (r) = a_{d} r^{d} +a_{-d} r^{-d}$.\\
	Moreover, by definition of $\mathcal{A}(\Omega )$, there is $q\in \mathbb{R}_+^*$ s.t. $\varphi \in \sum _{d\in [-q,q]} \mathcal{A}_{d} (\Omega )$. So\vspace{-1mm}
	\[\forall d\in {\textstyle \TZ*} \cap (q,\infty ),\quad  a_{d} = \lim _{r\rightarrow \infty } r^{-d} c_{d} (r) = \lim _{r\rightarrow \infty } \frac{2}{\Theta } \int _{0} ^{\Theta } \underbrace{r^{-d} \varphi (r,\theta )}_{\longrightarrow \,0}\,\sin (d\theta )\,\mathrm{d}\theta  =0. \vspace{-1mm}\]
	Similarly, looking at $r\rightarrow 0$ one gets: $\forall d\in \TZ* \cap (-\infty ,-q),\; a_{d} =0$. So $\varphi  = \sum _{d\in \TZ*\cap [-q,q]} a_{d} \,\mathrm{sgn} (d) \,\phiOmega_{d}$ (where $\mathrm{sgn} (d):=d/|d|$), which is a finite sum. Therefore, $\varphi $ is in the desired span.
	\item Any solution of the system vanishes in $\couche$, so point 2 easily follows from point 1.
\end{enumerate}
\end{demo}

\begin{textesansboite}
Let us now define the following linear forms $\sigma _{d}$ which satisfy $\sigma _{d} (\phi_{q} )=\delta _{d,q}$ for any $d,q\in \TZ*$ and which enable us to ``project'' any element of $\mathcal{A}(\Pi )$ on $\operatorname{Span}(\{\phi_{d} \mid d\in \TZ*\})$. These linear forms appear later as key singularity coefficients in the matching condition.
\end{textesansboite}

\begin{definition}{linear forms $\sigma _{d}$\titreligne}
\label{def: sigma}%
Let $d\in \TZ*$. For any $(q,k,P)\in \mathbb{R}\times \mathbb{N}\times \mathbb{R}[T]$ s.t. $\mathcal{P}(q,k,P)$ is true and $\varphi :z \mapsto  \Im[(\alpha z)^{q} \,\overline{\alpha z}^{k} P(\log (\alpha z))]$, let us define: \[\left\{\begin{array}{r@{\;=\;}ll}
	\sigma _{d} (\varphi ) &0 & \text{if } q\neq d \text{ or } k\neq 0\\[0.5mm]
	\sigma _{d} \big( \Im[(\alpha z)^{d} P(\log (\alpha z))] \big) & P(0) & \text{otherwise} 
\end{array}\right.\] It is well-defined by Lemma~\ref{unicite de P et Q}. By Lemma~\ref{decomposition de A en sev}, $\sigma _{d}$ can be extended into a linear form $\mathcal{A}(\Omega ) \rightarrow \mathbb{C}$. Finally, for any $\varphi \in \mathcal{A}(\Pi )$, $\sigma _{d} (\varphi ):=\sigma _{d} (\varphi _{|\Omega } )$.
\end{definition}

\begin{textesansboite}
Let us now build particular solutions of \eqref{eq: pb resolu dans A}. By linearity, it suffices to build particular solutions of three sub-problems. According to Definition~\ref{def: degre d'un operateur}, for any $D_{1} ,D_{2} \in \{\Pi ,\Omega ,\Gamma ,\couche\}$, we say that a linear map $F:\mathcal{A}(D_{1} )\rightarrow \mathcal{A}(D_{2} )$ has degree $d\in \mathbb{R}$ iff: $\forall q\in \mathbb{R},\, \forall \varphi \in \mathcal{A}_{q} (D_{1} ),\; F(\varphi )\in \mathcal{A}_{q+d} (D_{2} )$.
\end{textesansboite}

\newcommand{\smOmega}{\sm_\Omega}
\newcommand{\smCouche}{\sm_\couche}
\newcommand{\smGamma}{\sm_\Gamma}
\newcommand{\incOmega}{\inc_{\Delta}}
\newcommand{\incCouche}{\inc_{\partial_Y^2}}
\newcommand{\incGamma}{\inc_{\mathrm{N}}}
\begin{boitecarreesecable}{Proposition \nvnumpar{}: particular solutions of \eqref{eq: pb resolu dans A}\titreligne}
\label{operateurs R}%
Let us denote $\Aperp(\Pi ) := \{\varphi \in \mathcal{A}(\Pi ) \mid \forall d\in \TZ*,\; \sigma _{d} (\varphi )=0\}$, which is a supplementary of $\operatorname{Span}(\{\varphi _{d} \mid d\in \TZ*\})$ in $\mathcal{A}(\Pi )$.
\begin{enumerate}
	\item For any $\smOmega \in \mathcal{A}(\Omega )$ there exists a unique solution $\incOmega \in \Aperp(\Pi )$ of
	\begin{equation}{\label{eq: operateur R Omega Pi}}
	\left\{\begin{array}{r@{\;=\;}ll}
		\Delta  \inc & \smOmega &\text{ in }\Omega \\[0.5mm]
		\partial _{Y} ^2 \inc & 0 &\text{ in }\couche \\[0.5mm]
		\partial _{Y} \inc_{|Y=0^{-}} & 0 &\text{ on }\Gamma 
	\end{array}\right. \end{equation}
	The associated map $\ROmegaPi:\smOmega \in \mathcal{A}(\Omega ) \mapsto \incOmega \in \Aperp(\Pi )$ is linear and has degree 2.
	\item For any $\smCouche \in \mathcal{A}(\couche)$ there exists a unique solution $\incCouche \in \Aperp(\Pi )$ of
	\begin{equation}{\label{eq: operateur R couche Pi}}
	\left\{\begin{array}{r@{\;=\;}ll}
		\Delta  \inc & 0 &\text{ in }\Omega \\[0.5mm]
		\partial _{Y} ^2 \inc & \smCouche &\text{ in }\couche \\[0.5mm]
		\partial _{Y} \inc_{|Y=0^{-}} & 0 &\text{ on }\Gamma 
	\end{array}\right. \end{equation}
	The associated map $\RcouchePi:\smCouche \in \mathcal{A}(\couche) \mapsto \incCouche \in \Aperp(\Pi )$ is linear and has degree 0.
	\item For any $\smGamma \in \mathcal{A}(\Gamma )$ there exists a unique solution $\incGamma \in \Aperp(\Pi )$ of
	\begin{equation}{\label{eq: operateur R Gamma Pi}}
	\left\{\begin{array}{r@{\;=\;}ll}
		\Delta  \inc & 0 &\text{ in }\Omega \\[0.5mm]
		\partial _{Y} ^2 \inc & 0 &\text{ in }\couche \\[0.5mm]
		\partial _{Y} \inc_{|Y=0^{-}} & \smGamma &\text{ on }\Gamma 
	\end{array}\right. \end{equation}
	The associated map $\RGammaPi:\smGamma \in \mathcal{A}(\Gamma ) \mapsto \incGamma \in \Aperp(\Pi )$ is linear and has degree 0.
\end{enumerate}\fbi
\end{boitecarreesecable}

\begin{textesansboite}
Using Propositions~\ref{pb de Laplace homogene dans A} and \ref{operateurs R}, it is then easy to see that the set of solutions of \eqref{eq: pb resolu dans A} is $\ROmegaPi(\varphi _{\Omega } )+\RcouchePi(\varphi _{\couche} )+\RGammaPi(\varphi _{\Gamma } )+\operatorname{Span}(\{\phi_{d} \mid d\in \TZ*\})$. Moreover, the functions $\incOmega$, $\incCouche$ and $\incGamma$ in Proposition~\ref{operateurs R} have explicit expressions (see the proof below), which allows to compute them easily in practice.
\end{textesansboite}

\begin{demo}
Proposition~\ref{pb de Laplace homogene dans A} gives the uniqueness of the solutions $\incOmega$, $\incCouche$ and $\incGamma$, so only their existence remains to prove. This is done by a construction. For any $D\in \{\Pi ,\Omega ,\couche,\Gamma \}$, let $\mathcal{A}(D,\mathbb{R}) := \mathcal{A}(D)\cap \mathcal{C}^{0} (D,\mathbb{R})$. Since $\mathcal{A}(D)=\mathcal{A}(D,\mathbb{R}) \oplus  \mathrm{i}\mathcal{A}(D,\mathbb{R})$, it suffices to build the solutions when $(\smOmega,\smCouche,\smGamma)\in \mathcal{A}(\Omega ,\mathbb{R})\times \mathcal{A}(\couche,\mathbb{R})\times \mathcal{A}(D,\Gamma )$, and then extend it to any source term by complexification.
\begin{enumerate}
	\item According to Lemma~\ref{decomposition de A en sev}, it suffices to build $\incOmega$ when $\smOmega= \Im[(\alpha z)^{q} \,\overline{\alpha z}^{k} P_{\sm} (\log (\alpha z))]$ with $(q,k,P_{\sm} )\in \mathbb{R}\times \mathbb{N}\times \mathbb{R}[T]$ s.t. $\mathcal{P}(q,k,P_{\sm} )$ is true. First, ${\incOmega}_{|\couche} = 0$ because $\incOmega$ satisfies: \[\left\{\begin{array}{r@{\;=\;}ll}
		\partial _{Y} ^2 \incOmega & 0 &\text{ in }\couche \\[0.5mm]
		\partial _{Y} {\incOmega}_{|Y=0^{-}} & 0 &\text{ on }\Gamma \\[0.5mm]
		\incOmega & 0 &\text{ on }\borddroit
	\end{array}\right.\]
	Given that $\Delta  =4\partial _{z} \partial _{ \bar{z}}$, we have for any $\inc_{1} :z\mapsto  \Im[(\alpha z)^{q_{1}} \,\overline{\alpha z}^{k_{1}} P_{1} (\log (\alpha z))]$ that
	\[\Delta \inc_{1} = 4\,\Im[(\alpha z)^{q_{1} -1} \, k_{1} \overline{\alpha z}^{k_{1} -1} (q_{1} P_{1} +P_{1} ')(\log (\alpha z))].\]
	So taking $q_{1} :=q+1$, $k_{1} :=k+1$ and $P_{1} \in \mathbb{R}[T]$ a solution of $4k_{1} (q_{1} P_{1} +P_{1} ') = P_{\sm}$, we have $\Delta \inc_{1} = \smOmega$ in $\Omega $.
	\begin{itemize} 
		\item If $q_{1} =0$, $P_{1}$ is unique up to a constant, Moreover we can write
		\[\varphi _{1} = \Im\big[\overline{\alpha z}^{k_{1}} \big(P_{1} (\log (\alpha z))-P_{1} (0)\big)\big] - P_{1} (0)\, \Im[(\alpha z)^{k_{1}} ]\]
		where each term satisfy the property $\mathcal{P}$. So for any $d\in \TZ*$ different from $k_{1}$, we have $\sigma _{d} (\varphi _{1} )=0$, while $\sigma _{k_{1}} (\varphi _{1} )=-P_{1} (0)$ if $k_{1} \in \mathbb{N}\cap \TZ*$. Taking $P_{1} (0):=0$ thus gives: $\forall d\in \TZ*,\; \sigma _{d} (\inc_{1} )=0$.
		\item If $q_{1} \neq 0$, there is a unique solution $P_{1}$. Given that $q_{1} \neq 0$ and $k_{1} \neq 0$, Definition~\ref{def: sigma} implies that: $\forall d\in \TZ*,\; \sigma _{d} (\inc_{1} )=0$.
	\end{itemize}
	However, we cannot set ${\incOmega}_{|\Omega } = \inc_{1}$, because $\varphi _{1}$ does not vanish on $\Gamma $. Let us then introduce $\inc_{2} : z\mapsto  \Im[(\alpha z)^{q+k+2} P_{2} (\log (\alpha z))]$ with $P_{2} (T) := \RIma{q+k+2} \big(\Im[\alpha ^{q-k} P_{1} (T-\mathrm{i}\Theta )]\big)$ that satisfies by Lemma~\ref{resolution de Im alphad P de T moins iTheta egal Q}:
	\[\left\{\begin{array}{r@{\;=\;}ll}
		\Delta  \inc_{2} & 0 &\text{ in }\Omega  \\[0.5mm]
		\inc_{2} & x^{q+k+2} \Im[\alpha ^{q+k+2} P_{2} (\ln x-\mathrm{i}\Theta )] = x^{q+k+2} \Im[\alpha ^{q-k} P_{1} (\ln x-\mathrm{i}\Theta )] =\inc_{1}  &\text{ on }\Gamma 
	\end{array}\right.\]
	In addition Lemma~\ref{resolution de Im alphad P de T moins iTheta egal Q} implies that $q+k+2\in \TZ*\Rightarrow P_{2} (0)=0$, so: $\forall d\in \TZ*,\; \sigma _{d} (\inc_{2} )=0$. Finally we set ${\incOmega}_{|\Omega } := \inc_{1} -\inc_{2}$, which is in $\Aperp(\Pi )$.
	
	\item Similarly, it suffices to build $\incCouche$ for $\smCouche= x^{q} \,Q_{\sm} (\ln x,Y)$. We look for $\incCouche$ of the form $x^{q} \,Q_{\inc} (\ln x,Y)$ in $\couche$ with $Q_{\inc} \in \mathbb{R}[T,Y]$. Then necessarily we have \[\left\{\begin{array}{r@{\;=\;}l}
		\partial _{Y} ^2  Q_{\inc} & Q_{\sm} \\[0.5mm]
		\partial _{Y} Q_{\inc} (\cdot ,0)&0 \\[0.5mm]
		Q_{\inc} (\cdot ,-1) & 0
	\end{array}\right.\]
	This uniquely defines $Q_{\inc}$. Taking $P_{\inc} := \RIma{q}(Q_{\inc} (\cdot ,0))$ and ${\incCouche}_{|\Omega } := \Im[(\alpha z)^{q} P_{\inc} (\log (\alpha z))]$ then implies that $\Delta \incCouche =0$ in $\Omega $, $\inc$ continuous and $\incCouche \in \Aperp(\Pi )$.
	
	\item Again it suffices to consider $\smGamma= x^{q} \,P_{\sm} (\ln x)$ with $P_{\sm} \in \mathbb{R}[T]$. We take $\incGamma$ of the same form as in point 2. Then $\partial _{Y} {\incGamma}_{|Y=0^{-}} = x^{q} \,\partial _{Y} Q_{\inc} (\ln x,0)$ so it suffices to take the polynomial solutions of: \[\left\{\begin{array}{r@{\;=\;}l}
		\partial _{Y} ^2  Q_{\inc} & 0 \\[0.5mm]
		\partial _{Y} Q_{\inc} (\cdot ,0) & P_{\sm} \\[0.5mm]
		Q_{\inc} (\cdot ,-1) & 0\\[0.5mm]
		P_{\inc} & \RIma{q}(Q_{\inc} (\cdot ,0))
	\end{array}\right.\]
\end{enumerate}
Finally the linearity of Problems~\ref{eq: operateur R Omega Pi}--\ref{eq: operateur R Gamma Pi} and the uniqueness of $\incOmega$, $\incCouche$ and $\incGamma$ imply that $\ROmegaPi$, $\RcouchePi$ and $\RGammaPi$ are linear maps.
\end{demo}

\begin{textesansboite}
Moreover, we will need analogous operators in $\mathcal{A}(\Omega )$ to build the far fields in Theorem~\ref{cadre fonctionnel des champs lointains}. The proof is entirely similar to the one of Proposition~\ref{operateurs R}, so we omit it. Again, these operators have explicit expressions.
\end{textesansboite}

\renewcommand{\incGamma}{\inc_{\mathrm{D}}}
\begin{boitecarreesecable}{Proposition \nvnumpar{}: }
\label{operateurs R dans Omega}%
Let us denote $\Aperp(\Omega ) = \{\varphi \in \mathcal{A}(\Omega ) \mid \forall d\in \TZ*,\; \sigma _{d} (\varphi )=0\}$, which is a supplementary of $\operatorname{Span}(\{\phiOmega_{d} \mid d\in \TZ*\})$ in $\mathcal{A}(\Omega )$.
\begin{enumerate}
	\item For any $\smOmega \in \mathcal{A}(\Gamma )$ there exists a unique solution $\incOmega \in \Aperp(\Omega )$ of
	\begin{equation}{\label{eq: operateur R couche Omega}}
	\left\{\begin{array}{r@{\;=\;}ll}
		\Delta  \inc & \smOmega &\text{ in }\Omega \\[0.5mm]
		\inc & 0 &\text{ on }\Gamma 
	\end{array}\right. \end{equation}
	The associated map $\ROmegaOmega:\smGamma \in \mathcal{A}(\Omega ) \mapsto \incOmega \in \Aperp(\Omega )$ is linear and has degree 2.
	\item For any $\smGamma \in \mathcal{A}(\Gamma )$ there exists a unique solution $\incGamma \in \Aperp(\Omega )$ of
	\begin{equation}{\label{eq: operateur R Gamma Omega}}
	\left\{\begin{array}{r@{\;=\;}ll}
		\Delta  \inc & 0 &\text{ in }\Omega \\[0.5mm]
		\inc & \sm &\text{ on }\Gamma 
	\end{array}\right. \end{equation}
	The associated map $\RGammaOmega:\smGamma \in \mathcal{A}(\Gamma ) \mapsto \incGamma \in \Aperp(\Omega )$ is linear and has degree 0.
\end{enumerate}\fbi
\end{boitecarreesecable}

\begin{textesansboite}
To end this section, let us show that the spaces $\mathcal{A}$ are stable under some differential operators.
\end{textesansboite}

\begin{lemme}{}
\label{operateurs de derivation dans A}%
\begin{enumerate}
	\item $\partial _{x|\couche} ^2  :\varphi  \mapsto  \partial _{x} ^2  \varphi _{|\couche}$ maps $\mathcal{A}(\Pi )$ to $\mathcal{A}(\couche)$ and has degree $-2$ (see Definition~\ref{def: degre d'un operateur}).
	\item $\partial _{y|\Gamma ,y=0^{+}} : \varphi  \mapsto  \partial _{y} \varphi _{|\Gamma ,y=0^{+}}$ maps $\mathcal{A}(\Pi )$ to $\mathcal{A}(\Gamma )$ and has degree $-1$.
	\item $\Delta :\varphi \mapsto  \Delta \varphi _{|\Omega }$ maps $\mathcal{A}(\Omega )$ (and $\mathcal{A}(\Pi )$) to $\mathcal{A}(\Omega )$ and has degree $-2$.
\end{enumerate}\fbi
\end{lemme}

\begin{textesansboite}
Note that, when used on corner fields (which depend on $(X,Y)$), the first two operators will rather be denoted $\partial _{X|\couche} ^2 $ and $\partial _{Y|\Gamma ,Y=0^{+}}$.
\end{textesansboite}

\begin{demo}
It suffices to verify it when $\varphi $ has the form $\Im[(\alpha z)^{q} \,\overline{\alpha z}^{k} P(\log (\alpha z)) ]$ in $\Omega $ and $x^{q+k} \,Q(\ln x,Y)$ in $\couche$, with $q\in \mathbb{R}$, $k\in \mathbb{N}$, $P[T]$ and $Q\in \mathbb{R}[T,Y]$.
\begin{enumerate}
	\item We have $\partial _{x} ^2 \varphi (x,Y) = x^{q+k-2} ((q+k)(q+k-1)+2(q+k)\partial _{T} +\partial _{T} ^2 )\,Q(\ln x,Y)$ which is in $\mathcal{A}(\couche)$.
	\item Let $\psi : z\mapsto  (\alpha z)^{q} \,\overline{\alpha z}^{k} P(\log (\alpha z))$. From $\partial _{z} \psi  = \frac{1}{2} (\partial _{x} \psi -\mathrm{i}\partial _{y} \psi )$ and $\partial _{ \bar{z}} \psi  = \frac{1}{2} (\partial _{x} \psi +\mathrm{i}\partial _{y} \psi )$ it follows that $\partial _{y} \varphi  = \Im[\partial _{y} \psi ] = \Re[\partial _{z} \psi  -\partial _{ \bar{z}} \psi ]$. Therefore, \begin{align*}
		\partial _{y} \varphi _{|\Gamma ,y=0^{+}} &= \Re[\alpha  (\alpha x)^{q-1} \,\overline{\alpha x}^{k} (qP+P')(\log (\alpha x)) -(\alpha x)^{q} \, \overline{\alpha }\,k\,\overline{\alpha x}^{k-1} P(\log (\alpha x))] \\
		&= x^{q+k-1} \Re[\alpha ^{q-k} ((q-k)P+P')(\ln x-\mathrm{i}\Theta )] \ \ \in \ \mathcal{A}(\Gamma )
	\end{align*}
	\item Since $\Delta  =4\partial _{z} \partial _{ \bar{z}}$, we have $\Delta \varphi  = 4\,\Im[(\alpha z)^{q-1} \,k\overline{\alpha z}^{k-1(qP_{1}} +P_{1} ')(\log (\alpha z))]$ which is in $\mathcal{A}(\Omega )$.
\end{enumerate}\fbi
\end{demo}

\subsection{Asymptotic behaviors w.r.t. $\bfr$ of solutions of model problems}
\label{subsec: DA des champs lointains et de coin}

\begin{textesansboite}
In this section we give tools that will be used to compute the asymptotic behaviors of $\bfu_{p,\ell }$ and $S_{p,\ell }$ resp. when $\bfr\rightarrow 0$ and $\bfr\rightarrow \infty $ using the spaces $\mathcal{A}$. To do so, we will use series of elements of $\mathcal{A}$, which is made rigorous by the following definition.
\end{textesansboite}

\begin{definition}{the spaces $\Amaicro$\titreligne}
\label{def: A barre}%
Let $D\in \{\Pi ,\Omega ,\Gamma ,\couche\}$. We denote $\Amacro(D)$, resp. $\Amicro(D)$, the set of elements of $\fsum_{d\in \mathbb{R}} \mathcal{A}_{d} (D)$ whose support is included in the image of a sequence that tends to $\infty $, resp. $-\infty $. We write their elements as formal series according to Definition~\ref{def: serie formelle}.
\end{definition}

\begin{remarques}
\begin{itemize}
	\item The asymptotic of $\bfu_{p,\ell }$ when $\bfr\rightarrow 0$ involve increasing powers of $r$ so it will be expressed in $\Amacro(D)$. Similarly, $S_{p,\ell }$ when $\bfr\rightarrow \infty $ involves decreasing powers of $r$, so it will be expressed in $\Amicro(D)$.
	\item As seen in Section~\ref{subsec: preliminaires algebriques}, $\mathcal{A}(D)$ is included in $\Amacro(D)$ and $\Amicro(D)$. But elements of $\Amaicro(D)$ are not in general $D\rightarrow \mathbb{R}$ functions, as the formal series may diverge pointwise.
\end{itemize}\fbi
\end{remarques}

\begin{textesansboite}
Using Definitions~\ref{def: degre d'un operateur}, we can extend $\partial _{x|\couche} ^2 $, $\partial _{y|\Gamma ,y=0^{+}}$, $\RcouchePi$, $\RGammaPi$ and $\ROmegaPi$ to the spaces $\Amaicro$. We also use the notation $\langle .\rangle $ introduced in Definition~\ref{def: serie geometrique d'operateurs 1D}. E.g. $\langle  -k_{0} ^2 \,\ROmegaPi \rangle  = \sum_{n=0}^\infty  (-k_{0} ^2 \,\ROmegaPi)^{n}$. \phantomsection\label{def: sigma dans Amaicro}Moreover, we extend $\sigma _{d}$ to $\Amaicro(\Pi )$ for any $d\in \TZ*$, by setting $\sigma _{d} (\fsum_{d'} \varphi _{d'} ):=\sigma _{d} (\varphi _{d} )$ for any $\fsum_{d'} \varphi _{d'} \in \Amaicro(\Pi )$.
\end{textesansboite}

\begin{lemme}{}
\label{DA algebrique des champs lointains concatenes}%
Let $g^{0} \in \Amacro(\couche)$ and $h^{0} \in \Amacro(\Gamma )$. The solutions in $\Amacro(\Pi )$ of 
\begin{equation}{\label{eq: equations de bfuDA}}
\left\{\begin{array}{r@{\;=\;}ll}
\mu _{0} \Delta  \bfuDA +\omega ^2 \rho _{0} \bfuDA & 0 &\text{ in }\Omega \\[0.5mm]
\mu _{1} \partial _{Y} ^2  \bfuDA & g^{0} &\text{ in }\couche \\[0.5mm]
\mu _{1} \partial _{Y|Y=0^{-}} \bfuDA & h^{0} &\text{ on }\Gamma 
\end{array}\right.\end{equation}
are the formal series of the following form (where $\sigma _{d} (\bfuDA)$ vanishes when $d$ is small enough)
\begin{equation}{\label{eq: bfuDA}}
\bfuDA = \big\langle  {-}k_{0} ^2 \,\ROmegaPi \big\rangle  \bigg( \frac{1}{\mu _{1}} \RcouchePi(g^{0} ) + \frac{1}{\mu _{1}} \RGammaPi(h^{0} ) + \fsum_{d\in \TZ*} \sigma _{d} (\bfuDA) \,\phi_{d} \bigg) \end{equation}\fbeq[1.3]
\end{lemme}

\begin{demo}
Let $\bfv^{0} := ( \mathrm{id}+k_{0} ^2  \ROmegaPi)\bfuDA$. Since $\langle -k_{0} ^2  \ROmegaPi \rangle  = \sum_{n=0}^\infty  (-k_{0} ^2  \ROmegaPi )^{n}$, $ \mathrm{id}+k_{0} ^2  \ROmegaPi$ is the inverse of $\langle -k_{0} ^2  \ROmegaPi \rangle $, so $\bfuDA = \langle -k_{0} ^2  \ROmegaPi \rangle \bfv^{0}$. Moreover, Proposition~\ref{operateurs R} states that $\Delta \circ \ROmegaPi= \mathrm{id}$, $\partial _{Y|\couche} ^2  \circ \ROmegaPi =0$ and $\partial _{Y|\Gamma ,Y=0^{-}} \circ \ROmegaPi=0$, which imply resp. $(\mu _{0} \Delta  +\omega ^2 \rho _{0} )\circ \langle -k_{0} ^2  \ROmegaPi \rangle  = \mu _{0} \Delta  \circ ( \mathrm{id} +k_{0} ^2  \ROmegaPi)\circ \langle -k_{0} ^2  \ROmegaPi \rangle =\mu _{0} \Delta $,\; $\partial _{Y|\couche} ^2  \circ \langle -k_{0} ^2  \ROmegaPi \rangle  = \partial _{Y|\couche} ^2 $ and $\partial _{Y|\Gamma ,Y=0^{-}} \circ \langle -k_{0} ^2  \ROmegaPi \rangle  = \partial _{Y|\Gamma ,Y=0^{-}} ^2 $. Therefore, \eqref{eq: equations de bfuDA} is equivalent to
\begin{equation}{\label{eq: DA algebrique des champs lointains concatenes preuve}}
\left\{\begin{array}{r@{\;=\;}ll}
	\mu _{0} \Delta  \bfv^{0} & 0&\text{ in }\Omega  \\[0.5mm]
	\mu _{1} \partial _{Y} ^2  \bfv^{0} & g^{0} &\text{ in }\couche \\[0.5mm]
	\mu _{1} \partial _{Y} \bfv^{0} _{|Y=0^{-}} & h^{0} &\text{ on }\Gamma  
\end{array}\right.\end{equation}
Then Propositions~\ref{pb de Laplace homogene dans A} and \ref{operateurs R} imply that the solutions of \eqref{eq: DA algebrique des champs lointains concatenes preuve} in $\Amacro(\Pi )$ are the $ \frac{1}{\mu _{1}} \RcouchePi(g^{0} ) + \frac{1}{\mu _{1}} \RGammaPi(h^{0} ) + \fsum_{d\in \TZ*} c_{d} \,\phi_{d}$ where $c_{d}$ vanishes for small enough $d$. Finally we have $c_{d} =\sigma _{d} (\bfv^{0} ) = \sigma _{d} (\bfuDA)$ for any $d\in \TZ*$ because $\sigma _{d} \circ \ROmegaPi$, $\sigma _{d} \circ \RcouchePi$ and $\sigma _{d} \circ \RGammaPi$ vanish by Proposition~\ref{operateurs R}.
\end{demo}

\begin{lemme}{}
\label{DA algebrique des champs de coin}%
Let $F^{\infty } _{\Omega } \in \Amicro(\Omega )$ and $F^{\infty } _{\couche} \in \Amicro(\couche)$. The solutions in $\Amicro(\Pi )$ of
\begin{equation}{\label{eq: equations de SDA}}
\left\{\begin{array}{r@{\;=\;}ll}
\mu _{0} \Delta  \SDA & F^{\infty } _{\Omega } &\text{ in }\Omega  \\[0.5mm]
\mu _{1} \Delta  \SDA & F^{\infty } _{\couche} &\text{ in }\couche \\[0.5mm]
\mu _{0} \partial _{Y|Y=0^{+}} \SDA - \mu _{1} \partial _{Y|Y=0^{-}} \SDA & 0 &\text{ on }\Gamma 
\end{array}\right.\end{equation}
are the formal series of the following form (where $\sigma _{d} (\SDA)$ vanishes when $d$ is big enough)
\begin{equation}{\label{eq: SDA}}
\SDA = \left\langle -\RcouchePi\circ \partial _{X|\couche} ^2 ,\; \frac{\mu _{0} }{\mu _{1}} \RGammaPi\circ \partial _{Y|\Gamma ,Y=0^{+}} \right\rangle  \bigg( \frac{1}{\mu _{0}} \ROmegaPi (F^{\infty } _{\Omega } )+ \frac{1}{\mu _{1}} \RcouchePi (F^{\infty } _{\couche} ) + \!\fsum_{d\in \TZ*} \! \sigma _{d} (\SDA) \,\phi_{d} \bigg) \end{equation}\fbeq[1.3]
\end{lemme}

\begin{demo}
This is similar to Lemma~\ref{DA algebrique des champs lointains concatenes}. Let $R_{1} :=-\RcouchePi\circ \partial _{X|\couche} ^2 $, $R_{2} :=\frac{\mu _{0} }{\mu _{1}} \RGammaPi\circ \partial _{Y|\Gamma ,Y=0^{+}} $ and $\bfv^{\infty } := ( \mathrm{id}-R_{1} -R_{2} )\SDA$. By Definition~\ref{def: serie geometrique d'operateurs 1D} we have
\[\langle R_{1} ,R_{2} \rangle  = \sum_{n=0}^\infty  \sum _{(i_{1} ,\dots ,i_{n} )\in \{1,2\}^{n} } R_{i_{1}} \circ \cdots \circ R_{i_{n}} = \sum_{n=0}^\infty  (R_{1} +R_{2} )^{n} = ( \mathrm{id}-R_{1} -R_{2} )^{-1} .\]
So $\SDA = \langle R_{1} ,R_{2} \rangle \bfv^{\infty }$. Moreover, Proposition~\ref{operateurs R} implies that $\Delta \circ R_{1} =\Delta \circ R_{2} =0$, $\partial _{Y|\couche} ^2  \circ R_{1} =-\partial _{X|\couche} ^2 $, $\partial _{Y|\couche} ^2  \circ R_{2} =0$, $\partial _{Y|\Gamma ,Y=0^{-}} \circ R_{1} =0$ and $\partial _{Y|\Gamma ,Y=0^{-}} \circ R_{2} =\frac{\mu _{0} }{\mu _{1}} \partial _{Y|\Gamma ,Y=0^{+}}$. We deduce that $\Delta \circ \langle R_{1} ,R_{2} \rangle =\Delta $,
\[(\partial _{Y|\couche} ^2  + \partial _{X|\couche} ^2 )\circ \langle R_{1} ,R_{2} \rangle  = \partial _{Y|\couche} ^2  \circ ( \mathrm{id}-R_{1} )\circ \langle R_{1} ,R_{2} \rangle  = \partial _{Y|\couche} ^2  \circ ( \mathrm{id}-R_{1} -R_{2} )\circ \langle R_{1} ,R_{2} \rangle  =\partial _{Y|\couche} ^2 \] and similarly $(\mu _{0} \partial _{Y|\Gamma ,Y=0^{+}} - \mu _{1} \partial _{Y|\Gamma ,Y=0^{-}} )\circ \langle R_{1} ,R_{2} \rangle =\mu _{1} \partial _{Y|\Gamma ,Y=0^{-}}$. Therefore, \eqref{eq: equations de SDA} is equivalent to
\begin{equation}{\label{eq: DA algebrique des champs de coin preuve}}
\left\{\begin{array}{r@{\;=\;}ll}
	\mu _{0} \Delta  \bfv^{\infty } & F^{\infty } _{\Omega } &\text{ in }\Omega  \\[0.5mm]
	\mu _{1} \partial _{Y} ^2  \bfv^{\infty } & F^{\infty } _{\couche} &\text{ in }\couche \\[0.5mm]
	\mu _{1} \partial _{Y} \bfv^{\infty } _{|Y=0^{-}} & 0 &\text{ on }\Gamma  
\end{array}\right.\end{equation}
Then Propositions~\ref{pb de Laplace homogene dans A} and \ref{operateurs R} imply that the solutions of \eqref{eq: DA algebrique des champs de coin preuve} in $\Amacro(\Pi )$ are the $ \frac{1}{\mu _{0}} \ROmegaPi (F^{\infty } _{\Omega } )+ \frac{1}{\mu _{1}} \RcouchePi (F^{\infty } _{\couche} ) + \fsum_{d\in \TZ*} c_{d} \,\phi_{d}$ where $c_{d} =0$ for big enough $d$. Finally, $c_{d} = \sigma _{d} (\SDA)$ as for Lemma~\ref{DA algebrique des champs lointains concatenes}.
\end{demo}

\begin{paragraphesansboite}{Definition \nvnumpar{}:}
\label{def: petit o derivable}%
Let $d\in \mathbb{R}$ and $a\in \{0,\infty \}$. We define $o_{\partial }$, a kind of differentiable small $o$, as follows.
\begin{itemize}
	\item For any $\varphi :\Omega  \rightarrow \mathbb{R}$, we say that $\varphi  = o_{\partial } (r^{d} )$ when $r\rightarrow a$ if $\varphi $ is $\mathcal{C}^{\infty }$ in a vicinity of $r=a$ and $\forall (j,k)\in \mathbb{N}^2 ,\; \partial _{r} ^{j} \partial _{\theta } ^{k} \varphi  = o(r^{d-j} )$ uniformly w.r.t. $\theta $ when $r\rightarrow a$.
	\item For any $\varphi :\couche \rightarrow \mathbb{R}$, we say that $\varphi  = o_{\partial } (x^{d} )$ when $x\rightarrow a$ if $\varphi $ is $\mathcal{C}^{\infty }$ in a vicinity of $x=a$ and $\forall (j,k)\in \mathbb{N}^2 ,\; \partial _{r} ^{j} \partial _{Y} ^{k} \varphi  = o(x^{d-j} )$ uniformly w.r.t. $Y$ when $x\rightarrow a$.
	\item For any $\varphi :\Gamma \rightarrow \mathbb{R}$, we say that $\varphi  = o_{\partial } (x^{d} )$ when $x\rightarrow a$ if $\varphi $ is $\mathcal{C}^{\infty }$ in a vicinity of $x=a$ and $\forall j\in \mathbb{N},\; \partial _{r} ^{j} \varphi  = o(x^{d-j} )$ when $x\rightarrow a$.
	\item For any $\varphi :\Pi  \rightarrow \mathbb{R}$, we say that $\varphi  = o_{\partial } (\bfr^{d} )$ when $\bfr\rightarrow a$ if $\varphi _{|\Omega } = o_{\partial } (r^{d} )$ and $\varphi _{|\couche} = o_{\partial } (x^{d} )$.
\end{itemize}\fbi
\end{paragraphesansboite}

\begin{paragraphesansboite}{Definition \nvnumpar{}:}
\label{def: troncature dans Abarre}%
Let $d\in \mathbb{R}$, $D\in \{\Pi ,\Omega ,\Gamma ,\couche\}$ and $\varphi  = \fsum_{q\in \mathbb{R}} \varphi _{q}$ in $\Amacro(D)$ or $\Amicro(D)$ with $\varphi _{q} \in \mathcal{A}_{q} (D)$ for all $q\in \mathbb{R}$. We denote $T_{\leqslant d} (\varphi ):=\fsum_{q\in \mathbb{R},q\leqslant d} \varphi _{q}$ and $T_{\geqslant d} (\varphi ):=\fsum_{q\in \mathbb{R},q\geqslant d} \varphi _{q}$ the truncations of $\varphi $ below and above $d$.
\end{paragraphesansboite}

\begin{textesansboite}
Using Lemmas~\ref{DA algebrique des champs lointains concatenes}--\ref{DA algebrique des champs de coin} and Kondrat'ev's theory (involving weighted Sobolev spaces, Laplace's transform and the residue theorem), we proved the following theorems, giving asymptotic behaviors for solutions of model problems of the type of far-and-layer fields and corner fields. The proofs can be found in Appendix~\ref{sec: DA autour du coin}.
\end{textesansboite}

\begin{boitecarreesecable}{Theorem \nvnumpar{}:}
\label{DA des champs lointains concatenes}%
Let $\bfu\in H^{1} _{\mathrm{loc}} (\Pi )$, $f\in \mathcal{D}'(\Omega )$, $g\in L^2 _{\mathrm{loc}} (\couche)$ and $h\in L^2  _{\mathrm{loc}} (\Gamma )$ be such that
\begin{equation}{\label{eq: DA des champs lointains concatenes, hyp}}
\left\{\begin{array}{r@{\;=\;}ll}
	\mu _{0} \Delta  \bfu +\omega ^2 \rho _{0} \bfu & f &\text{ in }\Omega  \\[0.5mm]
	\mu _{1} \partial _{Y} ^2  \bfu & g &\text{ in }\couche \\[0.5mm]
	\mu _{1} \partial _{Y|Y=0^{-}} \bfu & h &\text{ on }\Gamma  \\[0.5mm]
	\bfu_{|y=0^{+}} - \bfu_{|Y=0^{-}} & 0 &\text{ on }\Gamma \\[0.5mm]
	\bfu & 0 &\text{ on }\bordgauche \cup \borddroit
\end{array}\right.\end{equation}
We assume that:
\begin{itemize}
	\item $f$ vanishes in the vicinity of the corner $(0,0)$,
	\item there exists $g^{0} \in \Amacro(\couche)$ s.t.: $\forall d\in \mathbb{R},\; g(x,Y) = T_{\leqslant d} (g^{0} )(x,Y) + o_{\partial } (x^{d} )$ when $x\rightarrow 0$,
	\item there exists $h^{0} \in \Amacro(\Gamma )$ s.t.: $\forall d\in \mathbb{R},\; h(x) = T_{\leqslant d} (h^{0} )(x) + o_{\partial } (x^{d} )$ when $x\rightarrow 0$,
	\item there is $\eta >0$, $u^{\mathrm{v}} \in H^{1} (\Omega \cap B(0,\eta ))$ and $\varphi \in \mathcal{A}(\Omega )$ s.t. $\bfu_{|\Omega \cap B(0,\eta )} =u^{\mathrm{v}} + \varphi $.
\end{itemize}
Then there is $\bfuDA \in \Amacro(\Pi )$ that has the form \eqref{eq: bfuDA} s.t.: \ $\forall d\in \mathbb{R},\ \bfu = T_{\leqslant d} (\bfuDA ) + o_{\partial } (\bfr^{d} )$ when $\bfr\rightarrow 0$.
\end{boitecarreesecable}

\begin{boitecarreesecable}{Theorem \nvnumpar{}:}
\label{DA des champs de coin}%
Let $S \in H^{1} _{\mathrm{loc}} (\Omega _{1} )$ and $F\in L^2 _{\mathrm{loc}} (\overline{\Omega _{1} })$ (i.e. $L^2 $ on any bounded subset of $\Omega _{1}$) such that 
\begin{equation}{\label{eq: DA des champs de coin, hyp}}
\left\{\begin{array}{r@{\;=\;}ll}
	\operatorname{div}(\mu  \nabla  S) & F &\text{ in }\Omega _{1} \\[0.5mm]
	S & 0 &\text{ on }\partial \Omega _{1}
\end{array}\right.\end{equation}
We assume that:
\begin{itemize}
	\item there exists $F^{\infty } _{\Omega } \in \Amicro(\Omega )$ such that $\forall d\in \mathbb{R},\; F_{|\Omega } = T_{\geqslant d} (F^{\infty } _{\Omega } ) + o_{\partial } (r^{d} )$ when $r\rightarrow \infty $,
	\item there exists $F^{\infty } _{\couche} \in \Amicro(\couche)$ such that $\forall d\in \mathbb{R},\; F_{|\couche} = T_{\geqslant d} (F^{\infty } _{\couche} ) + o_{\partial } (x^{d} )$ when $x\rightarrow \infty $,
	\item $S$ belongs to $\Vcoin+\khii \mathcal{A}(\Pi )$ (the space in which the corner fields will be build in Section~\ref{subsec: cadre fonctionnel des champs de coin}).
\end{itemize}
Then there is $\SDA \in \Amicro(\Pi )$ that has the form \eqref{eq: SDA} s.t.: \ $\forall d\in \mathbb{R},\ S = T_{\geqslant d} (\SDA ) +	o_{\partial } (\bfr^{d} )$ when $\bfr\rightarrow \infty $.
\end{boitecarreesecable}

\begin{textesansboite}
Note that, since $g^{0} \in \Amacro(\couche)$, we have $T_{\leqslant d} (g^{0} )\in \mathcal{A}(\couche)$, so the formula $g = T_{\leqslant d} (g^{0} ) +o_{\partial } (x^{d} )$ makes sense. The same applies to the truncations of $h^{0}$, $\bfu^{0}$, $F^{\infty } _{\Omega }$ and $F^{\infty } _{\couche}$ and $\SDA$.\\

A consequence of Theorems~\ref{DA des champs lointains concatenes}--\ref{DA des champs de coin} is Proposition~\ref{confimation des DA des champs} that states that for any $(p,\ell )$ there is $\bfuDA_{p,\ell } \in \Amacro(\Pi )$ and $\SDA_{p,\ell } \in \Amicro(\Pi )$ s.t.:\vspace{-3mm}
\[\forall d\in \mathbb{R}, \quad  \left\{\begin{array}{r@{\;=\;}ll} 
	\bfu_{p,\ell } & T_{\leqslant d} (\bfuDA_{p,\ell } ) + o_{\partial } (\bfr^{d} ) & \mathrm{when} \ \bfr\rightarrow 0\\[0.5mm]
	S_{p,\ell } & T_{\geqslant d} (\SDA_{p,\ell } ) + o_{\partial } (\bfr^{d} )& \mathrm{when} \ \bfr\rightarrow \infty 
\end{array} \right.\]
In the rest of this section, we will assume that such formal series exist. In addition, given the equations satisfied by $\bfu_{p,\ell }$ and $S_{p,\ell }$ (see \eqref{eq: eq vol de bfu p,l} and \eqref{eq: eq vol champs de coin}), \eqref{eq: bfuDA} and \eqref{eq: SDA} rewrite here as
\begin{align}
\label{eq: bfuDA p,l}
\!\bfuDA_{p,\ell } &= \big\langle  {-}k_{0} ^2 \,\ROmegaPi \big\rangle  \bigg( -\RcouchePi\circ (\partial _{x|\couche} ^2 +k_{1} ^2 )(\bfuDA_{p-2,\ell } ) + \frac{\mu _{0} }{\mu _{1}} \RGammaPi\circ \partial _{y|\Gamma ,y=0^{+}} (\bfuDA_{p-1,\ell } ) + \!\fsum_{d\in \TZ*} \! \sigma _{d} (\bfuDA_{p,\ell } ) \,\phi_{d} \bigg) \!\\
\nonumber
\!\SDA_{p,\ell } &= \left\langle \! -\RcouchePi\circ \partial _{X|\couche} ^2 , \frac{\mu _{0} }{\mu _{1}} \RGammaPi\circ \partial _{Y|\Gamma ,Y=0^{+}} \!\right\rangle  \!\bigg( \!- k_{0} ^2  \ROmegaPi (\SDA_{p-2,\ell |\Omega } )- k_{1} ^2  \RcouchePi (\SDA_{p-2,\ell |\couche} ) + \!\!\fsum_{d\in \TZ*} \! \sigma _{d} (\SDA_{p,\ell } ) \,\phi_{d} \bigg) 
\end{align}
where by convention $\bfuDA_{p,\ell }$ and $\SDA_{p,\ell }$ vanish when $p\in \mathbb{R} \setminus \NN$. Therefore $(\bfuDA_{p,\ell } )_{p,\ell }$ and $(\SDA_{p,\ell } )_{p,\ell }$ are uniquely defined by $(\sigma _{d} (\bfuDA_{p,\ell } ))_{d,p,\ell }$ and $(\sigma _{d} (\SDA_{p,\ell } ))_{d,p,\ell }$. When $d<0$ (resp. $d>0$), $\phi_{d}$ is non-variational for the far-and-layer fields (resp. corner fields), and Theorems~\ref{cadre fonctionnel des champs lointains} and \ref{cadre fonctionnel pour les champs de coin} show that $\sigma _{d} (\bfuDA_{p,\ell } )$ (resp. $\sigma _{d} (\SDA_{p,\ell } )$) can be fixed arbitrarily. The rest of this section is devoted to finding how to fix them in order to satisfy the matching conditions~\eqref{eq: hypothese de raccord}. On the contrary, when $d>0$ (resp. $d<0$), $\phi_{d}$ is variational for the far-and-layer fields (resp. corner fields) and the values of $\sigma _{d} (\bfuDA_{p,\ell } )$ (resp. $\sigma _{d} (\SDA_{p,\ell } )$) are uniquely defined once $(\sigma _{d'} (\bfuDA_{p,\ell } ))_{d'<0}$ (resp. $(\sigma _{d'} (\SDA_{p,\ell } ))_{d'>0}$) has been fixed. Ways to compute these values numerically will be investigated in a future work.
\end{textesansboite}

\subsection{Specifying of the matching conditions}
\label{subsec: raccordement algebrique}

\begin{textesansboite}
In this section we express $\fsum_{p,\ell } \varepsilon ^{p} \ln ^{\ell } \!\varepsilon  \cdot \bfuDA_{p,\ell }$ and $\fsum_{p,\ell } \varepsilon ^{p} \ln ^{\ell } \! \varepsilon \cdot \SDA_{p,\ell }$ in function of the $\sigma _{d} (\bfuDA_{p,\ell } )$ and $\sigma _{d} (\SDA_{p,\ell } )$, we then rewrite rigorously the matching conditions~\eqref{eq: hypothese de raccord}, and we finally show that they are equivalent to a set of equations on the coefficients $\sigma _{d} (.)$. Here $\fsum_{p,\ell } \varepsilon ^{p} \ln ^{\ell } \!\varepsilon  \cdot \bfuDA_{p,\ell }$ and $\fsum_{p,\ell } \varepsilon ^{p} \ln ^{\ell } \! \varepsilon \cdot \SDA_{p,\ell }$ are formal series that belong resp. to the spaces $\Aemacro$ and $\Aemicro$ defined below. \\

In this section, ``$\varepsilon $'' and ``$\ln \varepsilon $'' denote two algebraic indeterminates independent of each other (so they are not numbers).
\end{textesansboite}

\begin{definition}{}
We denote \[\Aemacro := \fsum_{p\in \NN} \,\varepsilon ^{p} \!\fsum_{d\in \NN-p} \mathcal{A}_{d} (\Pi )[\ln \varepsilon ] \qquad  \text{ and } \qquad  \Aemicro := \fsum_{p\in \NN} \,\varepsilon ^{p} \!\fsum_{d\in p-\NN} \mathcal{A}_{d} (\Pi )[\ln \varepsilon ].\]
According to Section~\ref{subsec: preliminaires algebriques}, we write their elements as formal series like
\[ \fsum_{p\in \NN} \;\fsum_{d\in \NN-p} \;\fsum_{\ell \in \mathbb{N}} \varepsilon ^{p} \ln ^{\ell } \!\varepsilon \cdot \varphi _{p,d,\ell ,} \qquad  \text{resp.}  \qquad  \fsum_{p\in \NN} \;\fsum_{d\in p-\NN} \;\fsum_{\ell \in \mathbb{N}} \varepsilon ^{p} \ln ^{\ell } \!\varepsilon \cdot \varphi _{p,d,\ell } .\]\fbeq[1.3]
\end{definition}

\begin{figure}[h]
\begin{center}
\includegraphics{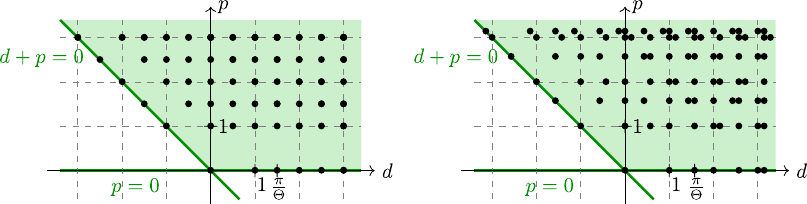}
\caption{Points of $\{(d,p) \mid p\in \NN \text{ and } d\in \NN-p\}$ for $\Theta = \frac{2\pi }{3}$ (on the left) and $\Theta =2$ (on the right)}
\label{fig:cone de NN}
\end{center}
\end{figure}

\begin{paragraphesansboite}{Ansatz \nvnumpar{}:}
\label{ansatz sur Aeps}%
We assume that \[\fsum_{p\in \NN,\ell \in \mathbb{N}} \varepsilon ^{p} \ln ^{\ell } \!\varepsilon  \cdot \bfuDA_{p,\ell } \in \Aemacro \qquad  \text{ and } \qquad  \fsum_{p\in \NN,\ell \in \mathbb{N}} \varepsilon ^{p} \ln ^{\ell } \!\varepsilon \cdot \SDA_{p,\ell } \in \Aemicro. \] Compared to the ansatz of Section~\ref{subsec: presentation du dvlpmt raccorde}, it adds that, for any $(p,\ell )\in \NN\times \mathbb{N}$, $\bfu_{p,\ell }$ has components only in the $\mathcal{A}_{d} (\Pi )$ s.t. $d\in \NN-p$ and $S_{p,\ell }$ has only in the $\mathcal{A}_{d} (\Pi )$ s.t. $d\in p-\NN$. This is necessary for the matching because e.g. we will see that for any $\varphi \in \varepsilon ^{p} \mathcal{A}_{d} (\Pi )[\ln \varepsilon ]\subset \Aemicro$ that is a term of $\fsum \varepsilon ^{p} \ln ^{\ell } \!\varepsilon  \,\SDA_{p,\ell }$, $\varphi ( \frac{\cdot }{\varepsilon } )$ must appear in $\fsum \varepsilon ^{p} \ln ^{\ell } \!\varepsilon  \,\bfuDA_{p,\ell }$ and $\varphi ( \frac{\cdot }{\varepsilon } ) \in \varepsilon ^{p-d} \mathcal{A}_{d} (\Pi )[\ln \varepsilon ]$, so $p-d\in \NN$.
\end{paragraphesansboite}

\begin{figure}[h]
\def\xmax{2.4}
\def\ymin{-0.3}
\def\ymax{2.3}
\begin{center}
\begin{tikzpicture}[scale=1]
	\draw [thin, gray, dashed] (-\xmax,\ymin) grid (\xmax,\ymax);
	\draw [->] (-\xmax,0) -- (\xmax,0) node [right] {$\deg_{\mathcal{A}}$};
	\draw [->] (0,\ymin) -- (0,\ymax+0.2) node [right] {$\deg_\varepsilon$};
	\draw (1,0) node{$+$} node[below]{\ \ 1};
	\draw (0,1) node{$+$} node[right]{1};
	\draw [->,>=latex,very thick,blue] (0,0) --  node[midway, above] {$-k_{0} ^2 \,\ROmegaPi$} (2,0);
	\draw [->,>=latex,very thick,magenta] (-0.9,0.9) -- (-2,2) node[above] {$-\varepsilon ^2  \RcouchePi\circ \partial _{x|\couche} ^2 $};
	\draw [->,>=latex,very thick,red] (0,0) -- node[near end,right] {$-\varepsilon ^2 k_{1} ^2 \RcouchePi$} (0,2);
	\draw [->,>=latex,very thick,green!55!black] (0,0) -- node[midway,left] {$\varepsilon  \frac{\mu _{0} }{\mu _{1}} \RGammaPi\circ \partial _{y|\Gamma ,y=0^{+}}$} (-1,1);
\end{tikzpicture}%
\hspace{10mm}\nolinebreak%
\begin{tikzpicture}[scale=1]
	\draw [thin, gray, dashed] (-\xmax,\ymin) grid (\xmax,\ymax);
	\draw [->] (-\xmax,0) -- (\xmax,0) node [right] {$\deg_{\mathcal{A}}$};
	\draw [->] (0,\ymin) -- (0,\ymax+0.2) node [right] {$\deg_\varepsilon$};
	\draw (1,0) node{$+$} node[above]{\ \ 1};
	\draw (0,1) node{$+$} node[right]{1};
	\draw [->,>=latex,very thick,blue] (0,0) -- node[midway, right] {\vspace{-4pt}$-\varepsilon ^2  k_{0} ^2 \,\ROmegaPi$} (2,2);
	\draw [->,>=latex,very thick,magenta] (-0.9,0) -- (-2,0) node[above] {$-\RcouchePi\circ \partial _{x|\couche} ^2 $};
	\draw [->,>=latex,very thick,red] (0,0) -- (0,2) node[below left] {$-\varepsilon ^2  k_{1} ^2 \RcouchePi$};
	\draw [->,>=latex,very thick,green!55!black] (0,0) node[below left] {$\frac{\mu _{0} }{\mu _{1}} \RGammaPi\circ \partial _{y|\Gamma ,y=0^{+}} \!\!$} -- (-1,0);
\end{tikzpicture}\par
\caption{Degrees of the elements of $\Rmacro$ (on the left) and $\Rmicro$ (on the right)}
\label{fig: degres des elements de Rmacro et Rmicro}
\end{center}
\end{figure}
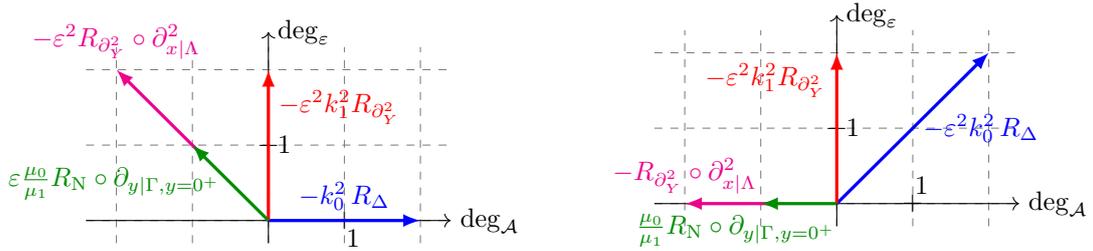

\begin{textesansboite}
Let us denote: \vspace{-3mm}
\[\textstyle \Rmacro := \left\{-k_{0} ^2 \,\ROmegaPi,\ -\varepsilon ^2  \RcouchePi\circ \partial _{x|\couche} ^2 ,\ -\varepsilon ^2 k_{1} ^2 \RcouchePi,\ \varepsilon  \frac{\mu _{0} }{\mu _{1}} \RGammaPi\circ \partial _{y|\Gamma ,y=0^{+}} \right\}\]
\[\textstyle \Rmicro := \left\{-\varepsilon ^2  k_{0} ^2 \,\ROmegaPi,\ -\RcouchePi\circ \partial _{X|\couche} ^2 ,\ -\varepsilon ^2  k_{1} ^2 \RcouchePi,\ \frac{\mu _{0} }{\mu _{1}} \RGammaPi\circ \partial _{Y|\Gamma ,Y=0^{+}} \right\}\]
where the above operators are defined as in \eqref{eq: operateur a degre selon p et d}. For any $R\in \Rmacro \cup \Rmicro$, we denote $(\deg _{\varepsilon } R,\deg _{\mathcal{A}} R ):=\deg  R$. See Figure~\ref{fig: degres des elements de Rmacro et Rmicro}. Thanks to Definition~\ref{def: serie geometrique d'operateurs 2D}, we can consider $\langle \Rmaicro\rangle $ which is well-defined on $\Aemaicro$ (it suffices to take $v:=(\pm 1,2)$ in Definition~\ref{def: serie geometrique d'operateurs 2D}). Moreover one can check that it maps $\Aemaicro$ into itself.
\end{textesansboite}

\begin{theoreme}{}
\label{champs totaux en fct de sigma}%
We have the following equalities in $\Aemacro$ and $\Aemicro$ respectively: \begin{align}
	\fsum_{p\in \NN,\ell \in \mathbb{N}} \varepsilon ^{p} \ln ^{\ell } \!\varepsilon  \cdot \bfuDA_{p,\ell } &= \langle \Rmacro\rangle \bigg(\fsum_{p\in \NN,\ell \in \mathbb{N}} \varepsilon ^{p} \ln ^{\ell } \!\varepsilon \, \sumdmacro \sigma _{d} (\bfuDA_{p,\ell } )\,\phi_{d} \bigg)
	\label{eq: champ lointain total}\\
	\fsum_{p\in \NN,\ell \in \mathbb{N}} \varepsilon ^{p} \ln ^{\ell } \!\varepsilon \cdot \SDA_{p,\ell } &= \langle \Rmicro\rangle \bigg(\fsum_{p\in \NN,\ell \in \mathbb{N}} \varepsilon ^{p} \ln ^{\ell } \!\varepsilon \, \sumdmicro \sigma _{d} (\SDA_{p,\ell } )\,\phi_{d} \bigg)
	\label{eq: champ de coin total}
\end{align}\fbeq
\end{theoreme}

\begin{demo}
We only prove \eqref{eq: champ lointain total}, as \eqref{eq: champ de coin total} is similar. We could do it by inductively by composing \eqref{eq: bfuDA p,l}, but we chose instead a proof similar to Lemma~\ref{DA algebrique des champs de coin} to avoid heavy calculations. Let  $(p,\ell )\in \NN\times \mathbb{N}$. Given the equations satisfied by $\bfu_{p,\ell }$ in \eqref{eq: eq vol de bfu p,l}, and the fact that $\bfu_{p,\ell } = T_{\leqslant d} (\bfuDA_{p,\ell } ) + o_{\partial } (\bfr^{d} )$ for any $d\in \mathbb{R}$, we have
\begin{equation}{\label{eq: champs totaux en fct de sigma, preuve}}
\left\{\begin{array}{r@{\;=\;}ll}
	\mu _{0} \Delta  \bfuDA_{p,\ell } +\omega ^2 \rho _{0} \bfuDA_{p,\ell } & 0 &\text{ in }\Omega  \\[0.5mm]
	\mu _{1} \partial _{Y} ^2  \bfuDA_{p,\ell } & -(\mu _{1} \partial _{x} ^2 +\omega ^2 \rho _{1} )\bfuDA_{p-2,\ell } &\text{ in }\couche \\[0.5mm]
	\mu _{1} \partial _{Y|Y=0^{-}} \bfuDA_{p,\ell } & \mu _{0} \,\partial _{y|y=0^{+}} \bfuDA_{p-1,\ell } &\text{ on }\Gamma 
\end{array}\right.\end{equation}
Let $\bfuDA_{\varepsilon } :=\fsum_{p,\ell } \varepsilon ^{p} \ln ^{\ell } \!\varepsilon \,\bfuDA_{p,\ell }$. Summing over $(p,\ell )$ \eqref{eq: champs totaux en fct de sigma, preuve} times $\varepsilon ^{p} \ln ^{\ell } \!\varepsilon $ yields
\begin{equation}{\label{eq: probleme de Helmholtz du champs total algebrique}}
\left\{\begin{array}{r@{\;=\;}ll}
	\mu _{0} \Delta  \bfuDA_{\varepsilon } +\omega ^2 \rho _{0} \bfuDA_{\varepsilon } & 0 &\text{ in }\Omega  \\[0.5mm]
	\mu _{1} \partial _{Y} ^2  \bfuDA_{\varepsilon } + \varepsilon ^2 (\mu _{1} \partial _{x} ^2 +\omega ^2 \rho _{1} )\bfuDA_{\varepsilon } &0&\text{ in }\couche \\[0.5mm]
	\mu _{1} \partial _{Y|Y=0^{-}} \bfuDA_{\varepsilon } - \varepsilon \,\mu _{0} \,\partial _{y|y=0^{+}} \bfuDA_{\varepsilon } &0 &\text{ on }\Gamma 
\end{array}\right.\end{equation}
Let $\bfv^{0} _{\varepsilon } := ( \mathrm{id}-\sum _{R\in \Rmacro} R)\bfuDA_{\varepsilon }$. Let us show that \eqref{eq: probleme de Helmholtz du champs total algebrique} implies \eqref{eq: champ lointain total}. It is easy to see that $ \mathrm{id}-\sum _{R\in \Rmacro} R$ maps $\Aemacro$ into itself, so $\bfv^{0} _{\varepsilon } \in \Aemacro$. By Definition~\ref{def: serie geometrique d'operateurs 2D}, we have
\[\langle \Rmacro \rangle  := \sum_{n=0}^\infty  \sum _{(R_{1} ,\dots ,R_{n} )\in \Rmacron} R_{1} \circ \cdots \circ R_{n} = \sum_{n=0}^\infty  \bigg( \sum _{R\in \Rmacro} R \bigg)^{n} = \bigg( \mathrm{id}-\sum _{R\in \Rmacro} R \bigg)^{-1} .\]
So $\bfuDA_{\varepsilon } = \langle \Rmacro\rangle \bfv^{0} _{\varepsilon }$. Moreover, Proposition~\ref{operateurs R} implies $\forall R\in \Rmacro \setminus \{-k_{0} ^2 \ROmegaPi\},\; \Delta \circ R=0$. Hence
\[(\mu _{0} \Delta +\omega ^2 \rho _{0} )\circ \langle \Rmacro\rangle  = \mu _{0} \Delta  \circ ( \mathrm{id}+k_{0} ^2 \ROmegaPi)\circ \langle \Rmacro\rangle  = \mu _{0} \Delta \circ \bigg( \mathrm{id}-\sum _{R\in \Rmacro} R \bigg)\circ \langle \Rmacro\rangle  = \mu _{0} \Delta \]
Similarly, one can check that $\big[\mu _{1} \partial _{Y|\couche} ^2  + \varepsilon ^2 (\mu _{1} \partial _{x|\couche} ^2 +\omega ^2 \rho _{1} )\big]\circ \langle \Rmacro\rangle  = \mu _{1} \partial _{Y|\couche} ^2  $ and $(\mu _{1} \partial _{Y|\Gamma ,Y=0^{-}} -\varepsilon \,\mu _{0} \,\partial _{y|\Gamma ,y=0^{+}} )\circ \langle \Rmacro\rangle = \mu _{1} \partial _{Y|\Gamma ,Y=0^{-}}$. Therefore $\bfv^{0} _{\varepsilon }$ satisfies
\[\left\{\begin{array}{r@{\;=\;}ll}
	\mu _{0} \Delta  \bfv^{0} _{\varepsilon } & 0 &\text{ in }\Omega  \\[0.5mm]
	\mu _{1} \partial _{Y} ^2  \bfv^{0} _{\varepsilon }  &0&\text{ in }\couche \\[0.5mm]
	\mu _{1} \partial _{Y|Y=0^{-}} \bfv^{0} _{\varepsilon } &0 &\text{ on }\Gamma 
\end{array}\right.\]
Then Proposition~\ref{pb de Laplace homogene dans A} implies that there are numbers $c_{p,\ell ,d}$ s.t. $\bfv^{0} _{\varepsilon } = \fsum_{p,\ell ,d} \varepsilon ^{p} \ln ^{\ell } \!\varepsilon \,c_{p,\ell ,d} \,\phi_{d}$. Finally we have $c_{p,\ell ,d} = \sigma _{d} (\bfuDA_{p,\ell } )$ for any $(p,\ell ,d)$ because Proposition~\ref{operateurs R} gives that $\forall R\in \Rmacro,\; \sigma _{d} \circ R=0$.
\end{demo}

\begin{paragraphesansboite}{Remark \nvnumpar{}:}
\label{rq: utilite des operateurs R}%
The elements of $\mathcal{A}(\Pi )$ present in \eqref{eq: champ lointain total}--\eqref{eq: champ de coin total} are the $R_{1} \circ \cdots \circ R_{n} (\phi_{d} )$ with: $\forall i,\, R_{i} \in \{-k_{0} ^2 \ROmegaPi,\, \RcouchePi\circ \partial _{x|\couche} ^2 ,\, -k_{1} ^2 \RcouchePi,\, \frac{\mu _{0} }{\mu _{1}} \RGammaPi\circ \partial _{y|\Gamma ,y=0^{+}} \}$. They cannot be described using only inductive sequences and their number grows exponentially with $n$. That is why we introduced the operators $\ROmegaPi$, $\RcouchePi$ and $\RGammaPi$. They allow a factorization of \eqref{eq: champ lointain total}--\eqref{eq: champ de coin total} into a geometric series of operators, that encapsulates all the complexity of the singularities, and the coefficents $\sigma _{d} (\dots )$, that need to be fixed for the different fields to match. This will allow us to reduce the matching condition to only fixing the coefficents $\sigma _{d} (\dots )$.
\end{paragraphesansboite}

\begin{textesansboite}
Now, we want to define $\Eps:\Aemacro\rightarrow \Aemicro$ as the scaling operator
\begin{equation}{\label{eq: E}}
\forall \varphi \in \Aemacro, \qquad  \Eps(\varphi ) = \left\{\begin{array}{ll}
	(X,Y) \mapsto  \varphi (\varepsilon X,\varepsilon Y) &\text{ in }\Omega  \\[0.5mm]
	(X,Y) \mapsto  \varphi (\varepsilon X,Y) &\text{ in }\couche
\end{array}\right.\end{equation}
However, $\varepsilon $ is an indeterminate. Let us first define $\Eps$ for any $\varphi \in \mathcal{A}_{d} (\Pi )$, $d\in \mathbb{R}$. If $\varphi _{|\Omega }$ has the form $\Im[(\alpha z)^{q} \,\overline{\alpha z}^{k} P(\log (\alpha z))]$ with $(q,k,P)\in \mathbb{R}\times \mathbb{N}\times \mathbb{R}[T]$ s.t. $q+k=d$ and $\mathcal{P}(q,k,P)$ is true, we set:
	\begin{equation}{\label{eq: E dans Omega}}
	\Eps(\varphi )_{|\Omega } := \sum _{j=0} ^{\deg P} \varepsilon ^{d} \,\frac{\ln ^{j} \!\varepsilon }{j!} \,\Im[(\alpha z)^{q} \,\overline{\alpha z}^{k} P^{(j)} (\log (\alpha z))]\in  \varepsilon ^{d} \mathcal{A}_{d} (\Omega )[\ln \varepsilon ].\end{equation}
And if $\varphi _{|\couche}$ has the form $x^{d} \,Q(\ln (x),Y)$ with $Q\in \mathbb{C}[T]$, we set: 
	\begin{equation}{\label{eq: E dans la couche}}
	\Eps(\varphi )_{|\couche} := \sum _{j=0} ^{\deg Q} \varepsilon ^{d} \,\frac{\ln ^{j} \!\varepsilon }{j!} \,x^{d} \,\partial _{T} ^{j} Q(\ln x,Y) \in  \varepsilon ^{d} \mathcal{A}_{d} (\couche)[\ln \varepsilon ].\end{equation}
Thanks to Lemmas~\ref{decomposition de A en sev} and \ref{unicite de P et Q}, it defines well $\Eps$ from $\mathcal{A}_{d} (\Pi )$ to $\varepsilon ^{d} \,\mathcal{A}_{d} (\Pi )[\ln \varepsilon ]$. Then we extend $\Eps$ to $\Aemacro$ by setting $\Eps(\fsum_{p,d,\ell } \varepsilon ^{p} \ln ^{\ell } \!\varepsilon \, \varphi _{p,d,\ell } ):=\fsum_{p,d} \sum _{\ell } \varepsilon ^{p} \ln ^{\ell } \!\varepsilon \,\Eps(\varphi _{p,d,\ell } )$ for any $(\varphi _{p,d,\ell } )$ s.t. $\varphi _{p,d,\ell } \in \mathcal{A}_{d} (\Pi )$ for any $(p,d,\ell )$.  One can check that is in $\Aemicro$.
\end{textesansboite}

\begin{remarques}
\begin{itemize}
	\item In practice, we will only use the informal definition of \eqref{eq: E}, but everything we will do can be checked using \eqref{eq: E dans Omega} and \eqref{eq: E dans la couche}.\vspace{-0.5mm}
	\item $\Eps$ is invertible and $\Eps^{-1}$ is roughly the scaling: $\Eps^{-1} (\varphi ) = \left\{\begin{array}{l@{\,}l}
		(x,y) \mapsto  \varphi ( \frac{x}{\varepsilon } , \frac{y}{\varepsilon } ) &\text{ in }\Omega  \\[0.5mm]
		(x,Y) \mapsto  \varphi ( \frac{x}{\varepsilon } ,Y) &\text{ in }\couche
	\end{array}\right.$
	\item For the first time, powers of $\ln \varepsilon $ naturally appear because of the power of $\ln r$. This explains why the presence of these powers in the ansatz is necessary from the beginning.
\end{itemize}\fbi
\end{remarques}

\begin{definition}{\vspace{-1.5mm}\titreligne}
\label{def: raccord algebrique}%
We rigorously rewrite the matching condition as: $\displaystyle \!\!\fsum_{p\in \NN,\ell \in \mathbb{N}} \varepsilon ^{p} \ln ^{\ell } \!\varepsilon \cdot \bfuDA_{p,\ell } = \Eps^{-1} \bigg(\fsum_{p\in \NN,\ell \in \mathbb{N}} \varepsilon ^{p} \ln ^{\ell } \!\varepsilon \cdot \SDA_{p,\ell } \bigg)$.\vspace{-1mm}
\end{definition}

\begin{textesansboite}
In order to rewrite the matching condition in terms of the coefficients $\sigma _{d} (\bfuDA_{p,\ell } )$ and $\sigma _{d} (\SDA_{p,\ell } )$, we define the following projectors of $\Aemaicro$:
\[\progphimaicro\bigg(\fsum_{p\in \NN} \fsum_{d\in \pm (\NN-p)} \fsum_{\ell \in \mathbb{N}} \varepsilon ^{p} \ln ^{\ell } \!\varepsilon \cdot  \varphi _{p,d,\ell } \bigg):= \fsum_{p\in \NN} \fsum_{d\in \TZ* \cap  \pm (\NN-p)} \fsum_{\ell \in \mathbb{N}}  \varepsilon ^{p} \ln ^{\ell } \!\varepsilon \cdot  \sigma _{d} (\varphi _{p,d,\ell } )\,\phi_{d}\]
for any $(\varphi _{p,d,\ell } )$ s.t. $\varphi _{p,d,\ell } \in \mathcal{A}_{d} (\Pi )$ for any $(p,d,\ell )$.
\end{textesansboite}

\begin{lemme}{}
\label{raccord algebrique sur les champs}%
The matching condition of Definition~\ref{def: raccord algebrique} is equivalent to each of these equalities:
\begin{align}
	\fsum_{p\in \NN,\ell \in \mathbb{N}} \varepsilon ^{p} \ln ^{\ell } \!\varepsilon  \sumdmacro \sigma _{d} (\bfu_{p,\ell } )\phi_{d} &= \progphimacro \circ  \Eps^{-1} \circ  \langle \Rmicro\rangle  \bigg(\fsum_{p\in \NN,\ell \in \mathbb{N}} \varepsilon ^{p} \ln ^{\ell } \!\varepsilon  \sumdmicro \sigma _{d} (\SDA_{p,\ell } )\phi_{d} \bigg)
	\label{eq: raccord sur les séries de sigma phi du pt de vue des champs lointains}\\
	\fsum_{p\in \NN,\ell \in \mathbb{N}} \varepsilon ^{p} \ln ^{\ell } \!\varepsilon  \sumdmicro \sigma _{d} (S_{p,\ell } )\phi_{d} &= \progphimicro \circ  \Eps \circ  \langle \Rmacro\rangle  \bigg(\fsum_{p\in \NN,\ell \in \mathbb{N}} \varepsilon ^{p} \ln ^{\ell } \!\varepsilon  \sumdmacro \sigma _{d} (\bfuDA_{p,\ell } )\phi_{d} \bigg)
\end{align}\fbeq
\end{lemme}

\begin{demo}
Let us prove only \eqref{eq: raccord sur les séries de sigma phi du pt de vue des champs lointains}. By Theorem~\ref{champs totaux en fct de sigma}, the matching condition of Definition~\ref{def: raccord algebrique} is equivalent to:
\begin{equation}{\label{eq: raccord algebrique presque prouve}}
\langle \Rmacro\rangle \bigg(\fsum_{p,\ell } \varepsilon ^{p} \ln ^{\ell } \!\varepsilon \, \fsum_{d} \sigma _{d} (\bfuDA_{p,\ell } )\phi_{d} \bigg) = \Eps^{-1} \circ  \langle \Rmicro\rangle \circ  \bigg(\fsum_{p,\ell } \varepsilon ^{p} \ln ^{\ell } \!\varepsilon \, \fsum_{d} \sigma _{d} (\SDA_{p,\ell } )\phi_{d} \bigg)\end{equation}
Moreover, Proposition~\ref{operateurs R} implies: $\forall R\in \Rmacro,\; \progphimacro\circ R=0$. Since $\langle \Rmacro\rangle $ is a sum of the identity and non-trivial products of elements of $\Rmacro$, we have $\progphimacro\circ \langle \Rmacro\rangle =\progphimacro\circ  \mathrm{id}=\progphimacro$. Thus applying $\progphimacro$ to \eqref{eq: raccord algebrique presque prouve} gives \eqref{eq: raccord sur les séries de sigma phi du pt de vue des champs lointains}.\\
It remains to show conversely that \eqref{eq: raccord sur les séries de sigma phi du pt de vue des champs lointains} implies \eqref{eq: raccord algebrique presque prouve}. Let us denote $A_{1}$ and $A_{2}$ the two sides of \eqref{eq: raccord algebrique presque prouve}. One can check that they are solutions of \eqref{eq: probleme de Helmholtz du champs total algebrique}. The proof of Theorem~\ref{champs totaux en fct de sigma} shows that such solutions have the following form
\[A_{i} = \langle \Rmacro\rangle \bigg(\fsum_{p,\ell ,d} \varepsilon ^{p} \ln ^{\ell } \!\varepsilon \,c_{i,p,\ell ,d} \,\phi_{d} \bigg), \quad  i\in \{1,2\}\]
for some numbers $c_{i,p,\ell ,d}$. In addition, we have $\progphimacro(A_{i} ) = \fsum_{p,\ell ,d} \varepsilon ^{p} \ln ^{\ell } \!\varepsilon \,c_{i,p,\ell ,d} \,\phi_{d}$, so $A_{i}$ is entirely determined by $\progphimacro(A_{i} )$. Since \eqref{eq: raccord sur les séries de sigma phi du pt de vue des champs lointains} states that $\progphimacro(A_{1} )=\progphimacro(A_{2} )$, we have $A_{1} =A_{2}$, i.e. \eqref{eq: raccord algebrique presque prouve}.
\end{demo}

\begin{remarque}
Since $\Rmacro$ and $\Rmicro$ only differ on powers of $\varepsilon $ due to the scaling, one could expect that $\Eps^{-1} \circ \langle \Rmicro\rangle $ is equal to $\langle \Rmacro\rangle \circ \Eps^{-1}$. This would simplify a lot \eqref{eq: raccord sur les séries de sigma phi du pt de vue des champs lointains}, which would become simply ``$\sigma _{d} (\bfuDA_{p,\ell } )=\sigma _{d} (\SDA_{p+d,\ell } )$'' for any $(d,p,\ell )$. However, $\ROmegaPi$ (resp. $\RcouchePi$, $\RGammaPi$) picks the particular solution of Equation~\ref{eq: operateur R Omega Pi} (resp. \ref{eq: operateur R couche Pi}, \ref{eq: operateur R Gamma Pi}) whose image by the $\sigma _{d}$ vanish. This means that they are the solutions in the kernels of $\progphimaicro$. It is easy to see that $\Eps^{\pm 1}$ does not map $\operatorname{Ker}\progphimaicro$ to $\operatorname{Ker}\progphimiacro$ (e.g. consider an element of $\Aemaicro$ equal to $\Im[(\alpha z)^{ \frac{\pi }{\Theta } } \log (\alpha z)]$ on $\Omega $). Therefore $\Eps^{-1} \circ \ROmegaPi$ selects other solutions than $\ROmegaPi\circ \Eps^{-1}$, and likewise for $\RcouchePi$ and $\RGammaPi$, which yields $\Eps^{-1} \circ \langle \Rmicro\rangle  \neq  \langle \Rmacro\rangle \circ \Eps^{-1}$. This observation explains the complexity of the matching relations in Theorem~\ref{condition de raccord}.
\end{remarque}

\begin{textesansboite}
For any $(p,\ell )$, let us define $\tau _{p,\ell } :\Aemaicro \rightarrow  \Amaicro(\Pi )$ by: 
$\forall \varphi =\fsum_{p'\!,\ell '} \varepsilon ^{p'} \ln ^{\ell '} \!\varepsilon \,\varphi _{p'\!,\ell '} \in \Aemaicro,\ \tau _{p,\ell } (\varphi ):= \varphi _{p,\ell }$. It allows us to define the matching coefficients, that are for any $(d,d',p,\ell )\in (\TZ*)^2 \times \NN\times \mathbb{N}$ the following complex numbers:\vspace{-1mm}
\begin{equation}{}
\cmacro_{d,d',p,\ell } := \sigma _{d} \circ  \tau _{p,\ell } \circ  \Eps^{-1} \circ  \langle \Rmicro\rangle  (\phi_{d'} ) \quad  \text{ and } \quad 
\cmicro_{d,d',p,\ell } := \sigma _{d} \circ  \tau _{p,\ell } \circ  \Eps \circ  \langle \Rmacro\rangle  (\phi_{d'} )
\end{equation}
In addition, for any $a,b\in \mathbb{R}$, we denote $\intvl{a,b}:=\{c\in \mathbb{R} \mid c-a\in \NN \text{ and } b-c\in \mathbb{N}\}$. It is a finite subset of $[a,b]$.\\

Theorem~\ref{condition de raccord} gives equations to concretely build the fields $\bfu_{p,\ell }$ and $S_{p,\ell }$ so that they match around the corner. It fixes their non-variational, which are determined by $\sigma _{d} (\bfuDA_{p,\ell } )$ when $d<0$ and $\sigma _{d} (\SDA_{p,\ell } )$ when $d>0$ (see Theorems~\ref{cadre fonctionnel des champs lointains} and \ref{cadre fonctionnel pour les champs de coin}). It also provides inductive formulas, depending on the fields with smaller $p$. Moreover these formulas have a convolutive structure w.r.t. $p$ and $\ell $.
\end{textesansboite}

\begin{theoreme}{}
\label{condition de raccord}%
The matching condition of Definition~\ref{def: raccord algebrique} is equivalent to the following set of equations:\vspace{-3mm}
\begin{subequations}
\label{eq: condition de raccord}
\begin{empheq}[left={\hspace{-5mm}\empheqlbrace\,}]{alignat=3}
	& \textstyle \forall (d,p,\ell )\in (-\TN){\times }\NN{\times }\mathbb{N}, \ \
	& \sigma _{d} (\bfuDA_{p,\ell } ) = \!\!\! \sum _{p'\in \intvl{0,p+d}} \! \sum _{\substack{d'\in \TZ*\\ p'-d'\in \intvl{0,p}}} \!\! \sum _{\ell '=0} ^{\ell } \,\cmacro_{d,d',p-p',\ell -\ell '} \cdot  \sigma _{d'} (\SDA_{p',\ell '} ) 
	\label{eq: condition de raccord 1} \\[-1mm]
	& \textstyle \forall (d,p,\ell )\in \TN\times \NN\times \mathbb{N}, \ \ 
	& \sigma _{d} (\SDA_{p,\ell } ) = \!\!\! \sum _{p'\in \intvl{0,p-d}} \! \sum _{\substack{d'\in \TZ*\\ p'+d'\in \intvl{0,p}}} \!\! \sum _{\ell '=0} ^{\ell } \,\cmicro_{d,d',p-p',\ell -\ell '} \cdot  \sigma _{d'} (\bfuDA_{p',\ell '} )
	\label{eq: condition de raccord 2}
\end{empheq}\end{subequations}\fbeq[1.5]
\end{theoreme}

\begin{demo}
Let $\Vmacro$ be the set of families of complex numbers $(\smacro_{p,d,\ell } )_{p\in \NN,d\in \TZ*\cap (\NN-p),\ell \in \mathbb{N}}$ s.t. for any $(p,d)$ only a finite number of the $\smacro_{p,d,\ell }$ are non-zero. We define $\Vmicro$ similarly by replacing ``$d\in \NN-p$'' with ``$d\in p-\NN$''. Let $\sMacro :=(\sigma _{d} (\bfuDA_{p,\ell } ))_{p,d,\ell } \in \Vmacro$ and $\sMicro :=(\sigma _{d} (\SDA_{p,\ell } ))_{p,d,\ell } \in \Vmicro$. We define the following linear maps:\vspace{-1mm}
\[\Pmacro : \Vmicro \rightarrow  \Vmacro,\ (\smicro_{p,d,\ell } )_{p,d,\ell } \ \mapsto \ \bigg( \sum _{p'\in \intvl{0,p+d}} \sum _{\substack{d'\in \TZ*\\ p'-d'\in \intvl{0,p}}} \!\sum _{\ell '=0} ^{\ell } \cmacro_{d,d',p-p',\ell -\ell '} \cdot  \smicro_{p',d',\ell '} \bigg)_{\!p,d,\ell } \vspace{-4mm}\]
\[\Pmicro : \Vmacro \rightarrow  \Vmicro,\ (\smacro_{p,d,\ell } )_{p,d,\ell } \ \mapsto \ \bigg( \sum _{p'\in \intvl{0,p-d}} \sum _{\substack{d'\in \TZ*\\ p'+d'\in \intvl{0,p}}} \!\sum _{\ell '=0} ^{\ell } \cmicro_{d,d',p-p',\ell -\ell '} \cdot  \smacro_{p',d',\ell '} \bigg)_{\!p,d,\ell }\]
These sums have a finite number of terms, so they are well-defined.\\

\underline{Step 1:} Let $(d,p,\ell )\in \TZ*\times \NN\times \mathbb{N}$. Applying $\sigma _{d} \circ \tau _{p,\ell }$ to the first equation of Lemma~\ref{raccord algebrique sur les champs} gives:
\begin{align*}
	\sigma _{d} (\bfuDA_{p,\ell } ) &= \sum _{p'\in \NN,\ell '\in \mathbb{N}} \ \sum _{d'\in \TZ*\cap (\NN-p')} \sigma _{d} \circ \tau _{p-p',\ell -\ell '} \circ \progphimacro\circ \Eps^{-1} \circ \langle \Rmicro\rangle (\phi_{d'} ) \cdot  \sigma _{d'} (\SDA_{p',\ell '} )\\
	&= \sum _{p'\in \NN} \ \sum _{d'\in \TZ*\!,p'+d'\in \NN} \ \sum _{\ell '\in \mathbb{N}} \cmacro_{d,d',p-p',\ell -\ell '} \cdot  \sigma _{d'} (\SDA_{p',\ell '} )
\end{align*}
because $\sigma _{d} \circ \tau _{p-p',\ell -\ell '} \circ \progphimacro = \sigma _{d} \circ \tau _{p-p',\ell -\ell '}$. Moreover, for any $R\in \Rmicro$, we have $\deg _{\varepsilon } \!R \in \mathbb{N}$ and $\deg _{\varepsilon } \!R-\deg _{\mathcal{A}} \!R\in \mathbb{N}$ (see Figure~\ref{fig: degres des elements de Rmacro et Rmicro}). So for $\cmacro_{d,d',p-p',\ell -\ell '}$ to be non-zero, we need $p-p'+d\in \mathbb{N}$ and $(p-p'+d)-(d-d')\in \mathbb{N}$. In addition, $\ell -\ell '\in \mathbb{N}$. Therefore, the matching condition in equivalent to  $\sMacro =\Pmacro\,\sMicro$. Similarly, it is equivalent to $\sMicro =\Pmicro\,\sMacro$.\\

\underline{Step 2:} We define the subspaces $\Vmacro_{\pm} = \{\smacro \in \Vmacro \mid \forall (p,d,\ell ),\; \smacro_{p,d,\ell } \neq 0 \Rightarrow  \pm d>0\}$ and likewise $\Vmicro_{\pm}$. So $\Vmacro = \Vmacro_{+} \oplus  \Vmacro_{-}$ and similarly for $\Vmicro$. Let $\sMacro_{\pm}$ be the components of $\sMacro$ on $\Vmacro_{\pm}$, and similarly for $\sMicro_{\pm}$. We decompose $\Pmacro$ and $\Pmicro$ on those subspaces, which gives in block matrix notation:\vspace{-1mm}
\begingroup
\def\reglerHauteurCase{\rule[-1.5ex]{0pt}{4ex}}
\def\reglerHauteurAcco{\rule{0pt}{2.2ex}}
\[\Pmacro \,=\, \raisebox{-2.2ex}{$\begin{array}{@{}c@{}c}
	\left(\begin{array}{c}
		\reglerHauteurCase \Pmacroo{}+ \\ \hline
		\reglerHauteurCase \Pmacroo{}-
	\end{array}\right)
	\! & \begin{array}{@{}c@{}}
		\left. \reglerHauteurAcco \right\} \Vmacro_+ \\[1ex]
		\left. \reglerHauteurAcco \right\} \Vmacro_-
	\end{array}\\[-1.9ex]
	\underbrace{\rule{5.5ex}{0pt}}_{\textstyle\Vmicro} & 
\end{array}$}
\qquad \text{ and } \qquad 
\Pmicro \,=\, \raisebox{-2.2ex}{$\begin{array}{@{}c@{}c@{}c}
	\left(\begin{array}{c}
		\reglerHauteurCase \Pmicroo{}+ \\ \hline
		\reglerHauteurCase \Pmicroo{}-
	\end{array}\right) 
	\;=\;& \left(\begin{array}{c|c}
		\reglerHauteurCase \Pmicroo++ & \Pmicroo-+ \\ \hline
		\reglerHauteurCase \Pmicroo+- & \Pmicroo--
	\end{array}\right)
	\! & \begin{array}{@{}c@{}}
		\left. \reglerHauteurAcco \right\} \Vmicro_+ \\[1ex]
		\left. \reglerHauteurAcco \right\} \Vmicro_-
	\end{array}\\[-1.9ex]
	& \underbrace{\rule{5ex}{0pt}}_{\textstyle\Vmacro_+} \underbrace{\rule{5ex}{0pt}}_{\textstyle\Vmacro_-} & 
\end{array}$} \vspace{-1mm}\]
\endgroup
The present theorem rewrites as: \begin{equation}{\label{eq: condition de raccord preuve}}
\left\{\begin{array}{r@{\;=\;}l}
	\sMacro_{-} & \Pmacroo{}- \,\sMicro \\[0.5mm]
	\sMicro_{+} & \Pmicroo{}+ \,\sMacro
\end{array}\right.\end{equation}
By Step~1, this is clearly a necessary condition for the matching. It remains to prove that it is sufficient. Let us show that \eqref{eq: condition de raccord preuve} implies one of the two conditions of Step~1, e.g. $\sMacro =\Pmacro\,\sMicro$.\\
Note that by Step~1 we have for any $\smacro \in \Vmacro$ and $\smicro \in \Vmicro$: $\smacro=\Pmacro\smicro \Leftrightarrow  \smicro=\Pmicro\smacro$. So $\Pmacro$ and $\Pmicro$ are inverses of each other. Hence $\Pmicro(\sMacro-\Pmacro\,\sMicro) = \Pmicro\,\sMacro-\sMicro$. Projecting this onto $\Vmicro_{+}$ we deduce:
\[\Pmicroo++(\sMacro_{+} -\Pmacroo{}+\,\sMicro) + \Pmicroo-+(\underbrace{\sMacro_{-} -\Pmacroo{}-\,\sMicro}_{=0} ) = \underbrace{\Pmicroo{}+\,\sMacro-\sMicro_{+} }_{=0} .\vspace{-5mm}\]
Let us show that $\Pmicroo++$ is injective.
\begin{sousdemo}
	Let $\smacro \in \Vmacro_{+} \setminus \{0\}$ and let us show that $\Pmicroo++ \smacro \neq 0$. Let $(p,d,\ell )$ be the smallest triplet for lexicographic order such that $\smacro_{p,d,\ell } \neq 0$. The term of $\Pmicroo++ \smacro$ of index $(p+d,d,\ell )$ is :\vspace{-2mm}
	\[\sum _{p'\in \intvl{0,p}} \sum _{\substack{d'\in \TN\\ p'+d'\in \intvl{0,p+d}}} \!\sum _{\ell '=0} ^{\ell } \cmicro_{d,d',p+d-p',\ell -\ell '} \cdot  \smacro_{p',d',\ell '} \,=\, \cmicro_{d,d,d,0} \cdot  \smacro_{p,d,\ell } .\vspace{-1mm}\]
	We claim that it is non-zero because $\cmicro_{d,d,d,0} =1$. Indeed, in the sum \vspace{-2mm}
	\[\langle \Rmicro\rangle  := \sum_{n=0}^\infty  \sum _{(R_{1} ,\dots ,R_{n} )\in \Rmicron } R_{1} \circ R_{2} \circ \dots \circ R_{n} \]
	the only term of degree $(0,0)$ w.r.t. $\mathcal{A}$ and $\varepsilon $ is the identity, which appears for $n=0$. So the component of $\Eps\circ \langle \Rmicro\rangle  (\phi_{d} )$ in $\varepsilon ^{d} \mathcal{A}_{d} (\Pi )[\ln \varepsilon ]$ is $\varepsilon ^{d} \phi_{d}$. Thus $\cmicro_{d,d,d,0} = \sigma _{d} \circ  \tau _{d,0} \circ  \Eps\circ \langle \Rmicro\rangle  (\phi_{d} ) = \sigma _{d} (\phi_{d} ) =1$. So $\Pmicroo++ \smacro \neq 0$.
\end{sousdemo}
We have proven that $\sMacro_{+} =\Pmacroo{}+\,\sMicro$. Given that $\sMacro_{-} = \Pmacroo{}- \,\sMicro$, we deduce $\sMacro =\Pmacro\,\sMicro$.
\end{demo}

\begin{paragraphesansboite}{Remarks \nvnumpar{}:}
\label{rq: raccord algebrique 2}%
\begin{itemize}
	\item Thanks to the tools of Section~\ref{subsec: operateurs de resolution dans A}, we can compute exactly and very quickly the coefficients $\cmacro_{d,d',p,\ell }$ and $\cmicro_{d,d',p,\ell }$. Moreover, these coefficients depend only on $\Theta $, $\omega $ and $(\mu _{0} ,\mu _{1} ,\rho _{0} ,\rho _{1} )$, but not on $\Omega _{1}$ nor precisely on the functions $\mu $ and $\rho $.
	\item In the sums of Theorem~\ref{condition de raccord}, the indexes $d$ and $d'$ satisfy $d-d'\in \mathbb{Z}\cap \TZ$. Indeed, on the one hand $d,d'\in \TZ*$. On the other we have in \eqref{eq: condition de raccord 1} that $p'\in \intvl{0,p+d}$ and $p'-d'\in \intvl{0,p}$, so $d-d' = (p+d-p')-(p-(p'-d'))\in \mathbb{N}-\mathbb{N}\subset \mathbb{Z}$ (and likewise in \eqref{eq: condition de raccord 2}). Note that the set $\mathbb{Z}\cap \TZ$ can be very small. If $\Theta  \in \pi \mathbb{Q}$, then $\mathbb{Z}\cap \TZ =b\mathbb{Z}$ where $\Theta =\pi  \frac{a}{b}$ with $(a,b)\in \mathbb{N}\times \mathbb{N}^*$ and $\mathrm{gcd} (a,b)=1$. Otherwise, $\mathbb{Z}\cap \TZ=\{0\}$.
\end{itemize}\fbi
\end{paragraphesansboite}

\section{Construction of the asymptotic expansion}
\label{sec: construction des champs}

\begin{textesansboite}
Equations~\ref{eq: eq vol champs lointains}--\ref{eq: eq vol champs de coin} and Theorem~\ref{condition de raccord} give the equations that the fields $u_{p,\ell }$, $U_{p,\ell }$ and $S_{p,\ell }$ must satisfy. In this section, we will build these fields according to those conditions. First of all, let us express the layer fields with the far fields, so that only two types of fields remain to build. Let $(\mathcal{U}_{n} )\in \mathbb{R}[Y]^{\mathbb{N}}$ be the unique sequence of polynomials s.t. for any $n\in \mathbb{N}^*$:
\begin{equation}{\label{eq: correcteurs de couche}}
\left\{\begin{array}{r@{\;=\;}ll}
	\mathcal{U}_{0} '' & 0 \\[0.5mm]
	\mathcal{U}_{0} '(0) & \frac{\mu _{0} }{\mu _{1}} \\[0.5mm]
	\mathcal{U}_{0} (-1) & 0 
\end{array}\right. \qquad  \text{ and } \qquad  \left\{\begin{array}{r@{\;=\;}ll}
	\mathcal{U}_{n} '' & -\mathcal{U}_{n-1} \\[0.5mm]
	\mathcal{U}_{n} '(0) & 0 \\[0.5mm]
	\mathcal{U}_{n} (-1) & 0 
\end{array}\right.\end{equation}\fbeq
\end{textesansboite}

\begin{lemme}{expression of the layer fields\titreligne}
\label{expression de Upl}%
Let us assume that the fields $u_{p,\ell }$ and $U_{p,\ell }$ are regular enough (we will check later that they are). \eqref{eq: eq vol champs de couche} implies for any $(p,\ell )\in \NN \times \mathbb{N}$ and $(x,Y)\in \couche$:\vspace{-1mm}
\begin{equation}{\label{eq: expression de Upl}}
U_{p,\ell } (x,Y) = \sum_{n=0}^\infty  (\partial _{x} ^2 +k_{1} ^2 )^{n} \partial _{y} u_{p-1-2n,\ell } (x,0) \cdot  \mathcal{U}_{n} (Y).\vspace{-1mm}\end{equation}
where this sum has a finite number of non-zero terms by the convention: $\forall p\in \mathbb{R} \setminus \NN,\forall \ell \in \mathbb{N},\; u_{p,\ell } :=0$.
\end{lemme}

\begin{demo}
There exists an increasing sequence $(p_{m} )$ s.t. $\NN=\{p_{m} \mid m\in \mathbb{N}\}$. So we can prove the result by induction on $p\in \NN$. For $p=0$, \eqref{eq: eq vol champs de couche} states: \[\left\{\begin{array}{r@{\;=\;}ll}
	\partial _{Y} ^2  U_{0,\ell } & 0 &\text{ in }\couche \\[0.5mm]
	\partial _{Y} U_{0,\ell } & 0 &\text{ on }\Gamma  \\[0.5mm]
	U_{0,\ell } & 0 &\text{ on }\borddroit
\end{array}\right.\] so $U_{0,\ell } =0$ for any $\ell $. It is coherent \eqref{eq: expression de Upl} (which is a sum of zeros in this case).\\
Next, for the inductive step, we assume that \eqref{eq: expression de Upl} holds for ranks smaller that $p$. \eqref{eq: eq vol champs de couche} gives: \[\left\{\begin{array}{r@{}l@{}l@{\,}l}
	\partial _{Y} ^2  U_{p,\ell } &\;=\; -(\partial _{x} ^2 +k_{1} ^2 )U_{p-2,\ell }
	&\;=\; -\sum_{n=0}^\infty  (\partial _{x} ^2 +k_{1} ^2 )^{n+1} \partial _{y} u_{(p-2)-1-2n,\ell } (x,0) \cdot  \mathcal{U}_{n} (Y) & \\[0.5mm]
	&&\;=\; -\sum_{n=1}^\infty  (\partial _{x} ^2 +k_{1} ^2 )^{\red{n}} \partial _{y} u_{p-1-2n,\ell } (x,0) \cdot  \mathcal{U}_{\red{n-1}} (Y) &\text{ in }\couche \\
	\partial _{Y} U_{p,\ell } &\;=\; \frac{\mu _{0} }{\mu _{1}} \,\partial _{y} u_{p-1,\ell } &&\text{ on }\Gamma  \\[0.5mm]
	U_{p,\ell } &\;=\; 0 &&\text{ on }\borddroit
\end{array}\right.\] It is easy to see that $\sum_{n=0}^\infty  (\partial _{x} ^2 +k_{1} ^2 )^{n} \partial _{y} u_{p-1-2n,\ell } (x,0) \cdot  \mathcal{U}_{n} (Y)$ is the only solution of this.
\end{demo}

\begin{textesansboite}
We saw in Section~\ref{sec: raccord de coin} that far and corner fields possess singularities when $r\rightarrow 0$, resp. $\bfr\rightarrow \infty $. So the usual variational frameworks are not sufficient to build these fields and we need to design new frameworks. It is done in Sections~\ref{subsec: cadre fonctionnel des champs lointains} and \ref{subsec: cadre fonctionnel des champs de coin}. In both sections we start by introducing the natural space $H_{\mathrm{var}}$ in which an ad hoc variational problem is well-posed. Then we define a bigger space $H := H_{\mathrm{var}} +\chi  \,\mathcal{A}(D) = \{u+\chi \varphi  \mid u\in H_{\mathrm{var}} , \varphi \in \mathcal{A}(D)\}$ that contains the singularities, where $\chi $ is a $\mathcal{C}^{\infty }$ truncation function in the vicinity of 0 (for $u_{p,\ell }$) or infinity (for $S_{p,\ell }$), and $D\in \{\Omega ,\Pi \}$. Next we determine the elements of $H_{\mathrm{var}} \cap  \chi  \mathcal{A}(D)$, which allows us to define on $H$ the linear forms $\sigma _{d}$ associated to the singularities. Finally we show that some model problems are well-posed in $H$.\\

Before we start, the following lemma is a tool to estimate the behavior at 0 and $\infty $ of functions of $\mathcal{A}$.
\end{textesansboite}

\begin{lemme}{}
\label{croissance des fcts de A}%
Let $a<b$ in $\mathbb{R}$, $n\in \mathbb{N}^*$, $(d_{i} ,\ell _{i} )_{i\in [\![1,n]\!]}$ be $n$ distinct elements of $\mathbb{R}\times \mathbb{N}$, $(f_{i} )_{i\in [\![1,n]\!]} \in (\mathcal{C}([a,b],\mathbb{C}) \setminus \{0\})^{n}$ and:\vspace{-2mm}
\[\varphi :(r,\theta )\in \mathbb{R}_+^*\times [a,b] \;\mapsto \;\sum_{i=1}^n r^{d_{i}} \ln ^{\ell _{i}} \! r\cdot  f_{i} (\theta ).\]
Then there is an interval $I\subset [a,b]$ with non-empty interior, $c\in \mathbb{R}_+^*$ and $r_{1} ,r_{2} \in \mathbb{R}$ such that:\vspace{-1mm}
\[\forall r\in (0,r_{1} ), \forall \theta \in I,\ \ |\varphi (r,\theta )|>c\, r^{\min _{j} d_{j} } \quad  \text{ and } \quad  \forall r\in (r_{2} ,\infty ), \forall \theta \in I,\ \ |\varphi (r,\theta )|>c\, r^{\max _{j} d_{j} } .\]\fbeq[1.8]
\end{lemme}

\begin{demo}
Let $j\in [\![1,n]\!]$ be s.t. $(d_{j} ,\ell _{j} )$ is maximal for the lexicographic order. Let $I\subset [a,b]$ be a non trivial interval on which $|f_{j} |$ is greater than a positive constant. Since $\sum _{i\neq j} r^{d_{i}} \ln ^{\ell _{i}} \! r\cdot  f_{i} (\theta ) =o(r^{d_{j}} \ln ^{\ell _{j}} r)$, we have when $r\rightarrow \infty $ and $\theta \in I$: $|\varphi (r,\theta )| \gtrsim  r^{d_{j}} \ln ^{\ell _{j}} r \gtrsim  r^{d_{j}}$. And we can similarly treat the vicinity of 0.
\end{demo}

\begin{textesansboite}
We can apply Lemma~\ref{croissance des fcts de A} to any $\varphi \in \mathcal{A}_{d} (\Omega )$ with $(a,b):=(0,\Theta )$, or to any $\varphi \in \mathcal{A}_{d} (\couche)$ with $(a,b):=(-1,0()$ (replacing the variables $(r,\theta )$ by $(x,Y)$). We can also apply it to $\partial _{r} \varphi $ and $\partial _{\theta } \varphi $ when $\varphi \in \mathcal{A}_{d} (\Omega )$ and to $\partial _{x} \varphi $ and $\partial _{Y} \varphi $ when $\varphi \in \mathcal{A}_{d} (\couche)$.
\end{textesansboite}

\begin{paragraphesansboite}{Definition:}
For any $D\in \{\Pi ,\Omega ,\Gamma ,\couche\}$ and $\varphi \in \mathcal{A}(D)$, we denote $\mindeg \varphi :=\sup \{d\in \mathbb{R} \mid \varphi \in \sum _{q\geqslant d} \mathcal{A}_{q} (D)\}$ and $\maxdeg \varphi :=\inf \{d\in \mathbb{R} \mid \varphi \in \sum _{q\leqslant d} \mathcal{A}_{q} (D)\}$.
\end{paragraphesansboite}

\subsection{Existence and uniqueness for far fields-like problems}
\label{subsec: cadre fonctionnel des champs lointains}

\begin{textesansboite}
We denote $H^{1/2} _{00} (\Gamma )$ the set of functions of $H^{1/2} (\Gamma )$ whose extension by 0 to $\partial \Omega =\Gamma \cup \{(0,0)\}\cup \bordgauche$ is in $H^{1/2} (\partial \Omega )$. Using Lax-Milgram theorem, it is easy to prove the following lemma.
\end{textesansboite}

\begin{lemme}{the Helmholtz problem in $H^{1} (\Omega )$\titreligne}
\label{pb de Helmholtz variationnel}%
Let $f\in (H^{1} (\Omega ))'$ and $g\in H^{1/2} _{00} (\Gamma )$. The following system has a unique solution in $H^{1} (\Omega )$. \[\left\{\begin{array}{r@{\;=\;}ll}
	\mu _{0} \Delta  u +\omega ^2 \rho _{0} u & f &\text{ in }\Omega  \\[0.5mm]
	u & g &\text{ on }\Gamma  \\[0.5mm]
	u & 0 &\text{ on }\bordgauche
\end{array}\right.\]\fbeq
\end{lemme}

\begin{textesansboite}
Let $\khiz$ be a radial function of $\mathcal{C}^{\infty } (\mathbb{R}^2 )$ equal to 1 in the vicinity of 0 and to 0 in the vicinity of infinity. The appropriate space to build the far fields is $H^{1} (\Omega )+\khiz \mathcal{A}(\Omega )$. One can check that it does not depend on the choice of $\khiz$.
\end{textesansboite}

\begin{paragraphesansboite}{Lemma \nvnumpar{}:}
\label{H1 inter A}%
$\displaystyle H^{1} (\Omega ) \cap  \khiz \mathcal{A}(\Omega ) = \khiz \sum _{d>0} \mathcal{A}_{d} (\Omega ).$\fbeq
\end{paragraphesansboite}

\begin{demo}
The inclusion $\supset$ is easy to check, so we focus on $\subset $. Let $\varphi \in \mathcal{A}(\Omega ) \setminus \{0\}$ be s.t. $\khiz \varphi \in H^{1} (\Omega )$. There is $d\in \mathbb{R}$, $\varphi _{1} \in \mathcal{A}_{d} (\Omega ) \setminus \{0\}$ and $\varphi _{2} \in \sum _{q>d} \mathcal{A}_{q} (\Omega )$ s.t. $\varphi =\varphi _{1} +\varphi _{2}$. If $\partial _{\theta } \varphi _{1}$ were null everywhere, then would so too $\varphi _{1}$ because $\varphi _{1|\bordgauche} =0$. But we assumed the contrary, so $\partial _{\theta } \varphi _{1} \neq 0$. Thus Lemma~\ref{croissance des fcts de A} implies that there is a non-trivial interval $I\subset [0,\Theta ]$ s.t. $|\partial _{\theta } \varphi _{1} (r,\theta )| \gtrsim  r^{d}$ when $r\rightarrow 0$ and $\theta \in I$. Finally $r^{-1} \partial _{\theta } (\khiz \varphi ) \in L^2 (\Omega )$ implies that $d>0$.
\end{demo}

\begin{paragraphesansboite}{Définition \nvnumpar{}:}
\label{def: sigma dans H1 + khi A}%
For any $u\in H^{1} (\Omega )+\khiz \mathcal{A}(\Omega )$ and $d\in -\TN$, we denote $\sigma _{d} (u):=\sigma _{d} (\varphi )$ where $\varphi \in \mathcal{A}(\Omega )$ is s.t. $u-\khiz \varphi \in H^{1} (\Omega )$ (see Definition~\ref{def: sigma} for $\sigma _{d} (\varphi )$). It does not depend on the choice of $\varphi $ thanks to Lemma~\ref{H1 inter A}.
\end{paragraphesansboite}

\begin{theoreme}{existence and uniqueness for a far fields-like model problem\titreligne}
\label{cadre fonctionnel des champs lointains}%
Let $f\in (H^{1} (\Omega ))'$, $g\in H^{1/2} _{00} (\Gamma )+\khiz \mathcal{A}(\Gamma )$ and $(s_{d} )\in \mathbb{C}^{-\TN}$ with finite support. The following system has a unique solution in $H^{1} (\Omega )+\khiz \mathcal{A}(\Omega )$.
\[\left\{\begin{array}{r@{\;=\;}ll}
	\mu _{0} \Delta  u +\omega ^2 \rho _{0} u & f &\text{ in }\Omega  \\[0.5mm]
	u & g &\text{ on }\Gamma  \\[0.5mm]
	u & 0 &\text{ on }\bordgauche\\[0.5mm]
	\sigma _{d} (u)& s_{d} & \ \forall d\in -\TN
\end{array}\right.\]\fbeq
\end{theoreme}

\begin{demo}
Let us show the existence first, and then the uniqueness.\\
\underline{Existence:} Let $\tilde{g} \in H^{1/2} _{00} (\Gamma )$ and $\varphi \in \mathcal{A}(\Gamma )$ be s.t. $g= \tilde{g} +\khiz \varphi $. We look for the solution in the form $u= \tilde{u}+\khiz \psi $ with $\tilde{u} \in H^{1} (\Omega )$ and $\psi \in \mathcal{A}(\Omega )$. Let
\[\psi ^{+} := \sum_{n=0}^\infty  (-k_{0} ^2 \ROmegaOmega)^{n} \bigg(\RGammaOmega(\varphi )+\sum _{d\in -\TN} s_{d} \, \phiOmega_{d} \bigg) \;\in \;\Amacro(\Omega )\]
(where $\phiOmega_{d} :=\phi_{d|\Omega }$) and $\psi :=T_{\leqslant 2} (\psi ^{+} )$. Using Proposition~\ref{operateurs R dans Omega} one can check by calculus the first system below (see also Lemma~\ref{DA algebrique des champs lointains concatenes} for a similar result). Then, the second system below derives from $\psi \in \mathcal{A}(\Omega )$, $\psi  = \psi ^{+} -( \mathrm{id}-T_{\leqslant 2} )(\psi ^{+} )$ and $\deg \Delta =-2$ (by Lemma~\ref{operateurs de derivation dans A}).
\[\left\{\begin{array}{@{\;}r@{\;=\;}l@{\;}l@{}}
	(\mu _{0} \Delta  +\omega ^2 \rho _{0} )\psi ^{+} & 0 &\text{ in }\Omega  \\[0.4mm]
	\psi ^{+} & \varphi  &\text{ on }\Gamma  \\[0.4mm]
	\psi ^{+} & 0 &\text{ on }\bordgauche\\[0.4mm]
	\sigma _{d} (\psi ^{+} )& s_{d} & \ \forall d\in -\TN
\end{array}\right. \quad  \mathrm{so} \quad  \left\{\begin{array}{@{\;}r@{}l@{}}
	(\mu _{0} \Delta  +\omega ^2 \rho _{0} )\psi  &\;\in \, \mathcal{A}(\Omega )\cap  \fsum\limits_{d>0} \mathcal{A}_{d} (\Omega ) = \sum\limits _{d>0} \mathcal{A}_{d} (\Omega ) \\[0.5mm]
	\psi _{|\Gamma } -\varphi  &\;\in \, \mathcal{A}(\Gamma )\cap  \fsum\limits_{d>2} \mathcal{A}_{d} (\Gamma ) = \sum\limits _{d>2} \mathcal{A}_{d} (\Gamma ) \\
	\psi _{|\bordgauche} &\;=\; 0 \\[0.5mm]
	\sigma _{q} (\psi  )&\;=\; s_{d} \;\ \forall d\in -\TN
\end{array}\right.\]
Therefore $(\mu _{0} \Delta  +\omega ^2 \rho _{0} )(\khiz \psi )\in L^2 (\Omega )$ and $\khiz \,(\psi _{|\Gamma } -\varphi )\in H^{1/2} _{00} (\Gamma )$. Finally, thanks to Lemma~\ref{pb de Helmholtz variationnel}, we can take $\tilde{u}$ as the unique solution in $H^{1} (\Omega )$ of: \[\left\{\begin{array}{r@{\;=\;}ll}
	\mu _{0} \Delta  \tilde{u} +\omega ^2 \rho _{0} \tilde{u} & f - (\mu _{0} \Delta +\omega ^2 \rho _{0} )(\khiz \psi ) &\text{ in }\Omega  \\[0.5mm]
	\tilde{u} & \tilde{g} -\khiz \,(\psi -\varphi )  &\text{ on }\Gamma  \\[0.5mm]
	\tilde{u} & 0 &\text{ on }\bordgauche
\end{array}\right.\]

\underline{Uniqueness:} Let $u$ be a homogeneous solution. There is $\tilde{u} \in H^{1} (\Omega )$ and $\psi \in \mathcal{A}(\Omega )$ such that $u= \tilde{u}+\khiz \psi $. Let $\varphi :=\Delta \psi +k_{0} ^2 \psi $. By Lemma~\ref{operateurs de derivation dans A}, $\varphi \in \mathcal{A}(\Omega )$. Let us show that $\mindeg \varphi >-2$.
\begin{sousdemo}
	We denote $d:=\mindeg \varphi $ and $\widetilde{\varphi } :=(\Delta +k_{0} ^2 )(\khiz \psi )$. By Lemma~\ref{croissance des fcts de A}, there is a non-trivial interval $I\subset [0,\Theta ]$ s.t. when $r\rightarrow 0$ and $\theta \in I$: $\widetilde{\varphi }(r,\theta ) =\varphi (r,\theta )\gtrsim  r^{d}$. Besides, $\widetilde{\varphi } = -(\Delta +k_{0} ^2 ) \tilde{u} \,\in  (H^{1} _{0} (\Omega ))'$. Let us test it with $\zeta _{q} :(r,\theta ) \mapsto  \widetilde{\varphi } (r,\theta )\cdot r^{q-2d-2} \,\big(1-\khiz(2^{1/q} r)\big)\,\chi (\theta )$ where $q>0$ and $\chi \in \mathcal{D}(0,\Theta ) \setminus \{0\}$ is everywhere non negative. Since $\zeta _{q} \in \mathcal{D}(\Omega )$:
	\[\langle  \widetilde{\varphi }, \zeta _{q} \rangle  := \int _{0} ^{\infty } \int _{0} ^{\Theta } \widetilde{\varphi } \cdot \zeta _{q} \cdot r\mathrm{d}\theta \,\mathrm{d}r \,\lesssim \, \|\zeta _{q} \|_{H^{1} (\Omega )} \qquad  \forall q>0\]
	Let us assume by contradiction that $d\leqslant -2$. It is not difficult to check that when $q\rightarrow 0^{+}$:
	\[\langle  \widetilde{\varphi }, \zeta _{q} \rangle  \gtrsim  \int _{c_{1} 2^{- \frac{1}{q} }} ^{c_{2}} r^{d} \,r^{q-d-2} \,r\mathrm{d}r \gtrsim  \frac{1}{q} \quad  \text{ and } \quad  \|\zeta _{q} \|_{H^{1}} \lesssim  \|\nabla  \zeta _{q} \|_{L^2 }  \lesssim  \bigg(\int _{c_{3} 2^{- \frac{1}{q} }} ^{c_{4}} r^{2(q-d-3)} \,r\mathrm{d}r \!\bigg)^{ \frac{1}{2} } \!\lesssim  \frac{1}{\sqrt{q}}\]
	with some constants $c_{i}$ and Poincaré's inequality. But it contradicts $\langle  \widetilde{\varphi }, \zeta _{q} \rangle  \lesssim  \|\zeta _{q} \|_{H^{1}}$. So $d>-2$.		
\end{sousdemo}
Now let us show that $\mindeg \psi >0$.
\begin{sousdemo}
	Let us assume the contrary. Then there is $d\leqslant 0$, $\psi _{1} \in \mathcal{A}_{d} (\Omega ) \setminus \{0\}$ and $\psi _{2} \in \sum _{q>d} \mathcal{A}_{q} (\Omega )$ s.t. $\psi =\psi _{1} +\psi _{2}$. Lemma~\ref{operateurs de derivation dans A} states that $\deg  \Delta =-2$, so $\Delta \psi _{1} \in \mathcal{A}_{d-2} (\Omega )$ and $\Delta \psi _{2} \in \sum _{q>d-2} \mathcal{A}_{q} (\Omega )$. In addition, $\varphi  \in \sum _{q>-2} \mathcal{A}_{q} (\Omega )$ and $\mathcal{A}(\Omega ) = \bigoplus_{q\in \mathbb{R}} \mathcal{A}_{q} (\Omega )$ according to Lemma~\ref{decomposition de A en sev}. So identifying the coordinate of $\Delta \psi +k_{0} ^2 \psi  =\varphi $ in $\mathcal{A}_{d-2} (\Omega )$ gives $\Delta \psi _{1} =0$.\\
	Moreover, $\psi _{|\Gamma } = - \tilde{u}_{|\Gamma } \in H^{1/2} (\Gamma )$, so $\psi _{|\Gamma } \in \sum _{q>0} \mathcal{A}_{q} (\Gamma )$. Hence: \[\left\{\begin{array}{r@{\;=\;}ll}
		\Delta  \psi _{1} & 0 &\text{ in }\Omega \\[0.5mm]
		\psi _{1} & 0 &\text{ on }\bordgauche \cup \Gamma 
	\end{array}\right.\]
	By Proposition~\ref{pb de Laplace homogene dans A}, it implies $d\in -\TN$ and $\psi _{1} =\sigma _{d} (\psi _{1} )\,\phiOmega_{d}$. However, $\sigma _{d} (\psi _{1} )=\sigma _{d} (\psi )=\sigma _{d} (u)=0$, which contradicts $\psi _{1} \neq 0$.
\end{sousdemo}
Finally, Lemma~\ref{H1 inter A} implies $\khiz \psi \in H^{1} (\Omega )$, so $u\in H^{1} (\Omega )$. This means that $u=0$ by Lemma~\ref{pb de Helmholtz variationnel}.
\end{demo}

\subsection{Existence and uniqueness for corner fields-like problems}
\label{subsec: cadre fonctionnel des champs de coin}

\begin{definition}{the variational space $\Vcoin$\titreligne}
\label{def: V}%
We define $\Vcoin := \{v\in H^{1} _{\mathrm{loc}} (\Omega _{1} ) \mid \nabla  v \in L^2 (\Omega _{1} ) \text{ and } v_{|\partial \Omega _{1}} = 0\}$ and the norm $\|v\|_{\Vcoin} := \|\nabla v\|_{L^2 (\Omega _{1} )}$.\\
Here ``$v_{|\partial \Omega _{1}} = 0$'' means that $\khiz( \frac{\cdot }{R} ) \,v\in H^{1} _{0} (\Omega _{1} )$ for any $R>0$, where $\khiz( \frac{\cdot }{R} ): (X,Y)\mapsto  \khiz\big( \frac{X}{R} \frac{,Y}{R} \big)$.
\end{definition}

\begin{paragraphesansboite}{Lemma \numpar{}:}
\label{inegalite de Hardy dans Omega1}%
Any $v\in \Vcoin$ satisfies $\|v\|_{H^{1} (\couche\cap \{X>\Rcoin\})} \lesssim  \|v\|_{\Vcoin}$ \,and \;$\displaystyle \Big\| \frac{v}{1+r} \Big\|_{L^2 (\Omega _{1} )} \lesssim  \|v\|_{\Vcoin}$.\fbeq[0.5]
\end{paragraphesansboite}

\newcommand{\coucheR}{A}
\begin{demo}
We denote $\coucheR :=\couche \cap \{X>\Rcoin\}$ and $B:=B(0,\Rcoin)$, and we recall that $\Rcoin>0$ is s.t. $\Omega _{1} \setminus B=\Pi  \setminus B$.\\
Poincaré's inequality gives for a.e. $X>\Rcoin$: $\|v(X,\cdot )\|_{L^2 (-1,0)} \; \lesssim  \; \|\partial _{Y} v(X,\cdot )\|_{L^2 (-1,0)}$. Integrating w.r.t. $X$ the square of this then gives: $\|v\|_{L^2 (\coucheR)} \lesssim  \|\partial _{Y} v\|_{L^2 (\coucheR)} \leqslant  \|v\|_{\Vcoin}$. So $\|v\|_{H^{1} (\coucheR)} \lesssim  \|v\|_{\Vcoin}$.\\
Next, a generalized Poincaré inequality on circular slices of $\Omega  \setminus B$ gives for a.e. $r>\Rcoin$: $\|v(r,\cdot )\|_{L^2 (0,\pi )} \lesssim  \|\partial _{\theta } v(r,\cdot )\|_{L^2 (0,\pi )} + |v(r,\theta =0)|$. Then integrating w.r.t. $r$ the square of this times $ \frac{1}{1+r}$ gives:\vspace{-2mm}
\[\int _{\Omega  \setminus B} \frac{v^2 }{(1+r)r} \,r\mathrm{d}\theta \,\mathrm{d}r \;\lesssim \; \int _{\Omega  \setminus B} \frac{(\partial _{\theta } v)^2 }{(1+r)r} \,r\mathrm{d}\theta \,\mathrm{d}r + \int _{\Gamma  \setminus B} \frac{v^2 }{1+r} \mathrm{d}r.\]
Therefore $\| \frac{v}{1+r} \|_{L^2 (\Omega  \setminus B)} ^2  \lesssim  \|\nabla  v\|_{L^2 (\Omega  \setminus B)} ^2 + \|v\|_{L^2 (\Gamma  \setminus B)} ^2 $. But we know that $\|v\|_{L^2 (\Gamma  \setminus B)} \lesssim  \|v\|_{H^{1} (A)} \lesssim  \|v\|_{\Vcoin}$, so $\| \frac{v}{1+r} \|_{L^2 (\Pi  \setminus B)} \lesssim \|v\|_{\Vcoin}$. Finally Poincaré's inequality in $B$ also gives $\| \frac{v}{1+r} \|_{L^2 (\Omega _{1} \cap  B)} \lesssim \|v\|_{\Vcoin}$.
\end{demo}

\begin{lemme}{the Poisson problem in $\Omega _{1}$\titreligne}
\label{pb variationnel de Poisson dans Omega1}%
Let $F:\Omega _{1} \rightarrow \mathbb{C}$ be s.t. $(1+r)F\in L^2 (\Omega _{1} )$, and $g\in L^2 (\Gamma \cap \{X>\Rcoin\})$. The following system is well-posed in $\Vcoin$.\vspace{-2mm}
\[\left\{\begin{array}{r@{\;=\;}ll}
	\operatorname{div}(\mu  \nabla  S) &F &\text{ in }\Omega _{1} \setminus (\Gamma \cap \{X>\Rcoin\})  \\[0.5mm]
	S &0 &\text{ on }\partial \Omega _{1} \\[0.5mm]
	S_{|Y=0^{+}} - S_{|Y=0^{-}} & 0 &\text{ on }\Gamma \cap \{X>\Rcoin\} \\[0.5mm]
	\mu _{0} \partial _{Y} S_{|Y=0^{+}} - \mu _{1} \partial _{y} S_{|Y=0^{-}} & g &\text{ on }\Gamma \cap \{X>\Rcoin\} 
\end{array}\right.\]\fbeq
\end{lemme}

\begin{demo}
The variational formulation of this problem is: \[ \forall v\in \Vcoin,\qquad  \int _{\Omega _{1}} \mu  \nabla S\cdot \nabla v = -\int _{\Omega _{1}} Fv+ \int _{\Gamma \cap \{X>\Rcoin\}} gv.\]
The left-hand side is coercive by definition of $\Vcoin$. Lemma~\ref{inegalite de Hardy dans Omega1} ensures that the right-hand side is continuous. Moreover it is easy to see that $\Vcoin$ is complete. So we can conclude using the Lax-Milgram theorem.
\end{demo}

\begin{textesansboite}
Let $\khii$ be a function of $\mathcal{C}^{\infty } (\mathbb{R}^2 )$ equal to 0 on $B(0,\Rcoin)$ and 1 in a vicinity of infinity. The proper space to build the corner fields is $\Vcoin +\khii \mathcal{A}(\Pi )$. One can check that it does not depend on the choice of $\khii$.
\end{textesansboite}

\begin{paragraphesansboite}{Lemma \nvnumpar{}:}
\label{V inter A}%
$\Vcoin\cap  \khii \mathcal{A}(\Pi ) = \big\{ \khii \varphi  \;\big|\; \varphi \in \mathcal{A}(\Pi ),\; \maxdeg(\varphi _{|\Omega } )<0\; \text{ and } \,\maxdeg(\varphi _{|\couche} )<- \frac{1}{2} \big\}.$
\end{paragraphesansboite}

\begin{demo}
The inclusion $\supset$ is easy to check, so we focus on $\subset $. Let $\varphi \in \mathcal{A}(\Pi )$ be s.t. $\khii \varphi \in \Vcoin$. The same method as Lemma~\ref{H1 inter A} shows that $\maxdeg \varphi _{|\Omega } <0$. And using that $\partial _{Y} (\khii \varphi ) \in L^2 (\couche)$ we likewise get $\maxdeg(\varphi _{|\couche} )<- \frac{1}{2}$.
\end{demo}

\begin{paragraphesansboite}{Definition \nvnumpar{}:}
\label{def: sigma dans V + khi A}%
For any $S\in \Vcoin +\khii \mathcal{A}(\Pi )$ and $d\in \TN$, we denote $\sigma _{d} (S):=\sigma _{d} (\varphi )$ where $\varphi \in \mathcal{A}(\Pi )$ is s.t. $S-\khii \varphi \in \Vcoin$. It does not depend on the choice of $\varphi $ thanks to Lemma~\ref{V inter A}.
\end{paragraphesansboite}

\begin{theoreme}{existence and uniqueness for a corner fields-like model problem\titreligne}
\label{cadre fonctionnel pour les champs de coin}%
Let $F:\Omega _{1} \rightarrow \mathbb{C}$ be s.t. $(1+r)F\in L^2 (\Omega _{1} )$, $\varphi :\Pi \rightarrow \mathbb{C}$ be s.t. $\varphi _{|\Omega } \in \mathcal{A}(\Omega )$ and $\varphi _{|\couche} \in \mathcal{A}(\couche)$, and $(s_{d} )\in \mathbb{C}^{\TN}$ with finite support. The following system has a unique solution in $V+\khii \mathcal{A}(\Pi )$.
\[\left\{\begin{array}{r@{\;=\;}ll}
	\operatorname{div}(\mu  \nabla  S) & F +\khii \varphi  &\text{ in }\Omega _{1} \\[0.5mm]
	S & 0 &\text{ on }\partial \Omega _{1} \\[0.5mm]
	\sigma _{d} (S)& s_{d} & \ \forall d\in \TN
\end{array}\right.\]\fbeq
\end{theoreme}

\begin{remarque}
$\varphi $ is not just an element of $\mathcal{A}(\Pi )$ because it may be discontinuous on $\Gamma $.
\end{remarque}

\begin{demo}
This proof is similar to Theorem~\ref{cadre fonctionnel des champs lointains}. Let us show the existence first, and then the uniqueness.\\
\underline{Existence:} We look for the solution in the form $S= \tilde{S}+\khii \psi $ with $\tilde{S} \in \Vcoin$ and $\psi \in \mathcal{A}(\Pi )$. More precisely we denote (using notation $\langle .\rangle $ from Definition~\ref{def: serie geometrique d'operateurs 1D})
\[\varphi ^{-} := \left\langle -\RcouchePi\circ \partial _{X|\couche} ^2 ,\; \frac{\mu _{0} }{\mu _{1}} \RGammaPi\circ \partial _{Y|\Gamma ,Y=0^{+}} \right\rangle  \bigg( \frac{1}{\mu _{0}} \ROmegaPi (\varphi _{|\Omega } )+ \frac{1}{\mu _{1}} \RcouchePi (\varphi _{|\couche} ) + \sum _{d\in \TN} s_{d} \,\phi_{d} \bigg) \,\in \, \Amicro(\Pi )\]
and $\psi :=T_{\geqslant -2} (\varphi ^{-} )$. By Lemma~\ref{DA algebrique des champs de coin}, we have the first system below. Then, similarly to the proof of Theorem~\ref{cadre fonctionnel des champs lointains}, one can check that it implies the second system below.
\[\left\{\begin{array}{r@{\;=\;}ll@{}}
	\mu _{0} \Delta  \varphi ^{-} & \varphi  &\text{ in }\Omega  \\[0.5mm]
	\mu _{1} \Delta  \varphi ^{-} & \varphi  &\text{ in }\couche \\[0.5mm]
	[\varphi ^{-} ]_{\Gamma } & 0 &\text{ on }\Gamma \\[0.5mm]
	[\mu \partial _{Y} \varphi ^{-} ]_{\Gamma } & 0 &\text{ on }\Gamma  \\[0.5mm]
	\varphi ^{-} & 0 &\text{ on }\bordgauche \cup \borddroit\\[0.5mm]
	\sigma _{d} (\varphi ^{-} ) & s_{d} & \ \forall d\in \TN
\end{array}\right. \quad  \mathrm{so} \quad  \left\{\begin{array}{@{\;}r@{}l@{}}
	(\mu _{0} \Delta  \psi  -\varphi )_{|\Omega } &\;\in \, \mathcal{A}(\Omega )\cap  \fsum\limits_{d<-4} \mathcal{A}_{d} (\Omega ) = \sum\limits _{d<-4} \mathcal{A}_{d} (\Omega ) \\[0.5mm]
	(\mu _{1} \Delta  \psi  -\varphi )_{|\couche} &\;\in \, \mathcal{A}(\couche)\cap  \fsum\limits_{d<-2} \mathcal{A}_{d} (\couche) = \sum\limits _{d<-2} \mathcal{A}_{d} (\couche) \\[-0.5mm]
	[\psi ]_{\Gamma } &\;=\; 0\\[0.5mm]
	[\mu \partial _{Y} \psi ]_{\Gamma } &\;\in \, \mathcal{A}(\Gamma )\cap  \fsum\limits_{d<-2} \mathcal{A}_{d} (\Gamma ) = \sum\limits _{d<-2} \mathcal{A}_{d} (\Gamma ) \\[-0.5mm]
	\psi _{|\bordgauche \cup \borddroit} &\;=\; 0 \\[0.5mm]
	\sigma _{q} (\psi )&\;=\; s_{d} \;\ \forall d\in \TN
\end{array}\right.\]
where $[{\dots }]_{\Gamma }$ stands for the jump on $\Gamma $. 
Therefore $(1+r)\,\big[\!\operatorname{div}\!\big(\mu  \nabla  (\khii \psi )\big)-\khii \varphi \big] \in L^2 (\Omega _{1} )$ and $[\mu \partial _{Y} (\khii \psi )]_{\Gamma } \in L^2 (\Gamma \cap \{X>\Rcoin\})$. Finally, thanks to Lemma~\ref{pb variationnel de Poisson dans Omega1}, we take $\tilde{S}$ as the unique solution in $V$ of: \[\left\{\begin{array}{r@{\;=\;}ll}
	\operatorname{div}(\mu  \nabla  \tilde{S}) & F+\khii \varphi - \operatorname{div}\!\big(\mu  \nabla  (\khii \psi )\big)&\text{ in }\Omega _{1} \setminus (\Gamma \cap \{X>\Rcoin\}) \\[0.5mm]
	\tilde{S} &0 &\text{ on }\partial \Omega _{1} \\[0.5mm]
	[\mu \partial _{Y} \tilde{S}]_{\Gamma } & -[\mu \partial _{Y} (\khii \psi )]_{\Gamma } &\text{ on }\Gamma \cap \{X>\Rcoin\} 
\end{array}\right.\]

\underline{Uniqueness:} Let $S$ be a homogeneous solution and let us show that $S=0$. There is $\tilde{S} \in \Vcoin$ and $\psi \in \mathcal{A}(\Omega )$ s.t. $S= \tilde{S}+\khii \psi $. The proof decomposes into the following steps, all proven by contradiction.
\begin{enumerate}
	\item \underline{$\maxdeg(\Delta \psi _{|\Omega } )<-2$\,:} Like in step 1 of the uniqueness proof of Theorem~\ref{cadre fonctionnel des champs lointains}, we test the inequality:
	\[\int _{\Omega } \Delta (\khii \psi )\cdot \zeta  =-\int _{\Omega } \Delta  \tilde{S} \cdot \zeta  =\int _{\Omega } \nabla  \tilde{S} \cdot \nabla  \zeta  \lesssim  \|\zeta \|_{V} \qquad  \forall \zeta  \in \mathcal{D}(\Omega )\]
	with $\zeta _{q} :(r,\theta ) \mapsto  \Delta (\khii \psi )(r,\theta )\cdot r^{-q-2d-2} \,\khiz(2^{-1/q} r)\,\chi (\theta )$ when $q\rightarrow 0^{+}$, where $d:=\maxdeg(\Delta \psi _{|\Omega } )$ and $\chi \in \mathcal{D}(0,\Theta ) \setminus \{0\}$ is everywhere non negative.
	\item \underline{$\maxdeg(\Delta \psi _{|\couche} )<- \frac{1}{2}$\,:} This step works like the previous one.
	\item \underline{$\maxdeg [\mu \partial _{Y} \psi ]_{\Gamma } <- \frac{1}{2}$\,:} Let $d:= \maxdeg  [\mu \partial _{Y} \psi ]_{\Gamma }$ and $\varphi := [\mu \partial _{Y} (\khii \psi )]_{\Gamma }$.	We have $\nabla  \tilde{S} \in L^2 (\Omega _{1} )$ and, by steps 1 and 2, $\operatorname{div}(\mu  \nabla  \tilde{S}) =-\operatorname{div}\!\big(\mu  \nabla  (\khii \psi )\big) \in L^2 (\Omega _{1} )$. Thus $\varphi  =-[\mu \partial _{Y} \tilde{S}]_{\Gamma } \in H^{-1/2} (\Gamma )$. Like previously, we test it with $\zeta _{q} :x \mapsto  \varphi (x) \cdot x^{-q-2d-1} \,\khiz(2^{-1/q} x)$ when $q\rightarrow 0^{+}$. If $d\geqslant  - \frac{1}{2}$, we get $\langle  \varphi , \zeta _{q} \rangle _{\Gamma } \gtrsim  \frac{1}{q}$ and $\|\zeta _{q} \|_{H^{1/2} (\Gamma )} \lesssim  \|\zeta _{q} \|_{H^{1} (\Gamma )} \lesssim  \frac{1}{\sqrt{q}}$, which contradicts $\langle  \varphi , \zeta _{q} \rangle _{\Gamma } \lesssim  \|\zeta _{q} \|_{H^{1/2} (\Gamma )}$.
	\item \underline{$\maxdeg(\psi _{|\couche} )<- \frac{1}{2}$\,:} Let $d:=\maxdeg(\psi _{|\couche} )$. There is $\psi _{1} \in \mathcal{A}_{d} (\couche) \setminus \{0\}$ and $\psi _{2} \in \sum _{q<d} \mathcal{A}_{q} (\couche)$ s.t. $\psi _{|\couche} =\psi _{1} +\psi _{2}$. By Lemma~\ref{operateurs de derivation dans A}, $\deg (\partial _{x|\couche} ^2 )=-2$, $\deg (\partial _{Y|\Gamma ,Y=0^{+}} )=-1$ and $\deg (\partial _{Y|\couche} ^2 )=\deg (\partial _{Y|\Gamma ,Y=0^{-}} )=0$. In addition, Proposition~\ref{decomposition de A en sev} states that $\mathcal{A}(D) = \bigoplus_{q\in \mathbb{R}} \mathcal{A}_{q} (D)$ for any $D\in \{\couche,\Gamma \}$. So by taking the coordinates of $\Delta \psi _{|\couche}$ in $\mathcal{A}_{d} (\couche)$ and of $[\mu \partial _{Y} \psi ]_{\Gamma }$ in $\mathcal{A}_{d} (\Gamma )$, we get if $d\geqslant  - \frac{1}{2}$: \[\left\{\begin{array}{r@{\;=\;}ll}
		\partial _{Y} ^2  \psi _{1} & 0 &\text{ in }\couche \\[0.5mm]
		\partial _{Y|\Gamma ,Y=0^{-}} \psi _{1} & 0 &\text{ on }\Gamma 
	\end{array}\right.\]
	But since $\psi _{|\borddroit} =0$, it implies $\psi _{1} =0$, which is contradictory.
	\item \underline{$\maxdeg(\psi _{|\Omega } )<0$\,:} We denote $d:=\maxdeg(\psi _{|\Omega } )$. There is $\psi _{1} \in \mathcal{A}_{d} (\Omega ) \setminus \{0\}$ and $\psi _{2} \in \sum _{q<d} \mathcal{A}_{q} (\Omega )$ s.t. $\psi _{|\Omega } =\psi _{1} +\psi _{2}$. Since $\deg (\Delta _{|\Omega } )=-2$, taking the coordinate of $\Delta \psi _{|\Omega }$ in $\mathcal{A}_{d-2} (\Omega )$ and of $\psi _{|\Gamma }$ in $\mathcal{A}_{d} (\couche)$ gives if $d\geqslant 0$: \[\left\{\begin{array}{r@{\;=\;}ll}
		\Delta  \psi _{1} & 0 &\text{ in }\Omega  \\[0.5mm]
		\psi _{1} & 0 &\text{ on }\Gamma 
	\end{array}\right.\]
	By Proposition~\ref{pb de Laplace homogene dans A}, it implies $d\in \TN$ and $\psi _{1} =\sigma _{d} (\psi _{1} )\phi_{q}$. But $\sigma _{d} (\psi _{1} )=\sigma _{d} (\psi )=\sigma _{d} (S)=0$.
\end{enumerate}
Finally, Lemma~\ref{V inter A} implies $\khii \psi \in V$, so $S\in V$. So the uniqueness in Lemma~\ref{pb variationnel de Poisson dans Omega1} implies $S=0$.
\end{demo}

\subsection{Construction of the fields}
\label{subsec: construction des termes}
	
\begin{boiterondesecable}{Definition \nvnumpar{}:}
\label{def: construction des termes}%
For any $(p,\ell )\in (\mathbb{R} \setminus \NN)\times \mathbb{N}$ we denote by convention $u_{p,\ell } =0$, $U_{p,\ell } =0$ and $S_{p,\ell } =0$. We define by induction on $p\in \NN$ that for any $\ell \in \mathbb{N}$:
\begin{itemize}
	\item $\displaystyle U_{p,\ell } :(x,Y)\in \couche \mapsto  \sum_{n=0}^\infty  (\partial _{x} ^2 +k_{1} ^2 )^{n} \partial _{y} u_{p-1-2n,\ell } (x,0) \cdot  \mathcal{U}_{n} (Y)$. \leavevmode\hfill\setlabeleq\label{eq: definition des champs de couche}
	\item $u_{p,\ell }$ is the unique solution in $H^{1} (\Omega )+\khiz \mathcal{A}(\Omega )$ of :
	\begin{equation}{\label{eq: definition des champs lointains}}
		\hspace*{-5mm}\left\{\begin{array}{r@{\;=\;}ll@{}}
		\mu _{0} \Delta  u_{p,\ell } +\omega ^2 \rho _{0} u_{p,\ell } & \fSource\,\delta _{p,0} \,\delta _{\ell ,0} &\text{ in }\Omega  \\[0.5mm]
		u_{p,\ell } & 0 &\text{ on }\bordgauche\\[0.5mm]
		u_{p,\ell } & U_{p,\ell } &\text{ on }\Gamma \\[-2mm]
		\sigma _{d} (u_{p,\ell } ) &\displaystyle \!\!\sum _{p'\in \intvl{0,p+d}} \sum _{\substack{d'\in \TZ*\\ p'-d'\in \intvl{0,p}}} \!\sum _{\ell '=0} ^{\ell } \,\cmacro_{d,d',p-p',\ell -\ell '} \cdot  \sigma _{d'} (\SDA_{p',\ell '} ) &\ \forall d\in -\TN
	\end{array}\right.\vspace{-2mm}\end{equation}
	\item $S_{p,\ell }$ is the unique solution in $\Vcoin+\khii \mathcal{A}(\Pi )$ of :
	\begin{equation}{\label{eq: definition des champs de coin}}
		\left\{\begin{array}{r@{\;=\;}ll}
		\operatorname{div}(\mu  \nabla  S_{p,\ell } ) & -\omega ^2 \rho  \,S_{p-2,\ell } &\text{ in }\Omega _{1} \\[0.5mm]
		S_{p,\ell } & 0 &\text{ on }\partial  \Omega _{1} \\[-2mm]
		\sigma _{d} (S_{p,\ell } ) &\displaystyle \!\!\sum _{p'\in \intvl{0,p-d}} \sum _{\substack{d'\in \TZ*\\ p'+d'\in \intvl{0,p}}} \!\sum _{\ell '=0} ^{\ell } \,\cmicro_{d,d',p-p',\ell -\ell '} \cdot  \sigma _{d'} (\bfuDA_{p',\ell '} ) &\ \forall d\in \TN
	\end{array}\right.\end{equation}\fbeq
\end{itemize}\fbi
\end{boiterondesecable}

\begin{textesansboite}
We will show that these fields are well-defined at the same time as the following proposition.
\end{textesansboite}

\begin{paragraphesansboite}{Proposition \nvnumpar{}:}
\label{confimation des DA des champs}%
For any $(p,\ell )\in \NN\times \mathbb{N}$ there exist $\bfuDA_{p,\ell } \in \Amacro(\Pi )$ and $\SDA_{p,\ell } \in \Amicro(\Pi )$ s.t., for any $d\in \mathbb{R}$, $\bfu_{p,\ell } = T_{\leqslant d} (\bfuDA_{p,\ell } ) +o_{\partial } (\bfr^{d} )$ when $\bfr\rightarrow 0$ and $S_{p,\ell } = T_{\geqslant d} (\SDA_{p,\ell } ) + o_{\partial } (\bfr^{d} )$ when $\bfr\rightarrow \infty $.
\end{paragraphesansboite}

\begin{demo}[Justification of Definition~\ref{def: construction des termes} and Proposition~\ref{confimation des DA des champs}]
Let us show by induction on $p$ that for any $\ell $:
\begin{itemize}
	\item[(H1)] $u_{p,\ell }$, $U_{p,\ell }$ and $S_{p,\ell }$ exist and are unique,
	\item[(H2)] Proposition~\ref{confimation des DA des champs} is true at rank $(p,\ell )$,
	\item[(H3)] $\forall m\in \mathbb{N},\; \partial _{x} ^{m} \partial _{y} u_{p,\ell |\Gamma } \in  H^{1} _{0} (\Gamma )+\khiz \mathcal{A}(\Gamma )$ (additional property that will be useful during the proof).
\end{itemize}
Since all fields are null for $p<0$, the initial case is trivial and only the inductive step remains to prove. Let $p\in \NN$. Let us assume (H1)--(H3) at any rank $p'<p$ and prove it at rank $p$.
\begin{enumerate}
	\item Existence and uniqueness of the fields:
	\begin{itemize}
 		\item \underline{$U_{p,\ell }$\,:} By (H3), the $\partial _{y} u_{p',\ell |\Gamma }$ with $p' <p$ are all in $H^{m} _{\mathrm{loc}} (\Gamma )$ for any $m\in \mathbb{N}$. So the functions $(\partial _{x} ^2 +k_{1} ^2 )^{n} \partial _{y} u_{p-1-2n,\ell |\Gamma }$ are continuous and Formula \eqref{eq: definition des champs de couche} is well-defined at any point of $\couche$.
 		\item \underline{$u_{p,\ell }$\,:} By (H2), $\SDA_{p',\ell '}$ exists for any $p'<p$ and $\ell '$, so the coefficients $\sigma _{d'} (\SDA_{p',\ell '} )$ are well-defined. Thus, \eqref{eq: definition des champs lointains} defines well $\sigma _{d} (u_{p,\ell } )$. And we have: $\sigma _{d} (u_{p,\ell } ) \neq 0 \Rightarrow  \intvl{0,p+d} \neq  \varnothing  \Rightarrow  d\geqslant -p$, so $(\sigma _{d} (u_{p,\ell } ))_{d\in -\TN}$ has finite support. Next, we apply Theorem~\ref{cadre fonctionnel des champs lointains}. Its hypotheses are satisfied, because (H3) implies:\vspace{-1mm}
 		\[U_{p,\ell |\Gamma } = \sum_{n=0}^\infty  (\partial _{x} ^2 +k_{1} ^2 )^{n} \partial _{y} u_{p-1-2n,\ell } (x,0) \cdot  \mathcal{U}_{n} (0) \,\in \,H^{1} _{0} (\Gamma )+\khiz \mathcal{A}(\Gamma ) \,\subset \, H^{1/2} _{00} (\Gamma )+\khiz \mathcal{A}(\Gamma ) \vspace{-1mm}\] 
 		\item \underline{$S_{p,\ell }$\,:} Similarly, (H2) implies that $(\sigma _{d} (S_{p,\ell } ))_{d\in \TN}$ is well defined and has finite support. Next, we apply Theorem~\ref{DA des champs de coin}. We need to check its hypotheses, i.e. there is $f:\Omega _{1} \rightarrow \mathbb{C}$ and $\varphi :\Pi \rightarrow \mathbb{C}$ s.t. $\omega ^2 \rho \,S_{p-2,\ell } = F+\khii \varphi $, $(1+r)F\in L^2 (\Omega _{1} )$, $\varphi _{|\Omega } \in \mathcal{A}(\Omega )$ and $\varphi _{|\couche} \in \mathcal{A}(\couche)$. By (H2), there is $\psi \in \mathcal{A}(\Pi )$ s.t. $S_{p-2,\ell } =\psi  +o_{\partial } (\bfr^{-4} )$. So it suffices to set $\varphi :=\omega ^2 \rho  \psi $ and $F:=\omega ^2 \rho  \,S_{p-2,\ell } -\khii \varphi $.
 	\end{itemize}
	
 	\item Asymptotic expansions:
	\begin{itemize}
 		\item \underline{$\bfu_{p,\ell }$\,:} We use Theorem~\ref{DA des champs lointains concatenes}. To do so, we must check that there is $g_{\Gamma } \in \Amacro(\Gamma )$ and $g_{\couche} \in \Amacro(\couche)$ s.t.:\vspace{-2mm}
		\[\forall d\in \mathbb{R}, \qquad  \left\{\begin{array}{r@{\;=\;}ll}
			\partial _{Y} ^2  U_{p,\ell } & T_{\leqslant d} (g_{\couche} )+o_{\partial } (x^{d} ) &\text{ in }\couche\\[0.5mm]
 			\partial _{Y} U_{p,\ell |Y=0^{-}} & T_{\leqslant d} (g_{\Gamma } )+o_{\partial } (x^{d} ) &\text{ on }\Gamma 
 		\end{array}\right.\] Given the definition of $U_{p,\ell }$ \eqref{eq: definition des champs de couche} , it suffices to show that, for any $p'<p$ there is $h\in \Amacro(\Gamma )$ s.t.: $\forall d\in \mathbb{R},\; \partial _{n} u_{p'\!,\ell |\Gamma } = T_{\leqslant d} (h)+o_{\partial } (x^{d} )$. But it derives from (H2).
 		\item \underline{$S_{p,\ell }$\,:} Similarly, we use Theorem~\ref{DA des champs de coin} thanks to (H2).
	\end{itemize}
	
	\item Let $m\in \mathbb{N}$. The asymptotic expansion of $\bfu_{p,\ell }$ implies: $\exists h\in \mathcal{A}(\Gamma ),\;\partial _{x} ^{m} \partial _{y} u_{p,\ell |\Gamma } = h+o_{\partial } (x^{1} )$. So there is $x_{0} \in \mathbb{R}_+^*$ s.t. $\partial _{x} ^{m} \partial _{y} u_{p,\ell |\Gamma } -h$ is $H^{1}$ on $\Gamma \cap \{x<x_{0} \}=(0,x_{0} )\times \{0\}$ and it vanishes at 0. To prove (H3) at rank $p$, it remains to show that $\partial _{x} ^{m} \partial _{y} u_{p,\ell } \in H^{1} ((\frac{x_{0} }{2} ,\infty )\times \{0\})$. To do so, it suffices to get $u \in H^{m+3} ((\frac{x_{0} }{2} ,\infty )\times (0, \frac{\delta }{2} ))$ with $\delta  :=\operatorname{dist}(\operatorname{supp}(f),\Gamma )$. But it follows from classical elliptic regularity because on one side $\mu _{0} \Delta u_{p,\ell } +\omega ^2 \rho _{0} u_{p,\ell } =0$ on $(\frac{x_{0} }{4} ,\infty )\times (0,\delta )$, and on the other $u_{p,\ell |\Gamma } =U_{p,\ell |\Gamma } \in H^{m+3} ((\frac{x_{0} }{4} ,\infty )\times \{0\})$ by \eqref{eq: definition des champs de couche} and (H3).
\end{enumerate}\fbi
\end{demo}

\begin{paragraphesansboite}{Proposition \nvnumpar{}:}
\label{borne sur l}%
$\forall p\in \NN, \exists n_{p} \in \mathbb{N},\forall \ell >n_{p} ,\ \ (u_{p,\ell } =0 \text{ and } U_{p,\ell } =0 \text{ and } S_{p,\ell } =0)$.
\end{paragraphesansboite}

\begin{demo}
For any $(d,d',p)\in (\TZ*)^2 \times \NN$, one has $\Eps^{-1} \circ  \langle \Rmicro\rangle  (\phi_{d'} )\in \Aemacro$, so by definition of $\Aemacro$ there is $n$ (depending of $d,d',p$) s.t. for any $\ell >n$:  $\cmacro_{d,d',p,\ell } := \sigma _{d} \circ  \tau _{p,\ell } \circ \Eps^{-1} \circ  \langle \Rmicro\rangle  (\phi_{d'} ) =0$. The same is true for the coefficients $\cmicro_{d,d',p,\ell }$. Finally the result follows by induction from Definition~\ref{def: construction des termes}.
\end{demo}

\begin{paragraphesansboite}{Proposition \nvnumpar{}:}
\label{confimation de l'hyp de raccord}%
Ansatz~\ref{ansatz sur Aeps} and the matching condition of Definition~\ref{def: raccord algebrique} are satisfied.
\end{paragraphesansboite}

\begin{demo}
First let us note that for any $(p,\ell )$
\begin{equation}{\label{eq: coincidence des sigma singuliers}}
\textstyle \sigma _{d} (u_{p,\ell } )=\sigma _{d} (\bfuDA_{p,\ell } ) \text{ when }  d\in -\TN \quad  \text{ and } \quad  \sigma _{d} (S_{p,\ell } )=\sigma _{d} (\SDA_{p,\ell } ) \text{ when }  d\in \TN. \end{equation}
where $\sigma _{d} (u_{p,\ell } )$, resp. $\sigma _{d} (S_{p,\ell } )$, is set by Definition~\ref{def: sigma dans H1 + khi A}, resp.~\ref{def: sigma dans V + khi A}, whereas $\sigma _{d} (\bfuDA_{p,\ell } )$ and $\sigma _{d} (\SDA_{p,\ell } )$ rest on the definition of $\sigma _{d}$ on $\Amaicro(\Pi )$ at page~\pageref{def: sigma dans Amaicro}. Indeed we have $u_{p,\ell } -T_{\leqslant 1} (\bfuDA_{p,\ell } )_{|\Omega } = o_{\partial } (r^{1} )$, so $\khiz \cdot (u_{p,\ell } -T_{\leqslant 1} (\bfuDA_{p,\ell } )_{|\Omega } )\in H^{1} (\Omega )$ which implies by Definition~\ref{def: sigma dans H1 + khi A} that $\sigma _{d} (u_{p,\ell } )=\sigma _{d} \big(T_{\leqslant 1} (\bfuDA_{p,\ell } )_{|\Omega } \big) =\sigma _{d} (\bfuDA_{p,\ell } )$. Likewise for $S_{p,\ell }$.\\
Now, given Proposition~\ref{borne sur l}, to prove Ansatz~\ref{ansatz sur Aeps} it suffices to check that, for any $(p,\ell )$, $\bfuDA_{p,\ell } \in \fsum_{d\in \NN-p} \mathcal{A}_{d} (\Pi )$ and $\SDA_{p,\ell } \in \fsum_{d\in p-\NN} \mathcal{A}_{d} (\Pi )$. Let us show it for $\bfuDA_{p,\ell }$ only, using induction.\vspace{1mm}
\begin{sousdemo}
	Since all fields vanish for $p<0$, only the inductive step is non-trivial. \eqref{eq: bfuDA p,l} states that
	\[\bfuDA_{p,\ell } =\sum_{n=0}^\infty  (-k_{0} ^2 \,\ROmegaPi)^{n} \bigg(\underbrace{\RcouchePi\circ (\partial _{x|\couche} ^2 +k_{1} ^2 )(\bfuDA_{p-2,\ell } ) + {\textstyle \frac{\mu _{0} }{\mu _{1}} }\RGammaPi\circ \partial _{y|\Gamma ,y=0^{+}} (\bfuDA_{p-1,\ell } )}_{:=A} + \!\! \fsum_{d\in \TZ*} \!\!\! \sigma _{d} (\bfuDA_{p,\ell } ) \,\phi_{d} \bigg). \vspace{-1mm}\]
	Since $\deg  \ROmegaPi \in \mathbb{N}$, it suffices to show that the big brackets belong to $\fsum_{d\in \NN-p} \mathcal{A}_{d} (\Pi )$. This is true for $A$ using the induction hypothesis and $\deg \RcouchePi=0$, $\deg (\partial _{x|\couche} ^2 )=-2$ and $\deg (\RGammaPi\circ \partial _{y|\Gamma ,y=0^{+}} )=-1$ (see Proposition~\ref{operateurs R} and Lemma~\ref{operateurs de derivation dans A}). It remains to show that $\forall d\in \TZ*\!,\; \sigma _{d} (\bfuDA_{p,\ell } ) \neq 0 \Rightarrow d\in \NN-p$. This last assertion holds because on one hand $\TN\subset \NN-p$, and on the other \eqref{eq: coincidence des sigma singuliers} and \eqref{eq: definition des champs lointains} imply: $\forall d\in -\TN,\; \sigma _{d} (\bfuDA_{p,\ell } ) \neq 0 \Leftrightarrow  \sigma _{d} (u_{p,\ell } ) \neq 0\Rightarrow  \intvl{0,p+d} \neq  \varnothing  \Rightarrow  d\in \NN-p$.
\end{sousdemo}
Finally \eqref{eq: definition des champs lointains}--\eqref{eq: definition des champs de coin} and \eqref{eq: coincidence des sigma singuliers} show that the matching relations of Theorem~\ref{condition de raccord} are satisfied, and we can apply Theorem~\ref{condition de raccord} thanks to Ansatz~\ref{ansatz sur Aeps}. Thus the matching condition of Definition~\ref{def: raccord algebrique} is satisfied.
\end{demo}

\subsection{Practical way to build the far fields}
\label{subsec: complements sur la construction}

\begin{textesansboite}
This section shows how to build directly the far fields without computing the layer and corner fields. Thanks to the explicit expression of the layer fields in \eqref{eq: definition des champs de couche}, the layer is replaced by boundary conditions on $\Gamma $, while the corner fields are replaced by corner conditions depending on corner profiles.
\end{textesansboite}

\begin{definition}{corner profiles\titreligne}
\label{def: correcteurs de coin}%
Let $d\in \TN$. We denote $(\mathcal{S}_{d,n} )_{n\in \mathbb{N}}$ the unique sequence of $\Vcoin+\khii \mathcal{A}(\Pi )$ s.t. for any $n\in \mathbb{N}^*$:
\[\left\{\begin{array}{r@{\;=\;}l@{\;}l}
	\operatorname{div}(\mu  \nabla  \mathcal{S}_{d,0} ) & 0 &\text{ in }\Omega _{1} \\[0.5mm]
	\mathcal{S}_{d,0} & 0 &\text{ on }\partial \Omega _{1} \\[0.5mm]
	\sigma _{q} (\mathcal{S}_{d,0} )& \delta _{d,q} & \ \forall q\in \TN
\end{array}\right. \quad  \text{ and } \quad  \ \left\{\begin{array}{r@{\;=\;}l@{\;}l}
	\operatorname{div}(\mu  \nabla  \mathcal{S}_{d,n} ) & -\omega ^2 \rho \,\mathcal{S}_{d,n-1} &\text{ in }\Omega _{1} \\[0.5mm]
	\mathcal{S}_{d,n} & 0 &\text{ on }\partial \Omega _{1} \\[0.5mm]
	\sigma _{q} (\mathcal{S}_{d,n} )& 0 & \ \forall q\in \TN
\end{array}\right.\]
And for any $(d,n)$ we denote $\mathcal{S}^{\infty } _{d,n}$ the element of $\Amicro(\Pi )$ s.t. $\forall d\in \mathbb{R},\; \mathcal{S}_{d,n} = T_{\geqslant d} (\mathcal{S}^{\infty } _{d,n} )+o_{\partial } (\bfr^{d} )$.\\
These objects are well-defined thanks to Theorems~\ref{cadre fonctionnel pour les champs de coin} and \ref{DA des champs de coin}. The proof is the same as for Definition~\ref{def: construction des termes} and Proposition~\ref{confimation des DA des champs}.
\end{definition}

\begin{textesansboite}
Using \eqref{eq: definition des champs de coin} and the uniqueness in Theorem~\ref{cadre fonctionnel pour les champs de coin}, one can easily show by induction on $p$:
\begin{equation}{\label{eq: expression de Spl avec les correcteurs}}
\forall (p,\ell )\in \NN\times \mathbb{N},\quad  S_{p,\ell } = \sum_{n=0}^\infty  \; \sum _{d\in \TN \cap (p-2n-\NN)} \sigma _{d} (S_{p-2n,\ell } )\,\mathcal{S}_{d,n} \end{equation}
Thus, the same holds replacing $S_{p,\ell }$ and $\mathcal{S}_{d,n}$ by resp. $\SDA_{p,\ell }$ and $\mathcal{S}^{\infty } _{d,n}$.\\
For any $(d,d',p,\ell )\in (-\TN)\times \TZ*\times \NN\times \mathbb{N}$, we introduce the corner coefficient $\cfinal_{d,d'\!,p,\ell }$, we have
\begin{subequations}
\begin{align}
\cfinal_{d,d'\!,p,\ell } &:= \! \sum _{\substack{(p_{1} ,p_{2} )\in (\mathbb{N}-d)\times (\mathbb{N}+d')\\ n\in \mathbb{N},\, p_{1} +p_{2} +2n=p}}
	\ \sum _{\substack{\red{d_{1} \in \TZ*}\\ p_{1} +d_{1} \in \mathbb{N}}} \ \sum _{\substack{d_{2} \in \TN\\ p_{2} -d_{2} \in \mathbb{N}}}
	\ \sum _{\substack{(\ell _{1} ,\ell _{2} )\in \mathbb{N}^2 \\ \ell _{1} +\ell _{2} =\ell }}
	\cmacro_{d,d_{1} ,p_{1} ,\ell _{1} }\cdot \sigma _{d_{1}} (\mathcal{S} ^{\infty } _{d_{2} ,n} )\cdot \cmicro_{d_{2} ,d'\!,p_{2} ,\ell _{2} }
	\label{eq: cfinal 1}\\[2mm]
&\,= \! \sum _{\substack{(p_{1} ,p_{2} )\in (\mathbb{N}-d)\times (\mathbb{N}+d')\\ n\in \mathbb{N},\, p_{1} +p_{2} +2n=p}}
	\ \sum _{\substack{\red{d_{1} \in -\TN}\\ p_{1} +d_{1} \in \mathbb{N}}} \ \sum _{\substack{d_{2} \in \TN\\ p_{2} -d_{2} \in \mathbb{N}}}
	\ \sum _{\substack{(\ell _{1} ,\ell _{2} )\in \mathbb{N}^2 \\ \ell _{1} +\ell _{2} =\ell }}
	\cmacro_{d,d_{1} ,p_{1} ,\ell _{1} }\cdot \sigma _{d_{1}} (\mathcal{S} ^{\infty } _{d_{2} ,n} )\cdot \cmicro_{d_{2} ,d'\!,p_{2} ,\ell _{2} }\\
&\quad  + \!\! \sum _{\substack{(p_{1} ,p_{2} )\in (\mathbb{N}-d)\times (\mathbb{N}+d')\\ p_{1} +p_{2} =p}}
	\ \sum _{\substack{\red{d_{1} \in \TN}\\ p_{1} +d_{1} \in \mathbb{N},\, p_{2} -d_{1} \in \mathbb{N}}}
	\ \sum _{\substack{(\ell _{1} ,\ell _{2} )\in \mathbb{N}^2 \\ \ell _{1} +\ell _{2} =\ell }}
	\cmacro_{d,d_{1} ,p_{1} ,\ell _{1} }\cdot \cmicro_{d_{1} ,d'\!,p_{2} ,\ell _{2} }
	\label{eq: cfinal 3}
\end{align}
\end{subequations}
(the two given formulas are equal because, when $d_{1} >0$, $\sigma _{d_{1}} (\mathcal{S}^{\infty } _{d_{2} ,n} )=\sigma _{d_{1}} (\mathcal{S}_{d_{2} ,n} ) =\delta _{d_{1} ,d_{2} } \,\delta _{n,0}$). Like in Remark~\ref{rq: raccord algebrique 2}, we have $d-d_{1} ,d_{2} -d' \in \mathbb{Z}\cap \TZ$ in \eqref{eq: cfinal 1}. Thus, if $d-d'\not\in \mathbb{Z}\cap \TZ$, line~\eqref{eq: cfinal 3} vanishes. Moreover, for any  $(p,\ell )$, we denote $\uDA _{p,\ell } := \bfuDA_{p,\ell |\Omega } \in \Amacro(\Omega )$, which satisfies: $\forall d\in \mathbb{R},\; u_{p,\ell } \underset{r\rightarrow 0}{=} T_{\leqslant d} (\uDA_{p,\ell } )+o_{\partial } (r^{d} )$.\vspace{-1mm}
\end{textesansboite}

\begin{theoreme}{direct construction of the far fields\titreligne}
\label{construction directe des u p,l}%
Let $(T_{n} )_{n\in \mathbb{N}}$ be the sequence of Taylor coefficients of the tangent: $\forall t\in (- \frac{\pi }{2} , \frac{\pi }{2} ),\; \tan  t = \sum_{n=0}^\infty  T_{n} t^{2n+1}$.\\
$(u_{p,\ell } )_{p\in \NN,\ell \in \mathbb{N}}$ is the unique family of $H^{1} (\Omega )+\khiz \mathcal{A}(\Omega )$ s.t. for any $(p,\ell )\in \NN\times \mathbb{N}$:
\[\left\{\begin{array}{r@{\;=\;}ll@{}}
	\mu _{0} \Delta  u_{p,\ell } +\omega ^2 \rho _{0} u_{p,\ell } & \fSource \,\delta _{p,0} \,\delta _{\ell ,0} &\text{ in }\Omega  \\[0.5mm]
	u_{p,\ell } & 0 &\text{ on }\bordgauche\\[0.5mm]
	u_{p,\ell } &\displaystyle \sum_{n=0}^\infty  \frac{\mu _{0} }{\mu _{1}} \,T_{n} \cdot  (\partial _{x} ^2 +k_{1} ^2 )^{n} \partial _{y} u_{p-1-2n,\ell } &\text{ on }\Gamma  \\[-1.5mm]
	\sigma _{d} (u_{p,\ell } )&\displaystyle \!\! \sum _{\substack{p'\in \NN\\ p-p'\in \NN+ \frac{2\pi }{\Theta } }} \; \sum _{\substack{d'\in \TZ*\\ p'+d'\in \intvl{0,p+d}}} \! \sum _{\ell '=0} ^{\ell } \cfinal_{d,d',p-p',\ell -\ell '} \cdot \sigma _{d'} (\uDA_{p',\ell '} ) & \ \forall d\in -\TN
\end{array}\right.\]\fbeq
\end{theoreme}

\newcommand{\bull}{\text{\tiny $\bullet$},\text{\tiny $\bullet$}}
\begin{demo}
Uniqueness follows from the uniqueness in Theorem~\ref{cadre fonctionnel des champs lointains}. So it suffices to prove that the far fields satisfy the equations above.\medskip

\underline{Boundary condition:} By \eqref{eq: definition des champs de couche} and \eqref{eq: definition des champs lointains}, we have
\[u_{p,\ell } (x,0) =U_{p,\ell } (x,0) = \sum_{n=0}^\infty  (\partial _{x} ^2 +k_{1} ^2 )^{n} \partial _{y} u_{p-1-2n,\ell } (x,0) \cdot  \mathcal{U}_{n} (0).\]
So we need to calculate $\mathcal{U}_{n} (0)$. Let $\mathrm{U} :(Y,t) \mapsto  \sum_{n=0}^\infty  \mathcal{U}_{n} (Y)\, t^{2n+1}$. The definition of $(\mathcal{U}_{n} )$ (reminded in \eqref{eq: construction directe des u p,l eq1} below) formally implies a differential equation on $\mathrm{U}$ given in \eqref{eq: construction directe des u p,l eq2}.\par
\vspace{2.5mm}
\begin{tabular}{@{}c@{}c@{}}
\begin{minipage}{0.5\linewidth}
\begin{equation}{\label{eq: construction directe des u p,l eq1}}
\forall n\in \mathbb{N},\ \left\{\begin{array}{r@{\;=\;}ll}
	\mathcal{U}_{n} '' & -\mathcal{U}_{n-1} \\[0.5mm]
	\mathcal{U}_{n} ' (0) & \frac{\mu _{0} }{\mu _{1}} \delta _{n,0} \\[0.5mm]
	\mathcal{U}_{n} (-1) & 0 
\end{array}\right.\end{equation}
\end{minipage}
&
\begin{minipage}{0.5\linewidth}
\begin{equation}{\label{eq: construction directe des u p,l eq2}}
\left\{\begin{array}{r@{\;=\;}ll}
\partial _{Y} ^2 \text{U}  & -t^2 \text{U}  \\[0.5mm]
\partial _{Y} \text{U}  _{|Y=0} & \mathcal{U}_{0} '(0)\cdot t = \frac{\mu _{0} }{\mu _{1}} t \\[0.5mm]
\text{U}  _{|Y=-1} & 0 
\end{array}\right.\end{equation}
\end{minipage}
\end{tabular}\par
\vspace{2.5mm}
Thus $\mathrm{U} (Y,t) = \frac{\mu _{0} }{\mu _{1}} \frac{\sin (t(Y+1))}{\cos (t)} :=\varphi (Y,t)$. This is formal, as we do not know whether the series $\mathrm{U}$ converges. However, there is a sequence of polynomial functions $(\Phi _{n} )$ s.t., for any $(Y,t)\in (-1,0)\times (- \frac{\pi }{2} , \frac{\pi }{2} )$, $\varphi (Y,t) = \sum_{n=0}^\infty  \Phi _{n} (Y)\, t^{2n+1}$ (because $\varphi $ is odd w.r.t. $t$). Since $\varphi $ satisfies \eqref{eq: construction directe des u p,l eq2}, $(\Phi _{n} )$ is a solution of \eqref{eq: construction directe des u p,l eq1}. But this solution is unique, so $(\Phi _{n} )=(\mathcal{U}_{n} )$ and $\mathrm{U} =\varphi $. Hence $\forall t\in (- \frac{\pi }{2} , \frac{\pi }{2} ),\; \sum_{n=0}^\infty  \mathcal{U}_{n} (0)\, t^{2n+1} = \frac{\mu _{0} }{\mu _{1}} \tan (t)$. That is to say $\forall n\in \mathbb{N},\; \mathcal{U}_{n} (0) = \frac{\mu _{0} }{\mu _{1}} \,T_{n}$.\medskip\smallskip

\underline{Corner condition:} For any $(d,p,\ell )\in \TN\times \mathbb{R}\times \mathbb{N}$, let $\widetilde{\mathcal{S}} ^{\infty } _{d,p,\ell } := \mathcal{S}^{\infty } _{d,p/2}$ if $\frac{p}{2} \in \mathbb{N}$ and $\ell =0$ and $\widetilde{\mathcal{S}} ^{\infty } _{d,p,\ell } :=0$ otherwise. For any family $(x_{p,\ell } )$, we denote $x_{\bull} :=(x_{p,\ell } )$ the family itself. We also denote $*$ the convolution product w.r.t. $(p,\ell )$. Let $d\in -\TN$. We have
\begin{align*}
	\sigma _{d} (u_{\bull} ) &= \sum _{d_{1} \in \TZ*} \cmacro_{d,d_{1} ,\bull} * \sigma _{d_{1}} (\SDA_{\bull} ) 
	\tag*{by \eqref{eq: definition des champs lointains}}\\
	&= \sum _{d_{1} \in \TZ*} \cmacro_{d,d_{1} ,\bull} * \sum _{d_{2} \in \TN} \sigma _{d_{1}} ( \widetilde{\mathcal{S}} ^{\infty } _{d_{2} ,\bull} )*\sigma _{d_{2}} (S_{\bull} ) \tag*{by \eqref{eq: expression de Spl avec les correcteurs}} \\
	&= \sum _{d_{1} \in \TZ*} \cmacro_{d,d_{1} ,\bull} * \sum _{d_{2} \in \TN} \sigma _{d_{1}} ( \widetilde{\mathcal{S}} ^{\infty } _{d_{2} ,\bull} )*\sum _{d'\in \TZ*} \cmicro_{d_{2} ,d'\!,\bull} * \sigma _{d'} (\bfuDA_{\bull} )
	\tag*{by \eqref{eq: definition des champs de coin}}\\
	&= \sum _{d'\in \TZ*} \bigg(\sum _{d_{1} \in \TZ*} \sum _{d_{2} \in \TN} \cmacro_{d,d_{1} ,\bull} * \sigma _{d_{1}} ( \widetilde{\mathcal{S}}_{d_{2} ,\bull} )* \cmicro_{d_{2} ,d'\!,\bull} \bigg)* \sigma _{d'} (\uDA_{\bull} )\\
	&= \sum _{d'\in \TZ*} \cfinalbis_{d,d'\!,\bull} * \sigma _{d'} (\uDA_{\bull} )
\end{align*}
with $\displaystyle \cfinalbis_{d,d',p,\ell } := \sum _{\substack{(p_{1} ,p_{2} ),n\in \mathbb{N}\\ p_{1} +p_{2} +2n=p}} \sum _{d_{1} \in \TZ*} \sum _{d_{2} \in \TN} \sum _{\substack{(\ell _{1} ,\ell _{2} )\in \mathbb{N}^2 \\ \ell _{1} +\ell _{2} =\ell }} \cmacro_{d,d_{1} ,p_{1} ,\ell _{1} }\cdot \sigma _{d_{1}} (\mathcal{S} ^{\infty } _{d_{2} ,n} )\cdot \cmicro_{d_{2} ,d'\!,p_{2} ,\ell _{2} }$.\\[1mm]
Moreover, for any $R\in \Rmaicro$, we have $\deg _{\varepsilon } \!R \in \mathbb{N}$ and $\deg _{\varepsilon } \!R \pm  \deg _{\mathcal{A}} \!R\in \mathbb{N}$ (see Figure~\ref{fig: degres des elements de Rmacro et Rmicro}). So by definition of the coefficients $\cmacro_{\dots }$ and $\cmicro_{\dots }$:
\begin{itemize}
	\item $\cmacro_{d,d_{1} ,p_{1} ,\ell _{1} } \neq 0 \Rightarrow  p_{1} +d\in \mathbb{N} \text{ and } p_{1} +d_{1} \in \mathbb{N}$,
	\item $\cmicro_{d_{2} ,d',p_{2} ,\ell _{2} } \neq 0 \Rightarrow  p_{2} -d_{2} \in \mathbb{N} \text{ and } p_{2} -d'\in \mathbb{N}$.
\end{itemize}
This implies that $\cfinalbis_{d,d',p,\ell } =\cfinal_{d,d',p,\ell }$ for any $(d,d',p,\ell )$.\\
Furthermore, those conditions on $(p_{1,} p_{2} ,d_{1} ,d_{2} )$ imply that, if $\cfinal_{d,d',p,\ell } \neq 0$, then:
\begin{itemize}
	\item $p=(p_{1} +d)+2n+(p_{2} -d_{2} )-d+d_{2} \in \mathbb{N}+\TN+\TN \subset \NN+ \frac{2\pi }{\Theta } $
	\item and $p+d-d'=(p_{1} +d)+2n+(p_{2} -d')\in \mathbb{N}$.
\end{itemize}
This and the property $(\sigma _{d'} (u_{p',\ell '} ) \neq 0 \Rightarrow  p'\in \NN \text{ and } p'+d'\in \NN)$ explain the sum indexes of the formula given for $\sigma _{d} (u_{p,\ell } )$ in Theorem~\ref{construction directe des u p,l}.
\end{demo}

\begin{exemple}
Using Theorem~\ref{construction directe des u p,l}, one can check that $u_{\pi /\Theta ,\ell }$ vanishes for any $\ell $. In addition $\frac{2\pi }{\Theta } > 1$, so the first non-zero far fields are $u_{0,0}$ and $u_{1,0}$ and they satisfy
\[\left\{\begin{array}{@{\;}r@{\;=\;}l@{\;\;}l@{}}
	\mu _{0} \Delta  u_{0,0} +\omega ^2 \rho _{0} u_{0,0} & \fSource &\text{ in }\Omega  \\[0.5mm]
	u_{0,0} & 0 &\text{ on }\bordgauche \\[0.5mm]
	u_{0,0} & 0 &\text{ on }\Gamma  \\[0.5mm]
	\sigma _{d} (u_{p,\ell } )& 0 & \ \forall d\in -\TN
\end{array}\right. \qquad  \text{ and } \quad  \left\{\begin{array}{@{\;}r@{\;=\;}l@{\;\;}l@{}}
	\mu _{0} \Delta  u_{1,0} +\omega ^2 \rho _{0} u_{1,0} & 0 &\text{ in }\Omega  \\[0.5mm]
	u_{1,0} & 0 &\text{ on }\bordgauche\\[0.5mm]
	u_{1,0} & \frac{\mu _{0} }{\mu _{1}} \partial _{y} u_{0,0} &\text{ on }\Gamma  \\[0.5mm]
	\sigma _{d} (u_{1,0} )& 0 & \ \forall d\in -\TN
\end{array}\right.\]\fbeq
\end{exemple}

\section{Error estimates}
\label{sec: estimations d'erreur}

\begin{textesansboite}
Let $\chi  \in \mathcal{C}^{\infty } (\mathbb{R}^2 )$ be equal to 1 on $B(0,1)$ and 0 outside $B(0,2)$, and, for any $\eta >0$, $\chi _{\eta } : x \mapsto  \chi ( \frac{x}{\eta } )$. We denote $\couche_{\varepsilon } := (-\varepsilon ,1)\times \mathbb{R}_+^*$ and $\Pi _{\varepsilon } :=\Omega \sqcup \Gamma \sqcup \couche_{\varepsilon }$ (defined similarly as $\Pi $ at page~\pageref{page: definition de Pi}). We define on $\Pi _{\varepsilon }$ the following variant of the far-and-layer fields $\bfu_{p,\ell }$ \[\bfue_{p,\ell } (x,y) := \left\{\begin{array}{ll}
	u_{p,\ell } (x,y) &\text{ in }\Omega  \\[0.5mm]
	U_{p,\ell } (x, \frac{y}{\varepsilon } )&\text{ in }\couche_{\varepsilon }
\end{array}\right.\]
We also denote, for any $p\in \NN$, $n_{p} :=\max \{\ell \in \mathbb{N} \mid u_{p,\ell } \neq 0 \text{ or } U_{p,\ell } \neq 0 \text{ or } S_{p,\ell } \neq 0\}$.\\

Let us define the approximate global field at order $\ptronc\in \mathbb{R}_+$ as follows for any $\varepsilon $ small enough
\[\forall (x,y)\in \Omega _{\varepsilon } , \quad  u_{\varepsilon ,\ptronc} (x,y) := (1-\chi _{\eta } (x,y)) \! \sum _{\substack{p\in \NN \cap [0,\ptronc]\\ \ell \in [\![0,n_{p} ]\!]}} \varepsilon ^{p} \ln ^{\ell } \!\varepsilon \, \bfue_{p,\ell } (x,y) \;+\; \chi _{\eta } (x,y) \! \sum _{\substack{p\in \NN \cap [0,\ptronc]\\ \ell \in [\![0,n_{p} ]\!]}} \varepsilon ^{p} \ln ^{\ell } \!\varepsilon \, S_{p,\ell } ({\textstyle \frac{x}{\varepsilon } , \frac{y}{\varepsilon } })\]
where $\eta :=\sqrt{\varepsilon }$. Note that $(1-\chi _{\eta } (x,y))\bfue_{p,\ell } (x,y)$ is well-defined on $\Omega _{\varepsilon }$ when $\eta >\varepsilon \Rcoin$ since $\Pi _{\varepsilon } \setminus B(0,\varepsilon \Rcoin) = \Omega _{\varepsilon } \setminus B(0,\varepsilon \Rcoin)$. We will see as a consequence of Theorem~\ref{estimation d'erreur globale} that $u_{\varepsilon ,\ptronc} \in H^{1} (\Omega _{\varepsilon } )$.\\

The matching zone is $\Omega _{\varepsilon } \cap \crown$ where $\crown$ is the annulus $\crown :=B(0,2\eta ) \setminus B(0,\eta )$. Letting $\eta =\sqrt{\varepsilon }$ makes the matching zone tend to 0 w.r.t. the far fields (because $\eta  \rightarrow 0$ when $\varepsilon \rightarrow 0$) and to infinity w.r.t. the corner fields (because $\frac{\eta }{\varepsilon } \rightarrow \infty $). Thanks to the matching assumption, we can state the following first error estimate concercing the error in the matching zone. We use the symbol $\lesssim $ for majorations valid up to a constant independent of $\varepsilon $.
\end{textesansboite}

\begin{lemme}{}
\label{erreur de raccord}%
Let $\ptronc \in \mathbb{R}_+$. For $\varepsilon $ small enough, we have
\[\bigg\| \sum _{p\in \NN\cap [0,\ptronc]} \sum _{\ell \in \mathbb{N}} \varepsilon ^{p} \ln ^{\ell } \!\varepsilon \,\big(\bfue_{p,\ell } (x,y)-S_{p,\ell } ({\textstyle \frac{x}{\varepsilon } , \frac{y}{\varepsilon } })\big)\bigg\|_{H^{1} (\Omega _{\varepsilon } \cap \crown)} \lesssim  \varepsilon ^{ \frac{\ptronc}{2} -1}.\]\fbeq[1.3]
\end{lemme}

\begin{demo}
We will compare $\bfue_{p,\ell }$ and $S_{p,\ell }$ in $\crown$ to their asymptotic expansions at 0, resp $\infty $. Let us denote:
\begin{itemize}
	\item $\bfu_{\varepsilon ,\ptronc} := \sum _{p\in \NN\cap [0,\ptronc]} \sum _{\ell \in \mathbb{N}} \varepsilon ^{p} \ln ^{\ell } \!\varepsilon \, \bfue_{p,\ell } $ in $\Pi _{\varepsilon }$
	\item $S_{\varepsilon ,\ptronc} := \sum _{p\in \NN\cap [0,\ptronc]} \sum _{\ell \in \mathbb{N}} \varepsilon ^{p} \ln ^{\ell } \!\varepsilon \, S_{p,\ell } ({\textstyle \frac{x}{\varepsilon } , \frac{y}{\varepsilon } })$ in $\Omega _{\varepsilon }$
	\item $\bfu_{\varepsilon ,\ptronc} ^{\mathcal{A}} := \sum _{p\in \NN\cap [0,\ptronc]} \sum _{\ell \in \mathbb{N}} \varepsilon ^{p} \ln ^{\ell } \!\varepsilon \cdot \left\{\begin{array}{@{}l@{}l}
		T_{\leqslant \ptronc-p} (\bfuDA_{p,\ell } )(x,y) &\text{ in }\Omega  \\[0.5mm]
		T_{\leqslant \ptronc-p} (\bfuDA_{p,\ell } )(x,{\textstyle \frac{y}{\varepsilon } }) &\text{ in }\couche_{\varepsilon } 
	\end{array}\right.$
	\item $S_{\varepsilon ,\ptronc} ^{\mathcal{A}} := \sum _{p\in \NN\cap [0,\ptronc]} \sum _{\ell \in \mathbb{N}} \varepsilon ^{p} \ln ^{\ell } \!\varepsilon \cdot T_{\geqslant p-\ptronc} (\SDA_{p,\ell } )({\textstyle \frac{x}{\varepsilon } , \frac{y}{\varepsilon } })$ in $\Pi _{\varepsilon }$.
\end{itemize}
(see Definition~\ref{def: troncature dans Abarre} for $T_{\bullet}$). We split the estimate into three parts that we will majorize separately:
\[\|\bfu_{\varepsilon ,\ptronc} -S_{\varepsilon ,\ptronc} \|_{H^{1} (\Omega _{\varepsilon } \cap \crown)} \leqslant  \|\bfu_{\varepsilon ,\ptronc} -\bfu_{\varepsilon ,\ptronc} ^{\mathcal{A}} \|_{H^{1} (\Omega _{\varepsilon } \cap \crown)} + \|\bfu_{\varepsilon ,\ptronc} ^{\mathcal{A}} -S_{\varepsilon ,\ptronc} ^{\mathcal{A}} \|_{H^{1} (\Omega _{\varepsilon } \cap \crown)} + \|S_{\varepsilon ,\ptronc} ^{\mathcal{A}} -S_{\varepsilon ,\ptronc} \|_{H^{1} (\Omega _{\varepsilon } \cap \crown)} .\]
In addition, we will split some of the norms $\|{\dots }\|_{H^{1} (\Omega _{\varepsilon } \cap \crown)}$ into $\|{\dots }\|_{H^{1} (\Omega \cap \crown)} +\|{\dots }\|_{H^{1} (\couche_{\varepsilon } \cap \crown)}$.
\begin{enumerate}
	\item \underline{$\|\bfu_{\varepsilon ,\ptronc} -\bfu_{\varepsilon ,\ptronc} ^{\mathcal{A}} \|_{H^{1} (\Omega \cap \crown)}$ :} By Proposition~\ref{confimation des DA des champs}, for any $(p,\ell )$ we have $u_{p,\ell } -T_{\leqslant \ptronc-p} (\bfuDA _{p,\ell } )=o_{\partial } (r^{\ptronc-p} )$ in $\Omega $ when $r\rightarrow 0$. By Definition~\ref{def: petit o derivable} of $o_{\partial }$, it implies $u_{p,\ell } -T_{\leqslant \ptronc-p} (\bfuDA_{p,\ell } )=O(r^{\ptronc-p} )=O(r^{\ptronc-p-1} )$ and $\nabla \big[u_{p,\ell } -T_{\leqslant \ptronc-p} (\bfuDA_{p,\ell } )\big]=O(r^{\ptronc-p-1} )$ uniformly in $\theta $. Thus:\vspace{-2mm} \begin{align*}
		\|\bfu_{\varepsilon ,\ptronc} -\bfu_{\varepsilon ,\ptronc} ^{\mathcal{A}} \|_{H^{1} (\Omega \cap \crown)}
		&\leqslant  \sum _{p\in \NN\cap [0,\ptronc]} \sum _{\ell =0} ^{n_{p}}  \varepsilon ^{p} \ln ^{\ell } \!\varepsilon \cdot \|u_{p,\ell } -T_{\leqslant \ptronc-p} (\bfuDA_{p,\ell } )\|_{H^{1} (\Omega \cap \crown)} \\
		&\lesssim  \sum _{p\in \NN\cap [0,\ptronc]} \sum _{\ell =0} ^{n_{p}}  \varepsilon ^{p} \ln ^{\ell } \!\varepsilon \cdot \left( \int _{B(0,2\eta )} r^{2(\ptronc-p-1)} \,\mathrm{d}r\,r\mathrm{d}\theta  \right)^{1/2} \\
		&\lesssim  \sum _{p\in \NN\cap [0,\ptronc]} \sum _{\ell =0} ^{n_{p}}  \varepsilon ^{p} \ln ^{\ell } \!\varepsilon  \cdot \eta ^{\ptronc-p} \\
		&\lesssim \; \eta ^{\ptronc-2} \qquad  \quad  \text{because }  \varepsilon ^{p} \lesssim  \eta ^{p} \text{ and } \ln ^{\ell } \!\varepsilon \lesssim \eta ^{-2}
	\end{align*}
	
	\item \underline{$\|\bfu_{\varepsilon ,\ptronc} -\bfu_{\varepsilon ,\ptronc} ^{\mathcal{A}} \|_{H^{1} (\couche_{\varepsilon } \cap \crown)}$ :} Similarly, $U_{p,\ell } -T_{\leqslant \ptronc-p} (\bfuDA_{p,\ell } )=o_{\partial } (x^{\ptronc-p} )$ in $\couche$ when $x\rightarrow \infty $. Thus
	\newlength{\widthGauche}%
	\settowidth{\widthGauche}{$\partial _{x} \big[U_{p,\ell } -T_{\leqslant \ptronc-p} (\bfuDA_{p,\ell } )\big]=O(x^{\ptronc-p-1} )$\,}%
	\begin{itemize}
		\item \parbox{\widthGauche}{$U_{p,\ell } -T_{\leqslant \ptronc-p} (\bfuDA_{p,\ell } )=O(x^{\ptronc-p} )$} so $U_{p,\ell } (x, \frac{y}{\varepsilon } ) -T_{\leqslant \ptronc-p} (\bfuDA_{p,\ell } )(x, \frac{y}{\varepsilon } )=O(x^{\ptronc-p} )$
		\item \parbox{\widthGauche}{$\partial _{x} \big[U_{p,\ell } -T_{\leqslant \ptronc-p} (\bfuDA_{p,\ell } )\big]=O(x^{\ptronc-p-1} )$} so $\partial _{x} \big[U_{p,\ell } (x, \frac{y}{\varepsilon } ) -T_{\leqslant \ptronc-p} (\bfuDA_{p,\ell } )(x, \frac{y}{\varepsilon } )\big] =O(\varepsilon ^{-1} x^{\ptronc-p-1} )$
		\item \parbox{\widthGauche}{$\partial _{Y} \big[U_{p,\ell } -T_{\leqslant \ptronc-p} (\bfuDA_{p,\ell } )\big]=O(x^{\ptronc-p} )$} so $\partial _{Y} \big[U_{p,\ell } (x, \frac{y}{\varepsilon } ) -T_{\leqslant \ptronc-p} (\bfuDA_{p,\ell } )(x, \frac{y}{\varepsilon } )\big]=O(\varepsilon ^{-1} x^{\ptronc-p} )$ 
	\end{itemize}
	which are all $O(\varepsilon ^{-1} x^{\ptronc-p-1} )$ (uniformly in $Y$). Hence:\vspace{-2mm} \begin{align*}
		\|\bfu_{\varepsilon ,\ptronc} -\bfu_{\varepsilon ,\ptronc} ^{\mathcal{A}} \|_{H^{1} (\couche_{\varepsilon } \cap \crown)}
		&\lesssim  \sum _{p\in \NN\cap [0,\ptronc]} \sum _{\ell =0} ^{n_{p}}  \varepsilon ^{p} \ln ^{\ell } \!\varepsilon \cdot \varepsilon ^{-1} \,\left( \int _{[0,2\eta ]\times [-\varepsilon ,0]} x^{2(\ptronc-p-1)} \,\mathrm{d}x\,\mathrm{d}y \right)^{1/2} \\
		&\lesssim  \sum _{p\in \NN\cap [0,\ptronc]} \sum _{\ell =0} ^{n_{p}}  \varepsilon ^{p} \ln ^{\ell } \!\varepsilon \cdot \varepsilon ^{-1} \cdot \varepsilon ^{ \frac{1}{2} } \eta ^{\ptronc-p- \frac{1}{2} } \\
		&\lesssim \; \eta ^{\ptronc-2} \qquad  \quad  \text{because }  \varepsilon =\eta ^2  \text{ and } \ln ^{\ell } \!\varepsilon \lesssim \eta ^{- \frac{1}{2} }
	\end{align*}
	
	\item \underline{$\|S_{\varepsilon ,\ptronc} -S_{\varepsilon ,\ptronc} ^{\mathcal{A}} \|_{H^{1} (\Omega \cap \crown)}$ :} Similarly $S_{p,\ell } -T_{\geqslant p-\ptronc} (\SDA _{p,\ell } )=o_{\partial } (r^{p-\ptronc} )$ in $\Omega $ when $r\rightarrow \infty $, so $S_{p,\ell } -T_{\geqslant p-\ptronc} (\SDA_{p,\ell } )=O(r^{p-\ptronc} ) \text{ and } \nabla \!_{(X,Y)}\big[S_{p,\ell } -T_{\geqslant p-\ptronc} (\SDA_{p,\ell } )\big]=O(r^{p-\ptronc-1} )=O(r^{p-\ptronc} )$ uniformly in $\theta $. Since $\nabla \!_{(x,y)} =\varepsilon ^{-1} \nabla \!_{(X,Y)}$, we deduce: \vspace{-2mm} \begin{align*}
		\|S_{\varepsilon ,\ptronc} -S_{\varepsilon ,\ptronc} ^{\mathcal{A}} \|_{H^{1} (\Omega \cap \crown)}
		&\leqslant  \sum _{p\in \NN\cap [0,\ptronc]} \sum _{\ell =0} ^{n_{p}}  \varepsilon ^{p} \ln ^{\ell } \!\varepsilon \cdot \Big\|\big[S_{p,\ell } -T_{\geqslant p-\ptronc} (\SDA_{p,\ell } )\big]({\textstyle \frac{x}{\varepsilon } , \frac{y}{\varepsilon } })\Big\|_{H^{1} (\Omega \cap \crown)} \\
		&\lesssim  \sum _{p\in \NN\cap [0,\ptronc]} \sum _{\ell =0} ^{n_{p}}  \varepsilon ^{p} \ln ^{\ell } \!\varepsilon \cdot  \varepsilon ^{-1} \left( \int _{\mathbb{R}^2  \setminus B(0,\eta )} \left( \frac{r}{\varepsilon } \right)^{2(p-\ptronc)} \,\mathrm{d}r\,r\mathrm{d}\theta  \right)^{1/2} \\
		&\lesssim  \sum _{p\in \NN\cap [0,\ptronc]} \sum _{\ell =0} ^{n_{p}}  \varepsilon ^{p} \ln ^{\ell } \!\varepsilon  \cdot \varepsilon ^{\ptronc-p-1} \cdot \eta ^{p-\ptronc+1} \\
		&\lesssim \; \eta ^{\ptronc-2} \qquad  \quad  \text{because }  \varepsilon =\eta ^2  \text{ and } \ln ^{\ell } \!\varepsilon \lesssim \eta ^{-1}
	\end{align*}
	
	\item \underline{$\|S_{\varepsilon ,\ptronc} -S_{\varepsilon ,\ptronc} ^{\mathcal{A}} \|_{H^{1} (\couche_{\varepsilon } \cap \crown)}$ :} $S_{p,\ell } -T_{\geqslant p-\ptronc} (\SDA _{p,\ell } )=o_{\partial } (X^{p-\ptronc} )$ in $\couche$ when $x\rightarrow 0$, so:
	\begin{itemize}
		\item $S_{p,\ell } -T_{\geqslant p-\ptronc} (\bfuDA_{p,\ell } )=O(X^{p-\ptronc} )$ 
		\item $\partial _{X} \big[S_{p,\ell } -T_{\geqslant p-\ptronc} (\SDA_{p,\ell } )\big]=O(X^{p-\ptronc-1} )=O(X^{p-\ptronc} )$
		\item $\partial _{Y} \big[S_{p,\ell } -T_{\geqslant p-\ptronc} (\SDA_{p,\ell } )\big]=O(X^{p-\ptronc} )$ 
	\end{itemize}
	uniformy in $Y$. Hence:\vspace{-2mm} \begin{align*}
		\|S_{\varepsilon ,\ptronc} -S_{\varepsilon ,\ptronc} ^{\mathcal{A}} \|_{H^{1} (\couche_{\varepsilon } \cap \crown)}
		&\leqslant  \sum _{p\in \NN\cap [0,\ptronc]} \sum _{\ell =0} ^{n_{p}}  \varepsilon ^{p} \ln ^{\ell } \!\varepsilon \, \Big\|\big[S_{p,\ell } -T_{\geqslant p-\ptronc} (\SDA_{p,\ell } )\big]({\textstyle \frac{x}{\varepsilon } , \frac{y}{\varepsilon } })\Big\|_{H^{1} (\couche_{\varepsilon } \cap \crown)} \\
		&\lesssim  \sum _{p\in \NN\cap [0,\ptronc]} \sum _{\ell =0} ^{n_{p}}  \varepsilon ^{p} \ln ^{\ell } \!\varepsilon \cdot \varepsilon ^{-1} \,\left( \int _{[\eta ,\infty [\times [-\varepsilon ,0]} \left( \frac{x}{\varepsilon } \right)^{2(p-\ptronc)} \,\mathrm{d}x\,\mathrm{d}y \right)^{1/2} \\
		&\lesssim  \sum _{p\in \NN\cap [0,\ptronc]} \sum _{\ell =0} ^{n_{p}}  \varepsilon ^{p} \ln ^{\ell } \!\varepsilon \cdot  \varepsilon ^{\ptronc-p-1} \cdot \varepsilon ^{ \frac{1}{2} } \eta ^{p-\ptronc+ \frac{1}{2} } \\
		&\lesssim \; \eta ^{\ptronc-2} \qquad  \quad  \text{because }  \varepsilon =\eta ^2  \text{ and } \ln ^{\ell } \!\varepsilon \lesssim \eta ^{- \frac{3}{2} }
	\end{align*}
	
	\item \underline{$\|\bfu_{\varepsilon ,\ptronc} ^{\mathcal{A}} -S_{\varepsilon ,\ptronc} ^{\mathcal{A}} \|_{H^{1} (\Omega _{\varepsilon } \cap \crown)}$ :} %
	\def\termeu{\tilde{\bfu}}%
	\def\termeS{\tilde{S}}%
	\def\Epsterme{\widetilde{\mathcal{H}}_\varepsilon}%
	Let us show that this norm vanishes. It suffices to prove that in $\Aemacro$: 
	\begin{equation}{\label{eq: exactitude du raccord algebrique tronque}}
	\fsum_{p\in \NN\cap [0,\ptronc]} \fsum_{\ell \in \mathbb{N}} \eps^{p} \ln ^{\ell } \!\eps \cdot T_{\leqslant \ptronc-p} (\bfuDA_{p,\ell } ) = \Eps^{-1} \bigg(\fsum_{p\in \NN\cap [0,\ptronc]} \fsum_{\ell \in \mathbb{N}} \eps^{p} \ln ^{\ell } \!\eps \cdot T_{\geqslant p-\ptronc} (\SDA_{p,\ell } ) \bigg) \end{equation}
	where $\Eps$ is defined in \eqref{eq: E}--\eqref{eq: E dans la couche}, and $\eps$ and $\ln \eps$ denote the algebraic indeterminates of $\Aemacro$ (we denote them differently from Sections~\ref{subsec: preliminaires algebriques} and \ref{subsec: raccordement algebrique} to avoid confusion with the real number $\varepsilon $). For any $(p,d)$ let $\termeu_{p,d}$ be the coordinate of $\fsum_{p,\ell } \eps^{p} \ln ^{\ell } \!\eps\, \bfuDA_{p,\ell }$ in $\mathcal{A}_{d} (\Pi )[\ln \eps]$ and $\termeS_{p,d}$ be the one of $\fsum_{p,\ell } \eps^{p} \ln ^{\ell } \!\eps\, \SDA_{p,\ell }$. Moreover, for any $d$ and $\varphi \in \mathcal{A}_{d} (\Pi )[\ln \eps]$, let $\Epsterme^{-1} (\varphi ):= \varepsilon ^{d} \Eps^{-1} (\varphi ) \in \mathcal{A}_{d} (\Pi )[\ln \eps]$. Then:\vspace{-2mm} \begin{align*}
		\text{\eqref{eq: exactitude du raccord algebrique tronque}}  &\Leftrightarrow  \fsum_{p\leqslant \ptronc} \fsum_{d\leqslant \ptronc-p} \eps^{p} \termeu_{p,d} = \Eps^{-1} \bigg(\fsum_{p\leqslant \ptronc} \fsum_{d\geqslant p-\ptronc} \eps^{p} \termeS_{p,d} \bigg)\\
		&\Leftrightarrow  \fsum_{p\leqslant \ptronc} \fsum_{p+d\leqslant \ptronc} \eps^{p} \termeu_{p,d} = \fsum_{p\leqslant \ptronc} \fsum_{p-d\leqslant \ptronc} \eps^{p-d} \Epsterme^{-1} ( \termeS_{p,d} )\\
		&\Leftrightarrow  \fsum_{p\leqslant \ptronc} \fsum_{p+d\leqslant \ptronc} \eps^{p} \termeu_{p,d} = \fsum_{p+d\leqslant \ptronc} \fsum_{p\leqslant \ptronc} \eps^{p} \Epsterme^{-1} ( \termeS_{p+d,d} )\\
		&\Leftrightarrow  \forall (p,d) \text{ s.t. }  p\leqslant \ptronc \text{ and } p+d\leqslant \ptronc, \quad  \termeu_{p,d} = \Epsterme^{-1} ( \termeS_{p+d,d} ) \stepcounter{equation}\tag{\theequation}\label{eq: exactitude du raccord algebrique tronque bis}
	\end{align*}
	because two formal series coincide iff their coordinates coincide one by one.\\
	But Proposition~\ref{confimation de l'hyp de raccord} shows the matching condition $ \fsum\limits_{p,d} \eps^{p} \termeu_{p,d} = \Eps^{-1} \Big(\fsum\limits_{p,d} \eps^{p} \termeS_{p,d} \Big) =\fsum\limits_{p,d} \eps^{p} \Epsterme^{-1} ( \termeS_{p+d,d} )$, which is equivalent to: $\forall (p,d),\;\termeu_{p,d} = \Epsterme^{-1} ( \termeS_{p+d,d} )$. Thus we get \eqref{eq: exactitude du raccord algebrique tronque bis}, and then \eqref{eq: exactitude du raccord algebrique tronque}.
\end{enumerate}\fbi
\end{demo}

\begin{theoreme}{global error estimate\titreligne}
\label{estimation d'erreur globale}%
For any $\ptronc \in \mathbb{R}_+$ we have\ $\big\|u_{\varepsilon } -u_{\varepsilon ,\ptronc} \big\|_{H^{1} (\Omega _{\varepsilon } )} = o(\varepsilon ^{ \frac{\ptronc}{2} -2} )$.\vspace{-1mm}
\end{theoreme}

\begin{demo}
Let $r_{\varepsilon ,\ptronc} := u_{\varepsilon ,\ptronc} -u_{\varepsilon }$. It satisfies, for some functions $f_{\varepsilon }$ and $g_{\varepsilon }$,
\[\left\{\begin{array}{r@{\;=\;}ll}
	\operatorname{div}(\mu _{\varepsilon } \nabla  r_{\varepsilon ,\ptronc} ) +\omega ^2 \rho _{\varepsilon } r_{\varepsilon ,\ptronc} & f_{\varepsilon } &\text{ in }\Omega _{\varepsilon } \\[0.5mm]
	r_{\varepsilon ,\ptronc} & 0 &\text{ on }\partial \Omega _{\varepsilon } \\[0.5mm]
	r_{\varepsilon ,\ptronc|y=0^{+}} - r_{\varepsilon ,\ptronc|y=0^{-}} & 0 &\text{ on }\Gamma \cap \{x>\eta \}\\[0.5mm]
	\mu _{0} \partial _{y} r_{\varepsilon ,\ptronc|y=0^{+}} - \mu _{1} \partial _{y} r_{\varepsilon ,\ptronc|y=0^{-}} & g_{\varepsilon } &\text{ on }\Gamma \cap \{x>\eta \}
\end{array}\right.\]
As for \eqref{eq: uniforme coercivite}, this problem is well-posed with a stability constant independent of $\varepsilon $:
\[\|r_{\varepsilon ,\ptronc} \|_{H^{1} (\Omega _{\varepsilon } )} \lesssim  \|f_{\varepsilon } \|_{L^2 (\Omega _{\varepsilon } )} + \|g_{\varepsilon } \|_{L^2 (\Gamma \cap \{x>\eta \}).}\]
To get our error estimate, it suffices to show that $\|f_{\varepsilon } \|_{L^2 } +\|g_{\varepsilon } \|_{L^2 } \lesssim  \varepsilon ^{ \frac{\ptronc}{2} -2} =\eta ^{\ptronc-4}$.

\medskip\underline{Estimate of $\|f_{\varepsilon } \|_{L^2 }$:} We denote:
\begin{itemize}
	\item $\bfu_{\varepsilon ,\ptronc} := \sum _{p\in \NN\cap [0,\ptronc]} \sum _{\ell \in \mathbb{N}} \varepsilon ^{p} \ln ^{\ell } \!\varepsilon \, \bfue_{p,\ell } $ in $\Pi _{\varepsilon }$
	\item $S_{\varepsilon ,\ptronc} := \sum _{p\in \NN\cap [0,\ptronc]} \sum _{\ell \in \mathbb{N}} \varepsilon ^{p} \ln ^{\ell } \!\varepsilon \, S_{p,\ell } ({\textstyle \frac{x}{\varepsilon } , \frac{y}{\varepsilon } })$ in $\Omega _{\varepsilon }$
	\item and $\mathcal{D}_{\varepsilon } :u\mapsto  \operatorname{div}(\mu _{\varepsilon } \nabla  u) +\omega ^2 \rho _{\varepsilon } u$ the differential operator of Helmholtz's equation.
\end{itemize}
Since $\mathcal{D}_{\varepsilon } u_{\varepsilon } =f$ by definition of $u_{\varepsilon }$, we have \begin{align*}
	f_{\varepsilon } &= \mathcal{D}_{\varepsilon } r_{\varepsilon ,\ptronc} \\
	&= \mathcal{D}_{\varepsilon } \big((1-\chi _{\eta } )\bfu_{\varepsilon ,\ptronc} + \chi _{\eta } S_{\varepsilon ,\ptronc} -u_{\varepsilon } \big) \\
	&= \big((1-\chi _{\eta } ) \mathcal{D}_{\varepsilon } \bfu_{\varepsilon ,\ptronc} -f\big)+ \chi _{\eta } \mathcal{D}_{\varepsilon } S_{\varepsilon ,\ptronc} +[\mathcal{D}_{\varepsilon } ,\chi _{\eta } ](S_{\varepsilon ,\ptronc} -\bfu_{\varepsilon ,\ptronc} )
\end{align*}
where $[.,.]$ is the commutator. Let us estimate these terms one by one.
\begin{enumerate}
	\item For $\varepsilon $ small enough, $1-\chi _{\eta }$ is equal to 0 in $B(0,\varepsilon \Rcoin)$ and 1 in $\operatorname{supp}(f)$, so $(1-\chi _{\eta } ) \mathcal{D}_{\varepsilon } \bfu_{\varepsilon ,\ptronc|\Omega } = (1-\chi _{\eta } )(\mu _{0} \Delta +\omega ^2 \rho _{0} ) \bfu_{\varepsilon ,\ptronc|\Omega } = (1-\chi _{\eta } )f=f$.\\
	Moreover, using that $\mu _{1} \partial _{Y} ^2  U_{p,\ell } = -(\mu _{1} \partial _{x} ^2 +\omega ^2 \rho _{1} )U_{p-2,\ell }$ for any $(p,\ell )$, we get in $\couche_{\varepsilon } \setminus B(0,\varepsilon \Rcoin)$:
	\begin{align*}
		\mathcal{D}_{\varepsilon } \bfu_{\varepsilon ,\ptronc} &= (\mu _{1} \Delta +\omega ^2 \rho _{1} )\bfu_{\varepsilon ,\ptronc} \\
		&= \sum _{p\in \NN\cap (\ptronc-2,\ptronc]} \sum _{\ell \in \mathbb{N}} \varepsilon ^{p} \ln ^{\ell } \!\varepsilon \,(\mu _{1} \partial _{x} ^2 +\omega ^2 \rho _{1} ) U_{p,\ell } \\
		&= \sum _{p\in \NN\cap (\ptronc-2,\ptronc]} \sum _{\ell \in \mathbb{N}} \varepsilon ^{p+ \frac{1}{2} } \ln ^{\ell } \!\varepsilon  \sum _{n\in \mathbb{N}} \mu _{1} \,(\partial _{x} ^2 +k_{1} ^2 )^{n+1} \partial _{y} u_{p-1-2n,\ell } (x,0) \cdot  \mathcal{U}_{n} ({\textstyle \frac{y}{\varepsilon } }) \tag*{by \eqref{eq: definition des champs de couche}}
	\end{align*}
	Since $\mathcal{U}_{n} \in L^{\infty } (-1,0)$ for any $n$ (see \eqref{eq: correcteurs de couche}), we deduce:
	\[\|(1-\chi _{\eta } ) \mathcal{D}_{\varepsilon } \bfu_{\varepsilon ,\ptronc} \|_{L^2 (\couche_{\varepsilon } )} \lesssim  \!\! \sum _{p\in \NN\cap (\ptronc-2,\ptronc]} \sum _{\ell =0} ^{n_{p}}  \sum _{n=0} ^{\lfloor {(p-1)/2}\rfloor }  \varepsilon ^{p+ \frac{1}{2} } \ln ^{\ell } \!\varepsilon  \cdot \big\|(\partial _{x} ^2 +k_{1} ^2 )^{n+1} \partial _{y} u_{p-1-2n,\ell } \big\|_{L^2 (\Gamma \cap \{x>\eta \})} .\]
	Let us show that, for any $p,\ell ,m$, $\|\partial _{x} ^{m} \partial _{y} u_{p,\ell |\Gamma } \|_{L^2 (\Gamma \cap \{x>\eta \})} \lesssim  \eta ^{-p-m-1}$ \setlabeleq\label{eq: estimation de derivees de u p,l sur Gamma}.
	\begin{itemize}
		\item The proof of Proposition~\ref{confimation des DA des champs} shows that $\partial _{x} ^{m} \partial _{y} u_{p,\ell } \in H^{1} _{0} (\Gamma )+\khiz \mathcal{A}(\Gamma )$ for any $p,\ell ,m$. So $\|\partial _{x} ^{m} \partial _{y} u_{p,\ell |\Gamma } \|_{L^2 (\Gamma \cap \{x>1\})} <\infty $.
		\item By Proposition~\ref{confimation des DA des champs}, $\partial _{y} u_{p,\ell |\Gamma } = \partial _{y} T_{\leqslant -p- \frac{1}{2} } (\bfuDA_{p,\ell } ) +o_{\partial } (x^{-p- \frac{3}{2} })$ when $x\rightarrow 0$, and Ansatz~\ref{ansatz sur Aeps} is satisfied, thus $T_{\leqslant -p- \frac{1}{2} } (\bfuDA_{p,\ell } )=0$. So: $\forall m\in \mathbb{N},\; \partial _{x} ^{m} \partial _{y} u_{p,\ell |\Gamma } =O(x^{-p- \frac{3}{2} -m})$ when $x\rightarrow 0$. It implies $\|\partial _{x} ^{m} \partial _{y} u_{p,\ell |\Gamma } \|_{L^2 (\Gamma \cap \{\eta <x<1\})} \lesssim  \eta ^{-p-m-1}$.
	\end{itemize}
	Hence:\vspace{-2mm} \begin{align*}
		\|(1-\chi _{\eta } ) \mathcal{D}_{\varepsilon } \bfu_{\varepsilon ,\ptronc} \|_{L^2 (\couche_{\varepsilon } )}
		&\lesssim  \sum _{p\in \NN\cap (\ptronc-2,\ptronc]} \sum _{\ell =0} ^{n_{p}}  \varepsilon ^{p+ \frac{1}{2} } \ln ^{\ell } \!\varepsilon  \sum _{n=0} ^{\lfloor {(p-1)/2}\rfloor }  \eta ^{-(p-1-2n)-2(n+1)-1} \\
		&\lesssim  \sum _{p\in \NN\cap (\ptronc-2,\ptronc]} \sum _{\ell =0} ^{n_{p}}  \varepsilon ^{p+ \frac{1}{2} } \ln ^{\ell } \!\varepsilon  \cdot \eta ^{-p-2} \\
		&\lesssim  \eta ^{\ptronc-4} \qquad  \quad  \text{because }  \ln ^{\ell } \!\varepsilon  \lesssim \eta ^{-1}
	\end{align*}
	So $\|(1-\chi _{\eta } )\mathcal{D}_{\varepsilon } \bfu_{\varepsilon ,\ptronc} -f\|_{L^2 (\Omega _{\varepsilon } )} \lesssim  \eta ^{\ptronc-4}$.
	
	\item For the second term of $f_{\varepsilon }$, \eqref{eq: definition des champs de coin} implies that $\mathcal{D}_{\varepsilon } S_{\varepsilon ,\ptronc} = \sum _{p\in \NN\cap (\ptronc-2,\ptronc]} \sum _{\ell \in \mathbb{N}} \varepsilon ^{p} \ln ^{\ell } \!\varepsilon \,\omega ^2 \rho  S_{p,\ell }$. But by Proposition~\ref{confimation des DA des champs} and Ansatz~\ref{ansatz sur Aeps}: $\forall (p,\ell ),\; S_{p,\ell } =O(r^{p} )$. Thus, using that $\rho $ is bounded and $\|S_{p,\ell } ( \frac{x}{\varepsilon } , \frac{y}{\varepsilon } ) \|_{L^2 (\Omega _{\varepsilon } \cap B(0,2\eta ))} = \varepsilon \, \|S_{p,\ell } \|_{L^2 (\Omega _{1} \cap B(0,2\eta /\varepsilon ))}$, we get \begin{align*}
		\|\chi _{\eta } \mathcal{D}_{\varepsilon } S_{\varepsilon ,\ptronc} \|_{L^2 (\Omega _{\varepsilon } )} &\lesssim  \sum _{p\in \NN\cap (\ptronc-2,\ptronc]} \sum _{\ell =0} ^{n_{p}}  \varepsilon ^{p} \ln ^{\ell } \!\varepsilon  \cdot \varepsilon \, \|S_{p,\ell } \|_{L^2 (\Omega _{1} \cap B(0,2\eta /\varepsilon ))} \\
		&\lesssim  \sum _{p\in \NN\cap (\ptronc-2,\ptronc]} \sum _{\ell =0} ^{n_{p}}  \varepsilon ^{p+1} \ln ^{\ell } \!\varepsilon \cdot \Big( \frac{\eta }{\varepsilon } \Big)^{p+1} \\
		&\lesssim  \eta ^{\ptronc-4} \qquad  \quad  \text{because }  \ln ^{\ell } \!\varepsilon  \lesssim \eta ^{-3}
	\end{align*}
	
	\item Using that $\|\Delta  \chi _{\eta } \|_{L^2 } +\|\nabla  \chi _{\eta } \|_{L^2 } =O(\eta ^{-2} )$, the last term of $f_{\varepsilon }$ satisfies \begin{align*}
		\|[\mathcal{D}_{\varepsilon } ,\chi _{\eta } ](S_{\varepsilon ,\ptronc} -\bfu_{\varepsilon ,\ptronc} )\|_{L^2 (\Omega _{\varepsilon } )}
		&\lesssim  \|\Delta  \chi _{\eta } \cdot (S_{\varepsilon ,\ptronc} -\bfu_{\varepsilon ,\ptronc} ) + 2 \nabla  \chi _{\eta } \cdot \nabla (S_{\varepsilon ,\ptronc} -\bfu_{\varepsilon ,\ptronc} )\|_{L^2 (\Omega _{\varepsilon } )} \\
		&\lesssim  \eta ^{-2} \|S_{\varepsilon ,\ptronc} -\bfu_{\varepsilon ,\ptronc} \|_{H^{1} (\Omega _{\varepsilon } \cap B(0,2\eta ) \setminus B(0,\eta ) )} \\
		&\lesssim  \eta ^{\ptronc-4} \qquad  \quad  \text{by Lemma~\ref{erreur de raccord}.} 
	\end{align*}\fbeq
\end{enumerate}

\underline{Estimate of $\|g_{\varepsilon } \|_{L^2 }$:} Without loss of generality, we assume $\partial _{y} \chi _{\eta |\Gamma } =0$. \eqref{eq: definition des champs de coin} implies that $\mu _{0} \partial _{Y} S_{p,\ell |Y=0^{+}} = \mu _{1} \partial _{Y} S_{p,\ell |Y=0^{-}}$ on $\Gamma \cap \{X>\Rcoin\}$ for any $(p,\ell )$. Thus we have on $\Gamma \cap \{x>\eta \}$: \begin{align*}
		g_{\varepsilon } &= \mu _{0} \partial _{y} r_{\varepsilon ,\ptronc|y=0^{+}} - \mu _{1} \partial _{y} r_{\varepsilon ,\ptronc|y=0^{-}} \\
		&= \chi _{\eta } \cdot \big(\mu _{0} \partial _{y} \bfu_{\varepsilon ,\ptronc|y=0^{+}} - \mu _{1} \partial _{y} \bfu_{\varepsilon ,\ptronc|y=0^{-}} \big) \\
		&= \chi _{\eta } \sum _{p\in \NN\cap (\ptronc-1,\ptronc]} \sum _{\ell \in \mathbb{N}} \varepsilon ^{p} \ln ^{\ell } \!\varepsilon \cdot \mu _{0} \partial _{y} u_{p,\ell |y=0^{+}} 
	\end{align*}
	By \eqref{eq: estimation de derivees de u p,l sur Gamma}, we have that $\|\partial _{y} u_{p,\ell |\Gamma } \|_{L^2 (\Gamma \cap \{x>\eta \})} \lesssim  \eta ^{-p-1}$. So $\|g_{\varepsilon } \|_{L^2 (\Gamma \cap \{x>\eta \})} \lesssim  \eta ^{\ptronc-4}$ since $\ln ^{\ell } \varepsilon \lesssim \eta ^{-2}$.
\end{demo}

\begin{remarque}
Theorem~\ref{estimation d'erreur globale} can be improved to $\|u_{\varepsilon } -u_{\varepsilon ,\ptronc} \|_{H^{1} (\Omega _{\varepsilon } )} = o(\varepsilon ^{\ptronc/2} )$. Indeed with the same notations as the above proof, one can show that, for any $(p,\ell )\in \NN\times \mathbb{N}$, $\|(1-\chi _{\eta } )\bfue_{p,\ell } \|_{H^{1} (\Omega _{\varepsilon } )} \lesssim  \eta ^{-p}$ and $\|\chi _{\eta } S_{p,\ell } \|_{H^{1} (\Omega _{\varepsilon } )} \lesssim  \eta ^{-p}$. Therefore:
\[\|r_{\varepsilon ,\ptronc} \|_{H^{1} (\Omega _{\varepsilon } )} \lesssim  \|r_{\varepsilon ,\ptronc+4} \|_{H^{1} (\Omega _{\varepsilon } )} + \!\! \sum _{p\in \NN\cap (\ptronc,\ptronc+4]} \sum _{\ell \in \mathbb{N}} \varepsilon ^{p} \ln ^{\ell } \!\varepsilon  \cdot \big(\|(1-\chi _{\eta } )\bfue_{p,\ell } \|_{H^{1} (\Omega _{\varepsilon } )} +\|\chi _{\eta } S_{p,\ell } \|_{H^{1} (\Omega _{\varepsilon } )}\big) = o(\varepsilon ^{ \frac{\ptronc}{2} }).\]\fbeq[1.5]
\end{remarque}

\phantomsection\label{demo du resultat principal}%
\begin{demo}[Proof of Theorem~\ref{resultat principal}]
It follows from Theorem~\ref{estimation d'erreur globale} applied at order $2\ptronc+4$ and from the fact that, for small enough $\varepsilon $, $\sum _{p\in \NN\cap [0,2\ptronc]} \sum _{\ell \in \mathbb{N}} \varepsilon ^{p} \ln ^{\ell } \!\varepsilon \, u_{p,\ell }$ and $u_{\varepsilon ,2\ptronc}$ coincide in $\Omega  \setminus B(0,\delta )$:
\begin{align*}
\bigg\|u_{\varepsilon } -\!\!\!\! \sum _{p\in \NN \cap [0,\ptronc]} \sum _{\ell \in \mathbb{N}} \!\varepsilon ^{p} \ln ^{\ell } \!\varepsilon \, u_{p,\ell } \bigg\|_{H^{1} (\Omega  \setminus B(0,\delta ))} &\lesssim  \|u_{\varepsilon } -u_{\varepsilon ,2\ptronc+4} \|_{H^{1} (\Omega _{\varepsilon } )} + \!\!\!\!\! \sum _{p\in \NN\cap (\ptronc,2\ptronc+4]} \sum _{\ell \in \mathbb{N}} \!\varepsilon ^{p} \ln ^{\ell } \!\varepsilon  \,\|u_{p,\ell } \|_{H^{1} (\Omega  \setminus B(0,\delta ))} \\[-2.5mm]
&= o(\varepsilon ^{\ptronc} )+ o(\varepsilon ^{\ptronc} ).
\end{align*}\fbeq[1.5]
\end{demo}

\begin{remarque}
An alternative way of stating the error estimate is: \[\bigg\|u_{\varepsilon } -\sum _{p\in \NN \cap [0,\red\ptronc)} \sum _{\ell \in \mathbb{N}} \varepsilon ^{p} \ln ^{\ell } \!\varepsilon \, u_{p,\ell } \bigg\|_{H^{1} (\Omega  \setminus B(0,\delta ))} = O(\varepsilon ^{\ptronc} \ln ^{n} \!\varepsilon  ) \qquad  \text{with}  \ n:=\max \{\ell \in \mathbb{N} \mid u_{\ptronc,\ell } \neq 0\}.\]\fbeq
\end{remarque}

\appendix
\section{Appendix: proof of Lemma~\ref{decomposition de A en sev}}
\label{sec: annexe A}

\begin{demonstration}%
\begin{textesansboite}
By Definition~\ref{def: espaces A}, we already know that the formulas to prove are true replacing $\bigoplus$ by $\sum $. So it suffices to show that those sums are direct. 
\begin{enumerate}
\item Let us show that $\mathcal{A}(\Omega ) = \bigoplus_{d\in \mathbb{R}} \mathcal{A}_{d} (\Omega )$. Let $n\in \mathbb{N}$ and, for any $j\in [\![1,n]\!]$, $d_{j} \in \mathbb{R}$ and $\varphi _{j} \in \mathcal{A}_{d_{j}} (\Omega )$. We assume that $d_{1} <d_{2} <\cdots <d_{n}$ and $\sum_{j=1}^n \varphi _{j} =0$. Let us show that $\forall j\in [\![1,n]\!],\; \varphi _{j} =0$.
Let $\theta \in (0,\Theta )$. By definition of $\mathcal{A}_{d_{j}} (\Omega )$, for any $j\in [\![1,n]\!]$, there is $P_{\theta ,j} \in \mathbb{C}[T]$ s.t.: $\forall r\in \mathbb{R}_+^*,\; \varphi _{j} (r,\theta )=r^{d_{j}} P_{\theta ,j} (\ln  r)$. Let us assume by contradiction that: $\exists j\in [\![1,n]\!],\; P_{\theta ,j} \neq 0$. Let $j_{0} :=\max \{j\in [\![1,n]\!] \mid P_{\theta ,j} \neq 0 \}$. Then $0 = \sum_{i=1}^n \lambda _{i} \varphi _{i} (r,\theta ) \sim  r^{d_{j_{0}}} P_{\theta ,j_{0}} (\ln r)$ when $r\rightarrow \infty $, so $P_{\theta ,j_{0}} =0$. This is contradictory, so: $\forall \theta \in (0,\Theta ),\forall j\in [\![1,n]\!],\forall r\in \mathbb{R}_+^*,\; \varphi _{j} (r,\theta )=0$.

\item By the same method, one can show that $\mathcal{A}(D) = \bigoplus_{d\in \mathbb{R}} \mathcal{A}_{d} (D)$ for any $D\in \{\Pi ,\Gamma ,\couche\}$.

\item Now we will prove \eqref{eq: decomposition de A en sev 2}. Let $d\in \mathbb{R}$ and $I:=\{(q,k)\in \mathbb{R}\times \mathbb{N} \mid q+k=d\}$. Separating real and imaginary parts, it is enough to show that:
\[\mathcal{A}_{d} (\Omega )\cap \mathcal{C}^{0} (\Omega ,\mathbb{R}) = \bigoplus_{(q,k)\in I} \big\{z\mapsto  \Im[(\alpha z)^{q} \,\overline{\alpha z}^{k} P(\log (\alpha z))] \ \big|\ P\in \mathbb{R}[T] \text{ and } \mathcal{P}(q,k,P)\big\}.\] 
Let, for any $(q,k)\in I$, $P_{q,k} \in \mathbb{R}[T]$ and $\varphi _{q,k} :z\mapsto  \Im[(\alpha z)^{q} \,\overline{\alpha z}^{k} P_{q,k} (\log (\alpha z))]$. We assume that $\mathcal{P}(q,k,P_{q,k} )$ holds for any $(q,k)\in I$, that the $P_{q,k}$ are all null except for a finite number, and that $\sum _{(q,k)\in I} \varphi _{q,k} =0$. Let us show by induction on $m:=\max _{(q,k)\in I} \deg  P_{q,k}$ that: $\forall (q,k)\in I,\; P_{q,k} =0$ (which implies in turn: $\forall (q,k),\; \varphi _{q,k} =0$). We initialize at $m=-\infty $, i.e. $(\forall (q,k),\; P_{q,k} =0)$, which is trivial. Thus, only the inductive step ($m\in \mathbb{N}$) remains to prove.\\
To do so, we first note that in $\Omega $:
\begin{equation}{\label{eq: preuve decomposition de A eq1}}
0 = r^{-d} \sum _{(q,k)\in I} \varphi _{q,k} (r,\theta ) = \sum  _{(q,k)\in I} \Im[e^{\mathrm{i}(q-k)(\theta -\Theta )} P_{q,k} (\ln r+\mathrm{i}(\theta -\Theta ))] \end{equation}
By applying $r\partial _{r}$, we deduce: $\displaystyle 0 = \sum  _{(q,k)\in I} \Im[e^{\mathrm{i}(q-k)(\theta -\Theta )} P_{q,k} ' (\ln r+\mathrm{i}(\theta -\Theta ))]$.\\
For any $(q,k)\in I$, let us define (note the switch of indexes at line 2): 
\[Q_{q,k} := \left\{\begin{array}{ll}
	P_{q,k} ' & \text{if } q\not\in \mathbb{N}\\[0.5mm]
	P_{q,k} '-P_{q,k} '(0) & \text{if } q\in \mathbb{N} \text{ and } q\leqslant k \\[0.5mm]
	P_{q,k} '- P_{k,q} '(0) & \text{if } q\in \mathbb{N} \text{ and } q>k
\end{array}\right.\]
Then $(\forall (q,k)\in I,\; \mathcal{P}(q,k,Q_{q,k} ))$. Since $\Im[(\alpha z)^{q} \,\overline{\alpha z}^{k} ] =-\Im[(\alpha z)^{k} \,\overline{\alpha z}^{q} ]$ (and it is null when $q=k$), the previous equality rewrites as
\[ 0 = \sum _{(q,k)\in I} \Im[(\alpha z)^{q} \,\overline{\alpha z}^{k} Q_{q,k} (\log (\alpha z))].\]
Then, by induction hypothesis: $\forall (q,k)\in I,\; Q_{q,k} =0$. This means for $P_{q,k}$ that:
\begin{itemize}
	\item If $q\not\in \mathbb{N}$, then $P_{q,k} '=0$.
	\item If $q\in \mathbb{N}$ and $q\neq k$, then $P_{q,k} ''=0$ and $P_{q,k} '(0)- P_{k,q} '(0)=0$ \ \setlabeleq\label{eq: preuve decomposition de A eq2}.
	\item If $q\in \mathbb{N}$ and $q=k$, then $P_{q,k} ''=0$.
\end{itemize}
So \eqref{eq: preuve decomposition de A eq1} reduces to : \[\begin{array}{r@{\;}r@{}l@{\hspace{-10ex}}l}
	0 &=\displaystyle \sum _{(q,k)\in I} & \Im[e^{\mathrm{i}(q-k)(\theta -\Theta )} P_{q,k} (0)] + \mathbf{1}_{q\in \mathbb{N}} \cdot  \Im[e^{\mathrm{i}(q-k)(\theta -\Theta )} P_{q,k} '(0)(\ln r+\mathrm{i}(\theta -\Theta ))] &\\[2.5mm]
	&=\displaystyle \sum _{(q,k)\in I} & \mathbf{1}_{(q\not\in \mathbb{N} \text{ or } q>k)} \cdot  \Im[e^{\mathrm{i}(q-k)(\theta -\Theta )} P_{q,k} (0)] & \text{by }  \mathcal{P}\\
	&+\ & \mathbf{1}_{(q\in \mathbb{N} \text{ and } q=k)} \cdot  \Im[e^{\mathrm{i}(q-k)(\theta -\Theta )} P_{q,k} '(0)(\ln r+\mathrm{i}(\theta -\Theta ))] &\\[2mm]
	&+\ & \mathbf{1}_{(q\in \mathbb{N} \text{ and } q>k)} \cdot  \Im[e^{\mathrm{i}(q-k)(\theta -\Theta )} P_{q,k} '(0)\, 2\mathrm{i}(\theta -\Theta )] & \text{by \eqref{eq: preuve decomposition de A eq2}}  \\[2.5mm]
	&=\displaystyle \sum _{(q,k)\in I} & \mathbf{1}_{(q\not\in \mathbb{N} \text{ or } q>k)} \cdot  P_{q,k} (0) \cdot \sin ((q-k)(\theta -\Theta )) &\\
	&+\ & \mathbf{1}_{(q\in \mathbb{N} \text{ and } q=k)} \cdot  P_{q,k} '(0)\cdot (\theta -\Theta ) &\\[2mm]
	&+\ & \mathbf{1}_{(q\in \mathbb{N} \text{ and } q>k)} \cdot  2P_{q,k} '(0) \cdot \sin \!\big((q-k)(\theta -\Theta )+ \frac{\pi }{2} \big)\cdot (\theta -\Theta ) &
\end{array}\] The functions of $\theta $ present here are linearly independent. So the coefficients $P_{q,k} (0)$ and $P_{q,k} '(0)$ are all zero. This concludes the proof.
\end{enumerate}\fbi
\end{textesansboite}
\end{demonstration}

\section{Appendix: asymptotic behaviors w.r.t. $\bfr$}
\label{sec: DA autour du coin}

In this section, we prove Theorems~\ref{DA des champs lointains concatenes} and \ref{DA des champs de coin}. We will use Sections~\ref{subsec: cadre fonctionnel des champs lointains} and \ref{subsec: cadre fonctionnel des champs de coin} which are after Theorems~\ref{DA des champs lointains concatenes} and \ref{DA des champs de coin} in the paper, but are independent of them.

\subsection{Proof of Theorem~\ref{DA des champs de coin}: asymptotic behavior for corner fields-like problems}
\label{subsec: DA des champs de coin}

\begin{figure}[h]
\begin{center}
\includegraphics{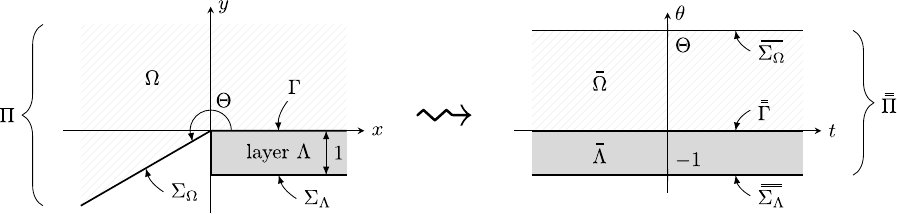}
\caption{The change of variables $(x,y) \rightsquigarrow (t,\theta )$ and the associated domains.}
\label{fig: sbd Pi}
\end{center}
\end{figure}

\begin{textesansboite}
The proof relies on the Kondrat'ev theory, usually used to analyse singularities of solutions of elliptic equations, see \cite{Kon67, KozMazRos97, MazNazPla00, CalCosDauVia06, Che12, Mak08}. We use it in a way that gives an expansion in $\Amacro(\Pi )$.\\

First we introduce the variables $(t,\theta )$, defined as $(\ln r,\theta )$ in $\Omega $ and $(\ln x,Y)$ in $\couche$. The pair $(t,\theta )$ lies in $\sbd \Pi  := \mathbb{R}\times (-1,\Theta )$. Moreover, we denote $\sbd \Omega $, $\sbd \couche$, $\sbd\bordgauche$, $\sbd \Gamma $ and $\sbd\borddroit$ the images of $\Omega $, $\couche$, $\bordgauche$, $\Gamma $ and $\borddroit$ by the change of variable $(x,y) \rightsquigarrow (t,\theta )$, see Figure~\ref{fig: sbd Pi}. The notation $\sbd{\raisebox{0pt}[1.5ex]{\dots }}$ is intended to remind the strip shape of $\sbd \Pi $. Finally, for any $u:\Pi \rightarrow \mathbb{C}$, we denote \[\sbd u:\, (t,\theta )\in \sbd \Pi  \;\mapsto \; \left\{\begin{array}{ll}
u(r=e^{t} ,\theta ) &\text{ in }\sbd \Omega \\[0.5mm]
u(x=e^{t} ,Y=\theta ) &\text{ in }\sbd \couche
\end{array}\right.\]\fbeq
\end{textesansboite}

\begin{definition}{Kondrat'ev spaces\titreligne}
\label{def: espaces de Kondratev}%
Let $H$ be a Hilbert space and $(s,\beta )\in \mathbb{R}^2 $. We define $K^{s} _{\beta } (\mathbb{R},H) := \{t\mapsto  e^{\beta t} u(t) \mid u\in H^{s} (\mathbb{R},H)\}$, equipped with the norm $\|u\|_{K^{s} _{\beta } (\mathbb{R},H)} := \|t\mapsto  e^{-\beta t} u(t)\|_{H^{s} (\mathbb{R},H)}$. We also denote $K^{\infty } _{\beta } (\mathbb{R},H) := \bigcap _{s\in \mathbb{R}} K^{s} _{\beta } (\mathbb{R},H)$, and $K^{s} _{\beta } (\sbd \Gamma ):=K^{s} _{\beta } (\mathbb{R},\mathbb{C})$ (identifying the line $\sbd \Gamma $ with $\mathbb{R}$).
\end{definition}

\begin{remarques}
\begin{itemize}
	\item If $s=m\in \mathbb{N}$, then $K^{m} _{\beta } (\mathbb{R},H) = \{u\in H^{m} _{\mathrm{loc}} (\mathbb{R},H) \mid \forall k\in [\![0,m]\!],\; e^{-\beta t} \partial _{t} ^{k} u(t) \in L^2 (\mathbb{R},H)\}$.
	\item Note that if $u\in K^{s} _{\beta } (\mathbb{R},H)$ and $u$ is zero in a vicinity of $-\infty $, then: $\forall \beta '\geqslant \beta ,\; u\in K^{s} _{\beta '} (\mathbb{R},H)$.
	\item Kondrat'ev's spaces are linked to $o_{\partial }$ (see Definition~\ref{def: petit o derivable}) in the following way. Let $\chi  \in \mathcal{C}^{\infty } (\mathbb{R})$ be equal to 0 in a vicinity of $-\infty $ and 1 in a vicinity of $+\infty $. Then for any $u:\Omega \rightarrow \mathbb{C}$:
	\[\left\{\begin{array}{l@{\quad\;}r@{\;\;\Rightarrow \;\;}l}
		\forall \beta \in \mathbb{R},&\chi \sbd{u}\in  \!\smash{\bigcap\limits _{m\in \mathbb{N}} } K^{\infty } _{\beta } (\mathbb{R},H^{m} (0,\Theta )) & \forall d>\beta ,\; u \! \underset{r\rightarrow \infty }{=} \! o_{\partial } (r^{d} )\\[0.5mm]
		\forall d\in \mathbb{R},&  u \! \underset{r\rightarrow \infty }{=} \! o_{\partial } (r^{d} ) & \exists t_{0} \in \mathbb{R},\,\forall \beta >d,\ \; \chi (\cdot +t_{0} )\,\sbd{u}\in  \!\smash{\bigcap\limits _{m\in \mathbb{N}} } K^{\infty } _{\beta } (\mathbb{R},H^{m} (0,\Theta ))
	\end{array}\right.\]
	And there are similar implications in $\couche$ and $\Gamma $.
\end{itemize}\fbi
\end{remarques}

\begin{textesansboite}
Let $H$ be a Hilbert space, $\beta \in \mathbb{R}$ and $\varphi \in K^{0} _{\beta } (\mathbb{R},H)$. For any $\lambda \in \mathbb{C}$ s.t. $\Re(\lambda )=\beta $ we define the (bilateral) Laplace transform of $\varphi $ at $\lambda $ as\vspace{-1.5mm}
\begin{equation}{\label{eq: transformee de Fourier-Laplace}}
\widehat{\varphi }(\lambda ) := \int _{-\infty } ^{+\infty } e^{-\lambda t} \varphi (t)\,\mathrm{d}t = \mathcal{F}\big[t\mapsto  e^{-\beta t} \varphi (t)\big](\Im(\lambda ))\end{equation}
where $\mathcal{F}$ denotes the Fourier transform. By properties of $\mathcal{F}$, we have $\widehat{\varphi } \in L^2 (\{\lambda \in \mathbb{C} \mid \Re(\lambda )=\beta \})$. If $\varphi $ depends on $t$ and $\theta $, $\widehat{\varphi }$ implies that we see $\varphi $ as a function from $\mathbb{R}$ to a space of functions of $\theta $.\\

To introduce the method, let us use the Laplace transform on $\sbd{S\,}$\!. Let $\khii$ be the truncation function introduced in Section~\ref{subsec: cadre fonctionnel des champs de coin} and $s:=\sbd{\khii S}$. Since $S\in \Vcoin+\khii \mathcal{A}(\Pi )$, there is $\beta \in \mathbb{R}$ s.t. $s\in K^{0} _{\beta } (\mathbb{R},L^2 (-1,\Theta ))$, so $\hat{s}(\lambda ,\theta )$ is well-defined for any $\lambda \in \mathbb{C}$ s.t. $\Re(\lambda )>\beta $. In addition we have for any $f\in L^2 (\Omega )$, $g\in L^2 (\couche )$, $h\in H^{-1/2} (\Gamma )$ and $u\in H^{1} _{\mathrm{loc}} (\Pi )$:
\begin{equation}{\label{eq: transfert du probleme de Poisson de Pi a sbdPi}}
\left\{\begin{array}{@{\,}r@{\;\;}l@{}}
	-\mu _{0} \Delta  u =f &\text{ in }\Omega  \\[0.5mm]
	-\mu _{1} \Delta  u =g &\text{ in }\couche \\[0.5mm]
	u_{|y=0^{+}} - u_{|y=0^{-}} = 0 &\text{ on }\Gamma  \\[0.5mm]
	\mu _{0} \partial _{y} u_{|y=0^{+}} - \mu _{1} \partial _{y} u_{|y=0^{-}} = h &\text{ on }\Gamma  \\[0.5mm]
	u = 0 &\text{ on }\bordgauche \cup \borddroit
\end{array}\right. \;\Leftrightarrow \; \left\{\begin{array}{@{\,}r@{\;\;}l@{}}
	-e^{-2t} \mu _{0} \Delta  \sbd u = \sbd f &\text{ in }\sbd \Omega  \\[0.5mm]
	-\mu _{1} \,(e^{-2t} (\partial _{t} ^2 -\partial _{t} )+\partial _{\theta } ^2 ) \sbd u =\sbd g &\text{ in }\sbd\couche \\[0.5mm]
	\sbd u_{|\theta =0^{+}} - \sbd u_{|\theta =0^{-}} = 0 &\text{ on }\sbd \Gamma  \\[0.5mm]
	e^{-t} \mu _{0} \partial _{\theta } \sbd u_{|\theta =0^{+}} - \mu _{1} \partial _{\theta } \sbd u_{|\theta =0^{-}} = \sbd h &\text{ on }\sbd \Gamma  \\[0.5mm]
	\sbd u = 0 &\text{ on }\sbd\bordgauche \cup \sbd\borddroit
\end{array}\right.\end{equation}
Taking $u:=\khii S$ in \eqref{eq: transfert du probleme de Poisson de Pi a sbdPi} and applying the Laplace transform yield that $\hat{s}$ satisfies:
\begin{equation}{\label{eq: pb verifie par la transformee de Laplace}}
\left\{\begin{array}{r@{\;=\;}ll}
	(\partial _{\theta } ^2 +\lambda ^2 ) \hat{s}(\lambda ,\theta ) & f & \text{if } \theta \in (0,\Theta )\\[0.5mm]
	\partial _{\theta } ^2  \hat{s}(\lambda ,\theta ) & g-\big[(\lambda +2)^2 -(\lambda +2)\big] \, \hat{s}(\lambda +2,\theta ) & \text{if } \theta \in (-1,0)\\[0.5mm]
	\hat{s}(\lambda ,0^{+} ) - \hat{s}(\lambda ,0^{-}) &0& \\[0.5mm]
	\mu _{1} \partial _{\theta } \hat{s}(\lambda ,0^{-}) & \mu _{0} \partial _{\theta } \hat{s}(\lambda +1,0^{+} ) & \\[0.5mm]
	\hat{s}(\lambda ,\Theta ) = \hat{s}(\lambda ,-1) & 0 & 
\end{array}\right.\end{equation}
for some functions functions $f$ and $g$ depending on $F$ and $\khii$. Solving this system w.r.t. $\theta $ allows us to extend $\hat{s}$ w.r.t. $\lambda $ further to the left in the complex plane, except at the $\lambda $ for which \eqref{eq: pb verifie par la transformee de Laplace} is ill-posed. These $\lambda $ are poles of $\hat{s}$ and they will be used in Proposition~\ref{etape de base du DA des champs de coin} to identify terms of the asymptotic expansion of $S$.\\

For any $m\in \mathbb{N}^*$, let us define the Hilbert spaces $\mathcal{H}^{m} := \{u\in H^{1} _{0} (-1,\Theta ) \mid u_{|(0,\Theta )} \in H^{m} (0,\Theta ) \text{ and } u_{|(-1,0)} \in H^{m} (-1,0)\}$ with the norm $\|u\|_{\mathcal{H}^{m}} ^2  := \|u\|_{H^{m} (-1,0)} ^2 + \|u\|_{H^{m} (0,\Theta )} ^2 $, and $\prodH{m}:=H^{m} (0,\Theta )\times H^{m} (-1,0)\times \mathbb{C}$. We need to solve problems of the following form with $u\in \mathcal{H}^{m+2}$ and $(f,g,a)\in \prodH{m}$, $m\in \mathbb{N}$:
\begin{equation}{\label{eq: EDO en theta a resoudre}}
\left\{\begin{array}{r@{\;=\;}ll}
	u''+\lambda ^2 u & f &\text{ on }(0,\Theta ) \\[0.5mm]
	u'' & g &\text{ on }(-1,0) \\[0.5mm]
	u(0^{+} ) - u(0^{-}) &0& \\[0.5mm]
	u'(0^{-}) & a & \\[0.5mm]
	u(\Theta ) = u(-1) & 0 & 
\end{array}\right. \end{equation}
Denoting $A(\lambda ):u \mapsto  \big((u''+\lambda ^2 u)_{|(0,\Theta )} ,\, u''_{|(-1,0)} ,\, u'(0^{-}) \big)$, \eqref{eq: EDO en theta a resoudre} is equivalent to $A(\lambda )(u) = (f,g,a)$. Note that for any $m\in \mathbb{N}$, $A(\lambda ) \in  \mathcal{L}(\mathcal{H}^{m+2} ,\prodH{m})$ where $\mathcal{L}({\dots })$ denotes the space of continuous linear maps between two normed vector spaces.
\end{textesansboite}

\newcommand{\pole}{q}
\begin{boitecarreesecable}{Lemma \nvnumpar{}:}
\label{resolution de l'EDO en theta}%
Let $m\in \mathbb{N}$.
\begin{enumerate}
	\item $A(\lambda ):\mathcal{H}^{m+2} \rightarrow \prodH{m}$ is invertible iff $\lambda \in  \mathbb{C} \setminus \TZ*$.
	\item $\lambda \mapsto A(\lambda )^{-1}$ is meromorphic from $\mathbb{C}$ to $\mathcal{L}(\prodH{m},\mathcal{H}^{m+2} )$. Its poles are in $\TZ*$ and are simple. Moreover, for any $\pole \in \TZ*$ and $v\in \prodH{m}$, the residue $(\Res_{\lambda =\pole} A(\lambda )^{-1} )(v)$ is proportional to the function $\theta  \in [-1,\Theta ] \mapsto  \mathbf{1}_{[0,\Theta ]} (\theta ) \cdot \sin (q\theta )$.
	\item Let $\beta _{1} < \beta _{2}$ be some reals. There is $C>0$ depending only on $(m,\beta _{1} ,\beta _{2} )$ s.t., for any $\lambda \in \mathbb{C}$ satisfying $\beta _{1} <\Re(\lambda )<\beta _{2}$ and $|\Im(\lambda )|>1$, we have: $\|A(\lambda )^{-1} \|_{\mathcal{L}(\prodH{m},\mathcal{H}^{m+2} )} \leqslant  C\, |\Im(\lambda )|^{m+2}$ \setlabeleq\label{eq: normes des operateurs resolvants R1, R2 et R3}.
\end{enumerate}
\end{boitecarreesecable}

\begin{demo}
\begin{enumerate}
	\item Let $\lambda \in \mathbb{C}$. An easy calculation gives that any element of $\operatorname{Ker} A(\lambda )$ must be proportional to $\theta  \in [-1,\Theta ] \mapsto  \mathbf{1}_{[0,\Theta ]} (\theta ) \cdot \sin (\lambda \theta )$. This function belongs to $\mathcal{H}^{m} \setminus \{0\}$ iff $\lambda \in \TZ*$. Therefore $\operatorname{Ker} A(\lambda ) \neq 0 \Leftrightarrow  \lambda \in \TZ*$.\\
	Moreover, it is easy to see that, for any $(f,g,a)\in \prodH{m}$, \eqref{eq: EDO en theta a resoudre} with $\lambda :=0$ has a unique solution in $\mathcal{H}^{m+2}$. Thus $A(0):\mathcal{H}^{m+2} \rightarrow \prodH{m}$ is invertible. Now, for any $\lambda \in \mathbb{C}$, we have $A(\lambda )=A(0)+\lambda ^2 B$ with $B:u\mapsto (u_{|(0,\Theta )} ,0,0)$. $B$ is a compact operator from $\mathcal{H}^{m+2}$ to $\prodH{m}$, so the Fredholm alternative holds for $A(\lambda )$. Therefore $A(\lambda ):\mathcal{H}^{m+2} \rightarrow \prodH{m}$ is invertible iff $\lambda \in \mathbb{C} \setminus \TZ*$.

	\item Since $\lambda \mapsto A(\lambda )$ is holomorphic on $\mathbb{C} \setminus \TZ*$, so too is $\lambda \mapsto A(\lambda )^{-1}$. Let us describe its behavior near the points $q\in \TZ*$ using \cite[Theorem~5.1.1, p.147]{KozMazRos97}. It depends on the ``Jordan chains'' of $A(q)$, which are the sequences $(u_{0} ,\dots ,u_{n} )\in (\mathcal{H}^{m} )^{n+1}$, $n\in \mathbb{N}$, s.t. 
	\begin{equation}{\label{eq: chaine de Jordan}}
	\forall k\in [\![0,n]\!],\qquad  \sum _{j=0} ^{k} \frac{1}{j!} \frac{\mathrm{d}^{j} A}{\mathrm{d}\lambda ^{j}} \Big|_{\lambda =q} (u_{k-j} )= 0 \end{equation}
	Let us compute these chains. Taking $k:=0$ in \eqref{eq: chaine de Jordan} gives $A(q)(u_{0} )=0$, so $u_{0}$ is proportional to $\theta \mapsto  \mathbf{1}_{[0,\Theta ]} (\theta ) \cdot \sin (q\theta )$ by step 1. In addition, if $n\geqslant 1$, taking $k:=1$ gives $A(q)(u_{1} )+ \frac{\mathrm{d}A}{\mathrm{d}\lambda } (q) (u_{0} )=0$. This implies on one hand $u_{1} '' +q^2 u_{1} +2qu_{0} =0$ in $(0,\Theta )$. On the other we get $u_{1|(-1,0)} '' =0$ and $u_{1} '(0^{-})=0$, so $u_{1|(-1,0)} =0$, which gives $u_{1} (0)=u_{1} (\Theta )=0$. Therefore
	\[0 \neq \int _{0} ^{\Theta } 2q|u_{0} |^2  = \int _{0} ^{\Theta } -(u_{1} ''+q^2 u_{1} ) \cdot  \overline{u_{0} }= \int _{0} ^{\Theta } -u_{1} \cdot (\underbrace{\overline{u_{0} ''+q^2 u_{0} }}_{=0} ) =0. \vspace{-2mm}\]
	This is absurd, so we must have $n=0$ for any Jordan chain of $A(q)$. Therefore, \cite[Theorem~5.1.1, p.147]{KozMazRos97} states that $A(\lambda )^{-1}$ has a simple pole at $q$, and $\operatorname{\operatorname{Im}} (\Res_{\lambda =q} A(\lambda )^{-1} )=\mathbb{C}u_{0}$.
	
	\item Let $(f,g,a)\in \prodH{m}$ and $u:=A(\lambda )^{-1} (f,g,a)\in \mathcal{H}^{m+2}$. We will write $\lesssim $ for inequalities valid up to a constant that depends on $(m,\beta _{1} ,\beta _{2} )$ but not on $(\lambda ,f,g,a)$. By Poincaré's inequality and integration by parts, we have\vspace{-2mm}
	\[\|u\|_{H^{1} (-1,0)} ^2  \lesssim  \int _{-1} ^{0} |u'|^2  = - \int _{0} ^{-1} g \bar{u}+ a\, \bar{u}(0)\lesssim  (\|g\|_{L^2 (-1,0)} +|a|) \cdot  \|u\|_{H^{1} (-1,0)}\]
	thus $\|u\|_{H^{1} (-1,0)} \lesssim  \|g\|_{L^2 (-1,0)} +|a|$. Since $u''_{|(-1,0)} =g$, we deduce $\|u\|_{H^{m+2} (-1,0)} \lesssim  \|g\|_{H^{m}} +|a|$.\\
	Then, let $u_{1} :\theta  \in [0,\Theta ] \mapsto u(0) \big(1- \frac{\theta }{\Theta } \big)$, $v:=u-u_{1}$ and $f_{\lambda } :=f-\lambda ^2 u_{1}$. We have $v(0)=v(\Theta )=0$ and $v''+\lambda ^2 v=f_{\lambda }$, so
	\[\int _{0} ^{\Theta } (-|v'|^2 +\lambda ^2 |v|^2 )=\int _{0} ^{\Theta } f_{\lambda } \bar{v}\]
	Dividing by $\lambda $ and taking the absolute value of the imaginary part, we get: $\frac{|\Im(\lambda )|}{|\lambda |^2 } \|v'\|_{L^2 } ^2  + |\Im(\lambda )| \|v\|_{L^2 } ^2  \leqslant  \frac{1}{|\lambda |} \|f_{\lambda } \|_{L^2 } \|v\|_{L^2 }$. Now we asume $|\Im(\lambda )|>1$ and $\beta _{1} <\Re(\lambda )<\beta _{2}$, so $|\lambda |\lesssim |\Im(\lambda )|$. Thus $\|v'\|_{L^2 } ^2  + |\lambda |^2  \|v\|_{L^2 } ^2  \leqslant  \|f_{\lambda } \|_{L^2 } \|v\|_{L^2 }$, which yields $\|v\|_{L^2 } \lesssim  |\lambda |^{-2} \|f_{\lambda } \|_{L^2 }$ and then $\|v'\|_{L^2 } \lesssim  |\lambda |^{-1} \|f_{\lambda } \|_{L^2 }$. Now from $v''=f_{\lambda } -\lambda ^2 v$, one can easily derive by induction on $m$ that $\|v\|_{H^{m+2}}  \lesssim  |\lambda |^{m} \|f_{\lambda } \|_{H^{m}}$. This implies $\|u\|_{H^{m+2} (0,\Theta )} \lesssim  |\lambda |^{m} \|f\|_{H^{m}} +|\lambda |^{m+2} |u(0)|$, with $|u(0)|\lesssim  \|u\|_{H^{m+2} (-1,0)}$.\\
	So finally $\|u\|_{\mathcal{H}^{m+2}} \lesssim  |\Im(\lambda )|^{m+2} \|(f,g,a)\|_{\prodH{m}}$.
\end{enumerate}\fbi
\end{demo}

\begin{textesansboite}
Proposition~\ref{etape de base du DA des champs de coin} is the base step of the asymptotic expansion. It will be applied to $s$ the rest of the asymptotic expansion at a given order, and the function $s_{\mathrm{r}}$ below will be the rest at the next order. Iterating this process provides an asymptotic expansion of $S$ at any order. Since Theorem~\ref{DA des champs de coin} implies infinite regularity on $S$, we work in the spaces $K^{\infty } _{\beta }$ and $\mathcal{H}^{m}$ for any $m$.
\end{textesansboite}

\begin{proposition}{}
\label{etape de base du DA des champs de coin}%
Let $d\in \mathbb{R}$, $m\in \mathbb{N} \setminus \{0,1\}$, $s\in  \bigcap\limits _{\beta >d} K^{\infty } _{\beta } (\mathbb{R},\mathcal{H}^{m} )$, $f\in  \!\bigcap\limits _{\beta >d-1} \! K^{\infty } _{\beta } (\mathbb{R},H^{m-2} (0,\Theta ))$, $g\in  \!\bigcap\limits _{\beta >d-1} \! K^{\infty } _{\beta } (\mathbb{R},H^{m-2} (-1,0))$, $h\in  \!\bigcap\limits _{\beta >d-1} \! K^{\infty } _{\beta } (\sbd \Gamma )$ and $d_{\mathrm{\min }} := \min (\TZ*\cap (d-1,d])$. We assume:\vspace{-2mm}
\[\left\{\begin{array}{rl}
	\mu _{0} \Delta  s = f  &\text{ in }\sbd \Omega \\[0.5mm]
	\mu _{1} \,(e^{-2t} (\partial _{t} ^2 -\partial _{t} )+\partial _{\theta } ^2 )s = g  &\text{ in }\sbd\couche \\[0.5mm]
	s_{|\theta =0^{+}} - s_{|\theta =0^{-}} =0  &\text{ on }\sbd \Gamma  \\[0.5mm]
	e^{-t} \mu _{0} \partial _{\theta } s_{|\theta =0^{+}} - \mu _{1} \partial _{\theta } s_{|\theta =0^{-}} = h &\text{ on }\sbd \Gamma  \\[0.5mm]
	s =0 &\text{ on }\sbd\bordgauche \cup \sbd\borddroit
\end{array}\right.\]
Then there is $s_{\mathrm{r}} \in  \!\bigcap\limits _{d-1<\beta <d_{\mathrm{\min }}}  \! K^{\infty } _{\beta } (\mathbb{R},\mathcal{H}^{m} )$ and $(c_{q} )\in \mathbb{C}^{\TZ*\cap  (d-1,d]}$ s.t. $\displaystyle s = \sum _{q\in \TZ*\cap  (d-1,d]} c_{q} \,\sbd{\phi_{q} } + s_{\mathrm{r}}$.\fbeq[0.5]
\end{proposition}

\newcommand{\poles}{\mathfrak{P}}
\begin{demo}
For any $\beta \in \mathbb{R}$, let $\mathbb{C}_{\Re>\beta } :=\{\lambda \in \mathbb{C} \mid \Re(\lambda )>\beta \}$. The Laplace transform of $s$ is well-defined and holomorphic from $\mathbb{C}_{\Re>d}$ to $\mathcal{H}^{m}$. We will show that $\hat{s}$ has a meromorphic extension on $\mathbb{C}_{\Re>d-1}$, with poles belonging in $\poles:=\TZ*\cap (d-1,d]$, then we will apply the residue theorem on a rectangle surrounding these poles. The functions $\phi_{q}$ will appear in the residues.\medskip

\underline{Step 1 :} $\hat{s}$ satisfies in $\mathbb{C}_{\Re>d}$:
\[\left\{\begin{array}{r@{}ll}
	\mu _{0} (\partial _{\theta } ^2 +\lambda ^2 ) \hat{s}(\lambda ,\theta ) &\;=\; \hat{f}(\lambda ,\theta ) & \text{if}  \ \theta \in (0,\Theta )\\[0.5mm]
	\mu _{1} \partial _{\theta } ^2  \hat{s}(\lambda ,\theta ) &\;=\; \hat{g}(\lambda ,\theta ) -\mu _{1} \,\big[(\lambda +2)^2 -(\lambda +2)\big] \, \hat{s}(\lambda +2,\theta ) := \tilde{g}(\lambda ,\theta ) & \text{if}  \ \theta \in (-1,0)\\[0.5mm]
	\hat{s}(\lambda ,0^{+} ) - \hat{s}(\lambda ,0^{-}) &\;=\; 0& \\[0.5mm]
	\mu _{1} \partial _{\theta } \hat{s}(\lambda ,0^{-}) &\;=\; \hat{h}(\lambda ) +\mu _{0} \partial _{\theta } \hat{s}(\lambda +1,0^{+} ) := \tilde{h}(\lambda ) & \\[0.5mm]
	\hat{s}(\lambda ,\Theta ) = \hat{s}(\lambda ,-1) &\;=\; 0 &
\end{array}\right.\] For any $\lambda \in  \mathbb{C}_{\Re>d-1} \setminus \TZ*$, let $\tilde{s} (\lambda ) := A(\lambda )^{-1} \big( \frac{1}{\mu _{0}} \hat{f}(\lambda ), \frac{1}{\mu _{1}} \tilde{g}(\lambda ), \frac{1}{\mu _{1}} \tilde{h}(\lambda )\big)$. Then $\tilde{s}$ is also solution of the above system, so $\tilde{s}$ and $\hat{s}$ coincide on $\mathbb{C}_{\Re>d} \setminus \TZ*$ by Lemma~\ref{resolution de l'EDO en theta}. Hence $\tilde{s}$ is an extension of $\hat{s}$ on $\mathbb{C}_{\Re>d-1} \setminus \poles$, that we will still denote $\hat{s}$. Lemma~\ref{resolution de l'EDO en theta} implies that this extension is meromorphic with simple poles and: $\forall q\in \poles,\; \Res_{q} \hat{s} =(\Res_{\lambda =q} A(\lambda )^{-1} ) \big( \frac{1}{\mu _{0}} \hat{f}(q), \frac{1}{\mu _{1}} \tilde{g}(q), \frac{1}{\mu _{1}} \tilde{h}(q)\big)$.\medskip\smallskip

\underline{Step 2 :} Let $t\in \mathbb{R}$, $b\in (d-1,d_{\mathrm{\min }} )$ and $k\in \mathbb{R}_+^*$. The residue theorem applied to $\lambda  \mapsto  e^{\lambda t} \hat{s}(\lambda )$ on the rectangle $[b,d+1]\times [-k,k] \subset \mathbb{C}$ gives the following equality in $\mathcal{H}^{m}$. \begin{eqnarray*}
	2\mathrm{i}\pi  \sum _{q\in \poles} \Res_{q} (\lambda \mapsto e^{\lambda t} \hat{s}(\lambda )) &=& \int _{-k} ^{k} e^{(d+1+\mathrm{i}\gamma )t} \hat{s}(d+1+\mathrm{i}\gamma )\,\mathrm{d}\gamma  - \int _{-k} ^{k} e^{(b+\mathrm{i}\gamma )t} \hat{s}(b+\mathrm{i}\gamma )\,\mathrm{d}\gamma  \\[-1mm]
	&& + \int _{b} ^{d+1} e^{(\beta -\mathrm{i}k)t} \hat{s}(\beta -\mathrm{i}k)\,\mathrm{d}\beta  - \int _{b} ^{d+1} e^{(\beta +\mathrm{i}k)t} \hat{s}(\beta +\mathrm{i}k)\,\mathrm{d}\beta 
\end{eqnarray*}
Note that when $k\rightarrow \infty $ the first integral tends to the inverse Laplace transform of $s$ (up to a constant). We will show that the last two integrals tend to 0, by proving the following: 
\begin{equation}{\label{eq: etape de base du DA des champs de coin preuve}}
\forall n\in \mathbb{N}, \qquad  \sup _{\beta \in [b,d+1]} \int _{\mathbb{R} \setminus [-1,1]} \gamma ^{n} \| \hat{s}(\beta +\mathrm{i}\gamma )\|_{\mathcal{H}^{m}} \,\mathrm{d}\gamma  < \infty  \vspace{-2mm}\end{equation}
\begin{sousdemo}
	Let $n\in \mathbb{N}$ and $\beta \in [b,d+1]$. \eqref{eq: normes des operateurs resolvants R1, R2 et R3} and the expressions of $\tilde{g}$ and $\tilde{h}$ yield
	\begin{align*}
		\!\int_{\mathbb{R} \setminus [-1,1]} \!\! \gamma ^{n} \| \hat{s}(\beta +\mathrm{i}\gamma )\|_{\mathcal{H}^{m}} \mathrm{d}\gamma 
		&\leqslant  \int_{\mathbb{R} \setminus [-1,1]} \!\! \gamma ^{n} \|A(\lambda )^{-1} \|_{\mathcal{L}(\prodH{m-2},\mathcal{H}^{m} )}  \Big\|\big({\textstyle  \frac{1}{\mu _{0}} \hat{f}(\beta +\mathrm{i}\gamma ), \frac{1}{\mu _{1}} \tilde{g}(\beta +\mathrm{i}\gamma ), \frac{1}{\mu _{1}} \tilde{h}(\beta +\mathrm{i}\gamma )}\big)\Big\|_{\prodH{m-2}} \mathrm{d}\gamma  \hspace{-20mm} \\
		&\lesssim  \int_{\mathbb{R} \setminus [-1,1]} \! \gamma ^{n+m+2} \,\Big( \| \hat{f} (\beta +\mathrm{i}\gamma )\|_{H^{m-2}} +\| \hat{g} (\beta +\mathrm{i}\gamma )\|_{H^{m-2}} +| \hat{h}(\beta +\mathrm{i}\gamma )| \\[-2mm]
		& \qquad  \qquad  \quad  \ +\; \| \hat{s}(\beta +\mathrm{i}\gamma +1)\|_{\mathcal{H}^{m}} + \gamma ^2 \| \hat{s}(\beta +\mathrm{i}\gamma +2)\|_{\mathcal{H}^{m}} \Big) \,\mathrm{d}\gamma  
	\end{align*}
	Let us treat the term with $f$, given that the others are similar.
	\begin{eqnarray*}
		\int _{\mathbb{R} \setminus [-1,1]} \gamma ^{n+m+2} \| \hat{f} (\beta +\mathrm{i}\gamma )\|_{H^{m-2}} \,\mathrm{d}\gamma  
		&\underset{\mathcal{\text{C.S.}}}{\leqslant }& \|\gamma ^{-1} \|_{L^2 (\mathbb{R} \setminus [-1,1])} \cdot  \Big\| \gamma ^{n+m+3} \| \hat{f}(\beta +\mathrm{i}\gamma )\|_{H^{m-2}} \Big\|_{L^2 (\mathbb{R} \setminus [-1,1])} \\[-0.5mm]
		&\lesssim & \|e^{-\beta t} \partial _{t} ^{n+m+3} f\|_{L^2 (\mathbb{R},H^{m-2} (0,\Theta ))} \\
		&\lesssim & \max _{\beta '\in \{b,d+1\}} \|e^{-\beta 't} \partial _{t} ^{n+m+3} f\|_{L^2 (\mathbb{R},H^{m-2} (0,\Theta ))} \\
		&=& \max _{\beta '\in \{b,d+1\}} \|f\|_{K^{n+m+3} _{\beta '} (\mathbb{R},H^{m-2} (0,\Theta ))}
	\end{eqnarray*}
	thanks to an interpolation between $b$ and $d+1$. We have majorized by a finite constant independent of $\beta $, so the supremum on $\beta $ is finite. Thus \eqref{eq: etape de base du DA des champs de coin preuve} is proven.
\end{sousdemo}
Therefore: $\displaystyle \int _{1} ^{\infty } \bigg\| \int _{b} ^{d+1} e^{(\beta \pm \mathrm{i}k)t} \hat{s}(\beta \pm \mathrm{i}k) \,\mathrm{d}\beta \,\bigg\|_{\mathcal{H}^{m}} \mathrm{d}k \lesssim  \int _{1} ^{\infty } \! \int _{b} ^{d+1} \| \hat{s}(\beta \pm \mathrm{i}k)\|_{\mathcal{H}^{m}} \,\mathrm{d}\beta \,\mathrm{d}k < +\infty $.\\
So by taking arbitrarily large $k$, we can make the terms $\int _{b} ^{d+1} e^{(\beta \pm \mathrm{i}\gamma )t} \hat{s}(\beta \pm \mathrm{i}k) \,\mathrm{d}\beta $ tend to 0 in $\mathcal{H}^{m}$. Hence:
\[\sum _{q\in \poles} \Res_{q} (\lambda \mapsto e^{\lambda t} \hat{s}(\lambda )) = \frac{1}{2\mathrm{i}\pi } \int _{\mathbb{R}} e^{(d+1+\mathrm{i}\gamma )t} \hat{s}(d+1+\mathrm{i}\gamma ) \,\mathrm{d}\gamma   - \frac{1}{2\mathrm{i}\pi } \int _{\mathbb{R}} e^{(b+\mathrm{i}\gamma )t} \hat{s}(b+\mathrm{i}\gamma ) \,\mathrm{d}\gamma .\]
Now, since $\hat{s}$ has only simple poles, the terms of the sum are equal to 
\[e^{qt} \Res_{q} \hat{s} = e^{qt} (\Res_{\lambda =q} A(\lambda )^{-1} )\big({\textstyle \frac{1}{\mu _{0}} \hat{f}(q), \frac{1}{\mu _{1}} \tilde{g}(q), \frac{1}{\mu _{1}} \tilde{h}(q)}\big) = e^{qt} c_{q} \mathbf{1}_{[0,\Theta ]} (\theta ) \sin (d\theta ) = c_{q} \,\sbd{\phi_{q} }(t,\theta )\]
for some constant $c_{q} \in \mathbb{C}$ by point 2 of Lemma~\ref{resolution de l'EDO en theta}. The first intergral is equal to $s(t)$ by inverse Laplace transform (because $s\in K_{d+1} ^{\infty } (\mathbb{R},\mathcal{H}^{m} )$). We define $s_{\mathrm{r}} (t)$ to be equal to the last integral. Thus we get the desired formula, and \eqref{eq: etape de base du DA des champs de coin preuve} shows that $s_{\mathrm{r}} \in K^{\infty } _{b} (\mathbb{R},\mathcal{H}^{m} )$.
\end{demo}

\begin{textesansboite}
Let $S$ be the function set in Theorem~\ref{DA des champs de coin} and $(\sigma _{d} (S))_{d\in \TN}$ the coefficients of Definition~\ref{def: sigma dans V + khi A}, which vanish for $d$ big enough. With the $\sigma _{d} (S)$, the non-variational part of $S$ (denoted $S^{\mathcal{A}}$ below) can be explicitly computed, so it remains to get an asymptotic expansion of its variational part (denoted $S^{\Vcoin}$ below). In order to apply Proposition~\ref{etape de base du DA des champs de coin}, we need that $\sbd{\khii S^{\Vcoin} }$ belongs to some space $K^{\infty } _{\beta }$ (since $\Omega _{1} \setminus B(0,\Rcoin)=\Pi  \setminus B(0,\Rcoin)$, $\khii S^{\Vcoin}$ can be seen as a function defined on $\Pi $, which allows us to consider $\sbd{\khii S^{\Vcoin} }$). We will use again the notations $\langle .\rangle $ and $T_{\geqslant d}$ of Definitions~\ref{def: serie geometrique d'operateurs 1D} and \ref{def: troncature dans Abarre}.
\end{textesansboite}

\begin{proposition}{}
\label{regularite de S variationnel}%
Let $d\in \mathbb{R}_-$,\vspace{-1mm}
\[S^{\mathcal{A}} := T_{\geqslant d} \Bigg[ \left\langle  -\RcouchePi\circ \partial _{X|\couche} ^2  ,\; \frac{\mu _{0} }{\mu _{1}} \RGammaPi\circ \partial _{Y|\Gamma ,Y=0^{+}} \right\rangle  \bigg( \frac{1}{\mu _{0}} \ROmegaPi (F^{\infty } _{\Omega } )+ \frac{1}{\mu _{1}} \RcouchePi (F^{\infty } _{\couche} ) + \fsum_{d\in \TN} \sigma _{d} (S) \,\phi_{d} \bigg) \Bigg] \vspace{-2.5mm}\]
in $\mathcal{A}(\Pi )$ and $S^{\Vcoin} := S-\khii S^{\mathcal{A}}$. There is $d\in \mathbb{R}_-$ s.t. $S^{\Vcoin} \in  \Vcoin$ and: $\forall m\in \mathbb{N},\;\sbd{\khii S^{\Vcoin} } \in K^{\infty } _{1/2} (\mathbb{R},\mathcal{H}^{m} )$.
\end{proposition}

\begin{demo}
We consider any $d\in \mathbb{R}_-$ and we will fix it later. The proof has three steps: showing that $S^{\Vcoin} \in \Vcoin$, showing that $S^{\Vcoin}$ is regular w.r.t $t$, and deducing that it is regular w.r.t. $\theta $. For any $r_{1} ,r_{2} \in \mathbb{R}$ s.t. $\Rcoin<r_{1} <r_{2}$, we will denote $\chi \in \mathcal{C}^{\infty } (\mathbb{R}^2 )$ a radial function equal to 0 on $B(0,r_{1} )$ and 1 on $\mathbb{R}^2  \setminus B(0,r_{2} )$ ($r_{1} ,r_{2} $ are implicit in this notation).\medskip\smallskip

\underline{Step 1:} Let $\varphi :\Pi \rightarrow \mathbb{C}$ be equal to $T_{\leqslant d} (F^{\infty } _{D} )$ on $D$ for any $D\in \{\Omega ,\couche\}$, and $F_{\mathrm{v}} = F-\khii \varphi $. If $d<-2$ then $(1+r)F_{\mathrm{v}} \in L^2 (\Omega _{1} )$. So the same construction than in Theorem~\ref{cadre fonctionnel pour les champs de coin} shows that for $d$ small enough, the problem 
\begin{equation}{\label{eq: pb verifie par S V}}
\left\{\begin{array}{r@{\;=\;}ll}
\operatorname{div}(\mu  \nabla  S^{\Vcoin} ) & f_{1} := F_{\mathrm{v}} + \khii \varphi  - \operatorname{div}\!\big(\mu  \nabla  (\khii \cdot S^{\mathcal{A}} )\big)&\text{ in }\Omega _{1} \setminus (\Gamma \cap \{X>\Rcoin\})  \\[0.5mm]
[\mu  \partial _{Y} S^{\Vcoin} ]_{\Gamma } & g_{1} := -[\mu \partial _{Y} (\khii S^{\mathcal{A}} )]_{\Gamma } &\text{ on }\Gamma \cap \{X>\Rcoin\} \\[0.5mm]
S^{\Vcoin} &0 &\text{ on }\partial \Omega _{1}
\end{array}\right.\end{equation}
has a unique solution in $V$ and that $S = S^{V} + \khii S^{\mathcal{A}}$. In addition, it also gives that $(\mu  \Delta  S^{\mathcal{A}} - \varphi )_{|D} \in  \sum _{d'<d} \mathcal{A}_{d'} (D)$ for any $D\in \{\Omega ,\couche\}$, and $[\mu \partial _{Y} S^{\mathcal{A}} ]_{\Gamma } \in  \sum _{d'<d} \mathcal{A}_{d'} (\Gamma )$. And we have $F_{\mathrm{v}} =o_{\partial } (\bfr^{d} )$ by hypothesis, so for $d$ small enough and up to increasing $\Rcoin$ (without loss of generality): $\forall m\in \mathbb{N},$
\begin{equation}{\label{eq: regularite de S V: snd membres et Kondratev}}
(\sbd{\chi  f_{1} })_{|\sbd \Omega } \in K^{\infty } _{-2} (\mathbb{R},H^{m} (0,\Theta )) \quad  \text{ and } \quad  (\sbd{\chi  f_{1} })_{|\sbd \couche} \in K^{\infty } _{-3/2} (\mathbb{R},H^{m} (-1,0)) \quad  \text{ and } \quad  \sbd{\chi  g_{1} }\in K^{\infty } _{-1/2} (\sbd \Gamma )\end{equation}

\underline{Step 2:} For any $\beta \in \mathbb{R}$, we denote $K^{0} _{\beta } (\sbd \Omega ):=K^{0} _{\beta } (\mathbb{R},L^2 (0,\Theta )) = \{(t,\theta )\mapsto  e^{\beta t} u(t,\theta ) \mid u\in L^2 (\sbd \Omega )\}$ and $K^{0} _{\beta } (\sbd\couche) :=K^{0} _{\beta } (\mathbb{R},L^2 (-1,0))$. And we define \[\sbd \Vcoin := \{w\in H^{1} _{\mathrm{loc}} (\sbd \Pi ) \mid \nabla  w_{|\sbd \Omega } \in L^2 (\sbd \Omega ),\; \partial _{t} w_{|\sbd\couche } \in K^{0} _{1/2} (\sbd\couche ),\; \partial _{\theta } w_{|\sbd\couche } \in K^{0} _{-1/2} (\sbd\couche )\, \text{ and } \,w_{|\sbd\bordgauche \cup \sbd\borddroit }=0 \}.\] Note that if $\Omega _{1}$ were equal to $\Pi $, $\sbd \Vcoin$ would simply be $\{\sbd u \mid u\in \Vcoin \}$. So it is the natural variational space for Poisson's problem transferred into $\sbd \Pi $. Looking at $\chi  S^{\Vcoin}$ as a function defined on $\Pi $, we will show by induction on $n$ that:
\[\forall n\in \mathbb{N},\; \forall \Rcoin<r_{1} <r_{2,\qquad}  \partial _{t} ^{n} \sbd{\chi  S^{\Vcoin} }\in \sbd \Vcoin.\]
Since $S^{\Vcoin} \in \Vcoin$, the initial case is trivial, so only the inductive step remains to prove. Let us assume the property at rank $n$ and show it at rank $n+1$. We will use the method of finite differences. Letting $\mu :=\mu _{0}$ on $\Omega $ and $\mu :=\mu _{1}$ on $\couche$, \eqref{eq: pb verifie par S V} implies: \[\left\{\begin{array}{r@{\;=\;}ll}
	\mu  \Delta (\chi  S^{\Vcoin} ) & f_{2} := \chi  f_{1} - 2\mu  \nabla  \chi  \cdot \nabla  S^{\Vcoin} -\mu  \Delta  \chi  \cdot S^{\Vcoin} &\text{ in }D,\ \forall D\in \{\Omega ,\couche\} \\[0.5mm]
	[\mu \partial _{Y} (\chi  S^{\Vcoin} )]_{\Gamma } & g_{2} := \chi  g_{1} &\text{ on }\Gamma  \\[0.5mm]
	\chi  S^{\Vcoin} &0 &\text{ on }\bordgauche \cup \borddroit
\end{array}\right.\]
By induction hypothesis applied to $(r_{1} ',r_{2} '):= (\frac{\Rcoin+r_{1} }{2} ,r_{1} )$, $\partial _{t} ^{k} \nabla  S^{\Vcoin}$ and $\partial _{t} ^{k} S^{\Vcoin}$ are $L^2 $ on $\{r_{1} <r<r_{2} \}$ for any $k\in [\![0,n]\!]$. So \eqref{eq: regularite de S V: snd membres et Kondratev} implies $\sbd{f_{2} }_{|\sbd \Omega } \in K^{n} _{-2} (\mathbb{R},L^2 (0,\Theta ))$ and $\sbd{f_{2} }_{|\sbd \couche} \in K^{n} _{-3/2} (\mathbb{R},L^2 (-1,0))$. Similarly $\sbd{g_{2} }\in K^{\infty } _{-1/2} (\sbd \Gamma )$. Next, we apply the change variables $(x,y) \rightsquigarrow (t,\theta )$ using \eqref{eq: transfert du probleme de Poisson de Pi a sbdPi}, and then $\partial _{t} ^{n}$. We get that $s:=\partial _{t} ^{n} (\sbd{\chi  S^{\Vcoin} })$ satisfies: 
\begin{equation}{\label{eq: regularite de S, equations de s}}
\hspace{-5mm}\left\{\begin{array}{@{\;}r@{\;=\;}l@{\;\,}l@{}}
	-e^{-2t} \mu _{0} \Delta s & f_{\Omega } :=e^{-2t} \partial _{t} ^{n} (e^{2t} \sbd{f_{2} }_{|\Omega } ) &\text{ in }\sbd \Omega  \\[0.5mm]
	-\mu _{1} \,(e^{-2t} (\partial _{t} ^2 -\partial _{t} )+\partial _{\theta } ^2 )s & f_{\couche} :=\partial _{t} ^{n} \sbd{f_{2} }_{|\couche} + \sum\limits _{k=0} ^{n-1} (-2)^{n-k} e^{-2t} \mu _{1} (\partial _{t} ^2 -\partial _{t} ) \partial _{t} ^{k} (\sbd{\chi  S^{V} }) &\text{ in }\sbd\couche \\[0.5mm]
	s_{|\theta =0^{+}} - s_{|\theta =0^{-}} & 0 &\text{ on }\sbd \Gamma  \\[0.5mm]
	e^{-t} \mu _{0} \partial _{\theta } s_{|\theta =0^{+}} - \mu _{1} \partial _{\theta } s_{|\theta =0^{-}} & g :=\partial _{t} ^{n} \sbd{g_{2} } &\text{ on }\sbd \Gamma  \\[0.5mm]
	s & 0 &\text{ on }\sbd\bordgauche \cup \sbd\borddroit \hspace{-7mm}
\end{array}\right. \end{equation}
with $f_{\Omega } \in K^{0} _{-2} (\Omega )$, $f_{\couche} \in K^{0} _{-3/2} (\sbd\couche )$ and $g\in K^{\infty } _{-1/2} (\sbd \Gamma )$ by induction hypothesis. The variational formulation of \eqref{eq: regularite de S, equations de s} is: $\forall \varphi \in \sbd \Vcoin ,$
\[ \int _{\sbd \Omega } \mu _{0} \nabla  s \cdot \nabla  \varphi +\int _{\sbd \couche} \mu _{1} \,(e^{-t} \,\partial _{t} s \,\partial _{t} \varphi + e^{t} \,\partial _{\theta } s \,\partial _{\theta } \varphi ) = \int _{\sbd \Omega } e^{2t} f_{\Omega } \varphi + \int _{\sbd\couche} e^{t} f_{\couche} \varphi + \int _{\sbd \Gamma } e^{t} g \varphi .\]
Let us denote $D_{\eta } \varphi (t,\theta ) := \frac{\varphi (t+\eta ,\theta )-\varphi (t,\theta )}{\eta }$ for any $\eta \in \mathbb{R}^*$ and any function $\varphi $. Taking $\varphi  := D_{-\eta } D_{\eta } \bar{s}$ and discretely integrating by parts $D_{-\eta }$, we have:
\begin{multline*}
	\int _{\sbd \Omega } \mu _{0} |D_{\eta } \nabla  s|^2  +\int _{\sbd\couche} \mu _{1} \,(D_{\eta } (e^{-t} \partial _{t} s) \cdot D_{\eta } \partial _{t} \bar{s}+ e^{t} \,|D_{\eta } \partial _{\theta } s|^2 ) \\
	= \int _{\sbd \Omega } e^{2t} f_{\Omega } \cdot D_{-\eta } D_{\eta } \bar{s}
	+ \int _{\sbd\couche} e^{t} f_{\couche} \cdot D_{-\eta } D_{\eta } \bar{s}
	+ \int _{\sbd \Gamma } D_{-\eta } D_{\eta } (e^{t} g) \cdot  \bar{s}
\end{multline*}
But for any functions $\varphi ,\psi $, we have $D_{\eta } (\varphi  \psi )=\varphi  \cdot  D_{\eta } \psi  + D_{\eta } \varphi  \cdot \psi $. Therefore:
\begin{itemize}
	\item $D_{\eta } (e^{-t} \partial _{t} s) \cdot D_{\eta } \partial _{t} \bar{s}=e^{-t} |D_{\eta } \partial _{t} s|^2  + \frac{e^{-\eta } -1}{\eta } e^{-t} \partial _{t} s\cdot D_{\eta } \partial _{t} \bar{s}$
	\item $D_{-\eta } D_{\eta } (e^{t} g) =e^{t} D_{-\eta } D_{\eta } g + 2 \frac{e^{\eta } -1}{\eta } e^{t} D_{\eta } g + \big(\frac{e^{\eta } -1}{\eta } \big)^2  e^{t} g$
\end{itemize}
Let us assume that $\eta $ is small enough so that $\frac{|e^{\eta } -1|}{\eta } <2$ and $\frac{|e^{-\eta } -1|}{\eta } <2$. Then:\vspace{-2mm}
\begin{multline*}
	\|D_{\eta } \nabla  s\|_{L^2 (\sbd \Omega )}^2  + \|D_{\eta } \partial _{t} s\|_{K^{0} _{1/2} (\sbd\couche )}^2  +\|D_{\eta } \partial _{\theta } s\|_{K^{0} _{-1/2} (\sbd\couche )}^2  + \int _{\sbd\couche} \mu _{1} {\textstyle \frac{e^{-\eta } -1}{\eta } }e^{-t} \partial _{t} s\cdot D_{\eta } \partial _{t} \bar{s}\\
	\lesssim  \|f_{\Omega } \|_{K^{0} _{-2} (\sbd \Omega )} \|D_{-\eta } D_{\eta } s\|_{L^2 (\sbd \Omega )} 
	+ \|f_{\couche} \|_{K^{0} _{-3/2} (\sbd\couche )} \|D_{-\eta } D_{\eta } s\|_{K^{0} _{1/2} (\sbd\couche )}
	+ \|g\|_{K^{2} _{-1/2} (\sbd \Gamma )} \|s\|_{K^{0} _{-1/2} (\sbd \Gamma )}
\end{multline*}
Then, moving the intergral to the right-hand side and majorizing some $\|D_{\eta } \cdot \|$ by $\|\partial _{t} \cdot \|$ or $\|\nabla  \cdot \|$, we get:
\begin{multline*}
	\|D_{\eta } \nabla  s\|_{L^2 (\sbd \Omega )}^2  + \|D_{\eta } \partial _{t} s\|_{K^{0} _{1/2} (\sbd\couche )}^2  +\|D_{\eta } \partial _{\theta } s\|_{K^{0} _{-1/2} (\sbd\couche )}^2  
	\lesssim  \|f_{\Omega } \|_{K^{0} _{-2} (\sbd \Omega )} \|D_{\eta } \nabla s\|_{L^2 (\sbd \Omega )} \\
	+ \|f_{\couche} \|_{K^{0} _{-3/2} (\sbd\couche )} \|D_{\eta } \partial _{t} s\|_{K^{0} _{1/2} (\sbd\couche )} 
	+ \|g\|_{K^{2} _{-1/2} (\sbd \Gamma )} \|s\|_{K^{0} _{-1/2} (\sbd \Gamma )}
	+ \|\partial _{t} s\|_{K^{0} _{1/2} (\sbd\couche )} \|D_{\eta } \partial _{t} s\|_{K^{0} _{1/2} (\sbd\couche )}
\end{multline*}
But, since $s_{|\sbd\borddroit} =0$, a Poincaré-type inequality gives $\|s\|_{K^{0} _{-1/2} (\sbd \Gamma )}^2  \lesssim  \|\partial _{\theta } s\|_{K^{0} _{-1/2} (\sbd\couche )} ^2 $. Finally, using Young's inequality $ab\lesssim  \frac{1}{\delta } a +\delta  b$ on norm products with $\delta $ small enough and moving the $\|D_{\eta } \partial _{t} s\|_{K^{0} _{1/2} (\sbd\couche )}^2 $ and $\|D_{\eta } \partial _{\theta } s\|_{K^{0} _{-1/2} (\sbd\couche )}^2 $ from the right-hand side to the left-hand one, we get:
\begin{align*}
	\|\nabla  \partial _{t} s\|_{L^2 (\sbd \Omega )}^2  &+ \|\partial _{t} ^2  s\|_{K^{0} _{1/2} (\sbd\couche )}^2  +\|\partial _{\theta } \partial _{t} s\|_{K^{0} _{-1/2} (\sbd\couche )}^2  \\
	&\leqslant  \limsup _{\eta \rightarrow 0} \|D_{\eta } \nabla  s\|_{L^2 (\sbd \Omega )}^2  + \|D_{\eta } \partial _{t} s\|_{K^{0} _{1/2} (\sbd\couche )}^2  +\|D_{\eta } \partial _{\theta } s\|_{K^{0} _{-1/2} (\sbd\couche )}^2  \\
	&\lesssim  \|f_{\Omega } \|_{K^{0} _{-2} (\sbd \Omega )}^2  + \|f_{\couche} \|_{K^{0} _{-3/2} (\sbd\couche )}^2  + \|g\|_{K^{1} _{-1/2} (\sbd \Gamma )}^2  + \|\partial _{\theta } s\|_{K^{0} _{-1/2} (\sbd\couche )}^2  + \|\partial _{t} s\|_{K^{0} _{1/2} (\sbd\couche )}^2  \\
	&<\infty 
\end{align*}
By definition of $\sbd \Vcoin$, it implies that $\partial _{t} s\in \sbd \Vcoin$ and completes the induction.\medskip\smallskip

\underline{Step 3:} Let $K^{0} _{1/2} (\sbd \Pi ) := K^{0} _{1/2} (\mathbb{R},L^2 (-1,\Theta ))$. For now we have proven that, for any $(n,i)\in (\mathbb{N}\times \{0,1\}) \setminus (\{(0,0)\})$ and $\Rcoin<r_{1} <r_{2}$, $\partial _{t} ^{n} \partial _{\theta } ^{i} \sbd{\chi  S^{\Vcoin} } \in  K^{0} _{1/2} (\sbd \Pi )$. It generalizes to the case $(n,i)=(0,0)$ thanks to a Poincaré inequality: $\|\sbd{\chi  S^{\Vcoin} }\|_{K^{0} _{1/2} (\sbd \Pi )} \lesssim  \|\partial _{\theta } \sbd{\chi  S^{\Vcoin} }\|_{K^{0} _{1/2} (\sbd \Pi )} \lesssim  \|\partial _{\theta } \sbd{\chi  S^{\Vcoin} }\|_{L^2 (\sbd \Omega )} + \|\partial _{\theta } \sbd{\chi  S^{\Vcoin} }\|_{K^{0} _{-1/2} (\sbd\couche)} < \infty $.\\
To treat higher-order $\theta $-derivatives, we start from the equality $\mu  \Delta (\chi  S^{\Vcoin} ) = \chi  f_{1} - 2\mu  \nabla  \chi  \cdot \nabla  S^{\Vcoin} -\mu  \Delta  \chi  \cdot S^{\Vcoin}$ in $\Omega \cup \couche$ proven in to step 2. Applying the change of variables $(x,y) \rightsquigarrow (t,\theta )$ gives by \eqref{eq: transfert du probleme de Poisson de Pi a sbdPi}: 
\begin{equation}{\label{eq: regularite de S, etape 3}}
\left\{\begin{array}{r@{\;=\;}ll}
	-e^{-2t} \mu _{0} \Delta \sbd{\chi  S^{\Vcoin} } & 
	\sbd{\chi  f_{1} } - 2\mu _{0} \sbd{\nabla  \chi  \cdot \nabla  S^{\Vcoin} } -\mu _{0} \sbd{\Delta  \chi  \cdot S^{\Vcoin} } &\text{ in }\sbd \Omega  \\[0.5mm]
	-\mu _{1} \,(e^{-2t} (\partial _{t} ^2 -\partial _{t} )+\partial _{\theta } ^2 )(\sbd{\chi  S^{\Vcoin} }) & 
	\sbd{\chi  f_{1} } - 2\mu _{1} \sbd{\nabla  \chi  \cdot \nabla  S^{\Vcoin} } -\mu _{1} \sbd{\Delta  \chi  \cdot S^{\Vcoin} } &\text{ in }\sbd\couche 
\end{array}\right.\end{equation}
Moreover, \eqref{eq: regularite de S V: snd membres et Kondratev} implies that, for any $m\in \mathbb{N}$, $\sbd{\chi  f_{1} }_{|\sbd \Omega } \in K^{\infty } _{-3/2} (\mathbb{R},H^{m} (0,\Theta ))$ and $\sbd{\chi  f_{1} }_{|\sbd \couche} \in K^{\infty } _{1/2} (\mathbb{R},H^{m} (-1,0))$. So deriving \eqref{eq: regularite de S, etape 3} w.r.t. $t$ and $\theta $ enough times gives by induction: $\forall (n,i)\in \mathbb{N}^2 , \forall \Rcoin<r_{1} <r_{2} ,\forall D\in \{\sbd \Omega ,\sbd\couche\},\; \partial _{t} ^{n} \partial _{\theta } ^{i} \sbd{\chi  S^{\Vcoin} }_{|D} \in  K^{0} _{1/2} (D)$. Finally, applying this to $\chi :=\khii$ completes the proof.
\end{demo}

\newcommand{\dTronc}[1]{d_{#1}}
\newcommand{\Sreste}[1]{S_{\mathrm{r},#1}}
\newcommand{\SAbarre}[1]{S^\infty_{#1}}
\newcommand{\SresteTilde}{\tilde{S}_{\mathrm{r},n+1}}
\newcommand{\nvCoefs}{\mathfrak{D}}
\phantomsection\label{preuve: DA des champs de coin 2}
\begin{demo}[Proof of Theorem~\ref{DA des champs de coin}]
For any $q\in \TN$, let $\sigma _{q} (\SDA):=\sigma _{q} (S)$. We denote, for any $d\in \mathbb{R}$, $T_{>d} :=  \mathrm{id}-T_{\leqslant d}$ and
\[\SAbarre{d} := \left\langle  -\RcouchePi\circ \partial _{x|\couche} ^2 ,\; \frac{\mu _{0} }{\mu _{1}} \RGammaPi\circ \partial _{y|\Gamma ,y=0^{+}} \right\rangle  \bigg( \frac{1}{\mu _{0}} \ROmegaPi (F^{\infty } _{\Omega } )+ \frac{1}{\mu _{1}} \RcouchePi (F^{\infty } _{\couche} ) +\! \fsum_{q\in \TZ*\cap (d,\infty )} \! \sigma _{q} (\SDA) \,\phi_{q} \bigg) \]
and, for any $n\in \mathbb{N}$, $\dTronc{n}:= \frac{1}{2} -n$ and $\Sreste{n} := \khii S- \khii T_{>\dTronc{n}} (\SAbarre{\dTronc{n}})$. Let $m\in \mathbb{N}$. We will show by induction on $n$ that there are coefficients $(\sigma _{q} (\SDA))_{q\in -\TN}$ s.t.: \[\forall n\in \mathbb{N},\; \forall \beta >\dTronc{n},\qquad  \sbd{\Sreste{n} }\in  K^{\infty } _{\beta } (\mathbb{R},\mathcal{H}^{m} ).\]
For now we only know $\sigma _{q} (\SDA)$ for $q\geqslant  \frac{\pi }{\Theta } >d_{0}$, and $\Sreste{n}$ involves $\sigma _{q} (\SDA)$ only when $q>\dTronc{n}$, so the $n$-th inductive step involves fixing $\sigma _{q} (\SDA)$ for all $q\in \TZ*\cap (\dTronc{n},\dTronc{n-1}]$.

\begin{description}
\item[Initial case:] By Proposition~\ref{regularite de S variationnel}, there is $d\in \mathbb{R}_-$ and $S^{\Vcoin}$ s.t. $S=S^{\Vcoin} +\khii T_{\geqslant d} (\SAbarre{\dTronc{0}})$ and $\sbd{\khii S^{\Vcoin} } \in K^{\infty } _{1/2} (\mathbb{R},\mathcal{H}^{m} )$. Hence
\begin{align*}
	\Sreste{0} &= \khii S- \khii T_{>\dTronc{0}} (\SAbarre{\dTronc{0}})\\
	&= \khii \big(S^{\Vcoin} + \khii T_{\geqslant d} (\SAbarre{\dTronc{0}}) \big) - \khii T_{>\dTronc{0}} (\SAbarre{\dTronc{0}})\\
	&= \khii S^{\Vcoin} + \khii \cdot (\khii -1)\cdot T_{\geqslant d} (\SAbarre{\dTronc{0}}) +\,\khii \cdot \big(T_{\geqslant d} (\SAbarre{\dTronc{0}}) -T_{>0} (\SAbarre{\dTronc{0}}) \big)
\end{align*}
So $\sbd{\Sreste{0}}\in K^{\infty } _{1/2} (\mathbb{R},\mathcal{H}^{m} )$. But $\sbd{\Sreste{0}}$ is null in a vicinity of $-\infty $, so: $\forall \beta >d_{0} = \frac{1}{2} ,\; \sbd{\Sreste{0}}\in K^{\infty } _{\beta } (\mathbb{R},\mathcal{H}^{m} )$.

\item[Inductive step:] We assume the property at rank $n$ and will show it at rank $n+1$. Let $\SresteTilde:= \khii S- \khii T_{>\dTronc{n+1}} (\SAbarre{\dTronc{n}})$ and $\nvCoefs:=\TZ*\cap (\dTronc{n+1},\dTronc{n}]$. $\SresteTilde$ is a variant of $\Sreste{n+1}$ that does not involve $\sigma _{q} (\SDA)$ for $q\in \nvCoefs$. We will apply Proposition~\ref{etape de base du DA des champs de coin} to $d:=\dTronc{n}$ and $s:=\sbd{\SresteTilde}$. To do so, we must check that: \[\left\{\begin{array}{r@{\;\,\in\;}l}
\mu _{0} \Delta  \sbd{\SresteTilde}_{|\sbd \Omega } & \bigcap\limits _{\beta >\dTronc{n+1}} \! K^{\infty } _{\beta } (\mathbb{R},H^{m-2} (0,\Theta )) \\
\mu _{1} \,(e^{-2t} (\partial _{t} ^2 -\partial _{t} )+\partial _{\theta } ^2 )\sbd{\SresteTilde}_{|\sbd\couche} & \bigcap\limits _{\beta >\dTronc{n+1}} \! K^{\infty } _{\beta } (\mathbb{R},H^{m-2} (-1,0)) \\
(e^{-t} \mu _{0} \partial _{\theta |\theta =0^{+}} - \mu _{1} \partial _{\theta |\theta =0^{-}} ) \sbd{\SresteTilde} & \bigcap\limits _{\beta >\dTronc{n+1}} \! K^{\infty } _{\beta } (\sbd \Gamma )
\end{array}\right.\]
Let us show only the first line, the others being similar.
\begin{sousdemo}
	Since $\sbd{\Delta \SresteTilde} = e^{-2t} \Delta \sbd{\SresteTilde}$, it suffices to have that: $\forall \beta >\dTronc{n+1}-2,\;\sbd{\Delta \SresteTilde} \in K^{\infty } _{\beta } (\mathbb{R},H^{m-2} (0,\Theta ))$. And this is true because we have in $\Omega $ in a vicinity of $r\rightarrow \infty $:\vspace{-1mm} \begin{align*}
		\mu _{0} \Delta \SresteTilde &= F - \mu _{0} \Delta  \big[T_{>\dTronc{n+1}} (\SAbarre{\dTronc{n}})\big] \tag*{by definition of $\SresteTilde$, and $\mu _{0} \Delta S=F\text{ in }\Omega $}\\
		&= T_{\geqslant \dTronc{n+1}-2} (F^{\infty } _{\Omega } )+o_{\partial } (r^{\dTronc{n+1}-2} ) \tag*{by hypothesis on $F$}\\
		&\ \ \ -\mu _{0} T_{>\dTronc{n+1}-2} (\Delta \SAbarre{\dTronc{n}}) \tag*{by $\deg \Delta =-2$ (Lemma~\ref{operateurs de derivation dans A})}\\
		&= (T_{\leqslant \dTronc{n+1}-2} -T_{<\dTronc{n+1}-2} )(F^{\infty } _{\Omega } ) + o_{\partial } (r^{\dTronc{n+1}-2} ) \tag*{by $\mu _{0} \Delta \SAbarre{\dTronc{n}}=F^{\infty } _{\Omega }$ (Lemma~\ref{DA algebrique des champs de coin})}
	\end{align*}\fbeq
\end{sousdemo}
Therefore, Proposition~\ref{etape de base du DA des champs de coin} states that there are coefficients $(c_{q} )_{q\in \nvCoefs}$, $s_{\mathrm{r}}$ and $d_{\mathrm{\min }} :=\min (\nvCoefs)$ s.t.: $\forall \beta \in (\dTronc{n+1},d_{\mathrm{\min }} ), \; s_{\mathrm{r}} \in K^{\infty } _{\beta } (\mathbb{R},\mathcal{H}^{m} )$ and $\sbd{\SresteTilde} = \sum _{q\in \nvCoefs} c_{q} \,\sbd{\phi_{q} } + s_{\mathrm{r}}$. Let $\sigma _{q} (\SDA) :=-c_{q}$ for any $q\in \TZ*\cap (\dTronc{n+1},\dTronc{n}]$. Then:
\[\sbd{\Sreste{n+1}} 
\;=\; \sbd{\SresteTilde} + \sum _{q\in \nvCoefs} \sigma _{q} (\SDA) \,\sbd{\khii \phi_{q} } 
\;=\; s_{\mathrm{r}} + \sum _{q\in \nvCoefs} c_{q} \,\sbd{(1-\khii )\phi_{q} } 
\;\in  \! \bigcap _{\beta \in (\dTronc{n+1},d_{\mathrm{\min }} )} \!\! K^{\infty } _{\beta } (\mathbb{R},\mathcal{H}^{m} )\]
Moreover, $\sbd{\Sreste{n+1}}$ is null in a vicinity of $-\infty $, so: $\forall \beta >\dTronc{n+1},\; \sbd{\Sreste{n+1}}\in K^{\infty } _{\beta } (\mathbb{R},\mathcal{H}^{m} )$. This concludes the induction.
\end{description}

To complete the proof of the theorem, we must show that: $\forall d\in \mathbb{R},\; S = T_{\geqslant d} (\SDA ) + o_{\partial } (\bfr^{d} )$ when $\bfr\rightarrow \infty $. We will do it in $\Omega $, but it works the same in $\couche$. Let $d\in \mathbb{R}$ and $n\in \mathbb{N}$ s.t. $\dTronc{n}<d$. We have in $\Omega $:
\[s:= \khii S -\khii T_{\geqslant d} (\SDA ) = \Sreste{n} + \khii \cdot T_{<d} \circ  T_{>\dTronc{n}} (\SDA ) \quad \text{ with }\quad  T_{<d} \circ  T_{>\dTronc{n}} (\SDA ) \in \sum\limits _{d'<d} \mathcal{A}_{d'} (\Omega ).\]
So there is $d'<d$ s.t.: $\forall m\in \mathbb{N},\; \sbd{s\,}\in K^{\infty } _{d'} (\mathbb{R},H^{m} (0,\Theta ))$. Thus for any $(i,j)\in \mathbb{N}^2 $:
\[\sbd{r^{-d'} (r\partial _{r} )^{i} \partial _{\theta } ^{j} s\,} = e^{-d't} \partial _{t} ^{i} \partial _{\theta } ^{j} \sbd{s\,} \in  H^2 (\mathbb{R},H^2 (0,\Theta )) \subset L^{\infty } (\sbd \Omega ).\]
So $r^{-d'+i} \partial _{r} ^{i} \partial _{\theta } ^{j} s$ is also bounded. By definition of $o_{\partial }$ (given in \ref{def: petit o derivable}), this concludes the proof.
\end{demo}

\subsection{Proof of Theorem~\ref{DA des champs lointains concatenes}: asymptotic behavior for far-and-layer fields-like problems}
\label{subsec: DA des champs lointains concatenes}

\begin{textesansboite}
This proof is very similar to Section~\ref{subsec: DA des champs de coin}. The main difference is that we look at the asymptotic expansion when $\bfr\rightarrow 0$ (i.e. $t\rightarrow -\infty $) instead of $\bfr\rightarrow \infty $ (i.e. $t\rightarrow +\infty $). So this time we use a truncation function $\khif \in \mathcal{C}^{\infty } (\Pi )$ that is equal to 1 in a vicinity of $\bfr=0$ and to 0 in a vicinity of infinity. In addition, we can assume that $\khif f=0$. Moreover the Laplace transform of $\sbd{\khif\bfu}$ is first defined in a left half-plane of the complex plan (instead of a right one), and then extended to the right.\\

Applying the change of variables $(x,y) \rightsquigarrow (t,\theta )$ on the equations satisfied by $\khif\bfu$ (that one can easily deduce from \eqref{eq: DA des champs lointains concatenes, hyp}) and then the Laplace transform yields that $s:=\sbd{\khif\bfu}$ satisfies:
\begin{equation}{\label{eq: pb verifie par la transformee de Laplace du champ lointain}}
\left\{\begin{array}{r@{\;=\;}ll}
	\mu _{0} (\partial _{\theta } ^2 +\lambda ^2 ) \hat{s}(\lambda ,\theta ) & \tilde{f}-\omega ^2 \rho _{0} \hat{s}(\lambda -2,\theta ) & \text{if } \ \theta \in (0,\Theta )\\[0.5mm]
	\partial _{\theta } ^2  \hat{s}(\lambda ,\theta ) & \tilde{g} & \text{if } \ \theta \in (-1,0)\\[0.5mm]
	\hat{s}(\lambda ,0^{+} ) - \hat{s}(\lambda ,0^{-}) &0& \\[0.5mm]
	\partial _{\theta } \hat{s}(\lambda ,0^{-}) & \tilde{h} & \\[0.5mm]
	\hat{s}(\lambda ,\Theta ) = \hat{s}(\lambda ,-1) & 0 & 
\end{array}\right.\end{equation}
for some functions $\tilde{f}, \tilde{g}, \tilde{h}$ depending on $g,h,\khif$. This system has the form of \eqref{eq: EDO en theta a resoudre}, so Lemma~\ref{resolution de l'EDO en theta} gives the tools to solve it and to extend $\hat{s}$ to the right (by steps of 2 here). This is stated in Proposition~\ref{etape de base du DA des champs lointains}, whose proof is very similar to Proposition~\ref{etape de base du DA des champs de coin}.
\end{textesansboite}

\begin{proposition}{}
\label{etape de base du DA des champs lointains}%
Let $d\in \mathbb{R}$, $m\in \mathbb{N} \setminus \{0,1\}$, $s\in  \bigcap\limits _{\beta <d} K^{\infty } _{\beta } (\mathbb{R},\mathcal{H}^{m} )$, $\tilde{f}\in  \!\bigcap\limits _{\beta <d+2} \! K^{\infty } _{\beta } (\mathbb{R},H^{m-2} (0,\Theta ))$, $\tilde{g}\in  \!\bigcap\limits _{\beta <d+2} \! K^{\infty } _{\beta } (\mathbb{R},H^{m-2} (-1,0))$, $\tilde{h}\in  \!\bigcap\limits _{\beta <d+2} \! K^{\infty } _{\beta } (\sbd \Gamma )$ and $d_{\mathrm{\max }} = \max (\TZ*\cap [d,d+2))$. We assume: \vspace{-2mm}
\[\left\{\begin{array}{rl}
	\mu _{0} \Delta  s+e^{2t} \omega ^2 \rho _{0} s = \tilde{f} &\text{ in }\sbd \Omega \\[0.5mm]
	\partial _{\theta } ^2 s = \tilde{g}  &\text{ in }\sbd\couche \\[0.5mm]
	s_{|\theta =0^{+}} - s_{|\theta =0^{-}} =0  &\text{ on }\sbd \Gamma  \\[0.5mm]
	\mu _{1} \partial _{\theta } s_{|\theta =0^{-}} = \tilde{h} &\text{ on }\sbd \Gamma  \\[0.5mm]
	s =0 &\text{ on }\sbd\bordgauche \cup \sbd\borddroit
\end{array}\right.\] Then there is $s_{\mathrm{r}} \in  \!\bigcap\limits _{d_{\mathrm{\max }} <\beta <d+2}  \! K^{\infty } _{\beta } (\mathbb{R},\mathcal{H}^{m} )$ and $(c_{q} )\in \mathbb{C}^{\TZ*\cap [d,d+2)}$ s.t. $\displaystyle s = \sum _{q\in \TZ*\cap [d,d+2)} c_{q} \,\sbd{\phi_{q} } + s_{\mathrm{r}}$.\fbeq[0.5]
\end{proposition}

\begin{paragraphesansboite}{Assumption \nvnumpar{}:}
\label{hyp: u a support compact}%
Replacing $\bfu$ by $\khii \bfu$, we assume without loss of generality that $\bfu$ has a compact support.
\end{paragraphesansboite}

\begin{textesansboite}
For any $q\in -\TN$, let $\sigma _{q} (\bfu) := \sigma _{q} (\bfu_{|\Omega } )$, where $\sigma _{q} (\bfu_{|\Omega } )$ is set in Definition~\ref{def: sigma dans H1 + khi A}. This quantity vanishes when $q$ is small enough. Proposition~\ref{regularite de u variationnel} is the analogue of Proposition~\ref{regularite de S variationnel}.
\end{textesansboite}

\begin{proposition}{}
\label{regularite de u variationnel}%
Let $d\in \mathbb{R}_+$,\vspace{-1mm}
\[\bfu^{\mathcal{A}} := T_{\leqslant d} \Bigg[ \big\langle  {-}k_{0} ^2 \,\ROmegaPi \big\rangle  \bigg( \frac{1}{\mu _{1}} \RcouchePi(g^{0} ) + \frac{1}{\mu _{1}} \RGammaPi(h^{0} ) + \fsum_{q\in -\TN} \sigma _{q} (\bfu) \,\phi_{q} \bigg)\Bigg] \vspace{-2mm}\]
 in $\mathcal{A}(\Pi )$ and $\bfu^{\mathrm{v}} := \bfu-\khif \bfu^{\mathcal{A}}$. There is $d\in \mathbb{R}_+$ s.t. $\bfu^{\mathrm{v}} _{|\Omega } \in  H^{1} (\Omega )$ and: $\forall m\in \mathbb{N},\;\sbd{\khif \bfu^{\mathrm{v}} } \in K^{\infty } _{-1} (\mathbb{R},\mathcal{H}^{m} )$.
\end{proposition}

\newcommand{\uAbarreRef}{\bfu^{\Abarre}_-}
\newcommand{\dregu}{0}
\begin{demo}
We consider any $d\in \mathbb{R}_+$ and we will fix it later. The proof has five steps: writing the equations satisfied by $\bfu^{\mathrm{v}}$, and then showing that $\bfu^{\mathrm{v}}$ is regular in the layer, that it is $H^{1}$ in $\Omega $, that it is regular in $\Omega $ w.r.t $t$, and that it is regular in $\Omega $ w.r.t. $\theta $. For any $r_{1} ,r_{2} \in \mathbb{R}$ s.t. $r_{1} <r_{2} <r_{f}$, we denote $\chi \in \mathcal{C}^{\infty } (\mathbb{R}^2 )$ a radial function equal to 1 on $B(0,r_{1} )$ and 0 on $\mathbb{R}^2  \setminus B(0,r_{2} )$ ($r_{1} ,r_{2} $ are implicit in this notation).\\

\underline{Step 1:} Similarly to the proofs of theorems~\ref{cadre fonctionnel des champs lointains} and \ref{cadre fonctionnel pour les champs de coin} and Proposition~\ref{regularite de S variationnel}, one can show that:
\[\left\{\begin{array}{@{\,}r@{\;=\;}l@{\,}l@{}}
	(\mu _{0} \Delta  +\omega ^2 \rho _{0} )\bfu^{\mathrm{v}} & f_{1} := f-(\mu _{0} \Delta  +\omega ^2 \rho _{0} )(\khif \bfu^{\mathcal{A}} ) &\text{ in }\Omega  \\[0.5mm]
	\mu _{1} \partial _{Y} ^2 \bfu^{\mathrm{v}} & g_{1} :=g-\mu _{1} \partial _{Y} ^2 (\khif \bfu^{\mathcal{A}} ) &\text{ in }\couche \\[0.5mm]
	\bfu^{\mathrm{v}} _{|y=0^{+}} - \bfu^{\mathrm{v}} _{|y=0^{-}} & 0 &\text{ on }\Gamma  \\[1mm]
	\mu _{1} \partial _{Y} \bfu^{\mathrm{v}} _{|Y=0^{-}} & h_{1} :=h-\mu _{1} \partial _{Y} (\khif \bfu^{\mathcal{A}} ) &\text{ on }\Gamma  \\[1mm]
	\bfu^{\mathrm{v}} & 0 &\text{ on }\bordgauche \cup \borddroit
\end{array}\right. \mathrm{with} \ \; \left\{\begin{array}{@{}r@{\;\in \;}l@{}}
	(\mu _{0} \Delta  +\omega ^2 \rho _{0} )\bfu^{\mathcal{A}} & \sum\limits _{d'>d} \! \mathcal{A}_{d'} (\Omega ) \\[0.5mm]
	\partial _{Y} ^2 \bfu^{\mathcal{A}} _{|\couche} - T_{\leqslant d} (g^{0} ) & \sum\limits _{d'>d} \! \mathcal{A}_{d'} (\couche)  \\[0.5mm]
	\partial _{Y} \bfu^{\mathcal{A}} _{|Y=0^{-}} - T_{\leqslant d} (h^{0} ) & \sum\limits _{d'>d} \! \mathcal{A}_{d'} (\Gamma )
\end{array}\right.\]
Moreover, given that $\chi  f=0$ for any $r_{1} <r_{2} <r_{f}$, one can check that for $d$ big enough and up to decreasing $r_{f}$ (without loss of generality), we have $\sbd{\chi  f_{1}} \in K^{\infty } _{\dregu} (\mathbb{R},H^{m} (0,\Theta ))$, $\sbd{\chi  g_{1}} \in K^{\infty } _{\dregu} (\mathbb{R},H^{m} (-1,0))$ and $\sbd{\chi  h_{1}} \in K^{\infty } _{\dregu} (\sbd \Gamma )$ for any $m\in \mathbb{N}$.\\

\underline{Step 2:} Let us show that: $\forall r_{1} <r_{2} <r_{f} , \forall m\in \mathbb{N},\; \sbd{\chi  \bfu^{\mathrm{v}} _{|\couche} } \in K^{\infty } _{0} (\mathbb{R},H^{m} (-1,0))$. By assumption, for any $x>0$, $\bfu^{\mathrm{v}}$ in $H^{1}$ in a vicinity of $\{x\}\times [-1,0] \subset \couche$, so $\bfu^{\mathrm{v}} _{|\{x\}\times [-1,0]}$ is well-defined in $L^2 (\{x\}\times [-1,0])$. Moreover, $g_{1}$ and $h_{1}$ are $o_{\partial } (x^{d} )$, so they are differentiable on $(0,r_{f} )\times [-1,0]$ (up to decreasing $r_{f}$). Hence, $Y\in (-1,0)\mapsto  \bfu^{\mathrm{v}} (x,Y)$ is $\mathcal{C}^2 $ and the ODE it satisfies (see step 1) is explicitly solvable:
\[\forall (x,Y)\in (0,r_{f} )\times [-1,0], \qquad  \bfu^{\mathrm{v}} (x,Y) = \int _{-1} ^{0} \big( (Y-Y')^{+} -Y-1 \big)\, g_{1} (x,Y')\,\mathrm{d}Y' +h_{1} (x)\cdot (Y+1).\]
Since $\sbd{\chi  g_{1}} \in K^{\infty } _{\dregu} (\mathbb{R},H^{m} (-1,0))$ for any $m$ and $\sbd{\chi  h_{1}} \in K^{\infty } _{\dregu} (\sbd \Gamma )$, we deduce $\sbd{\chi  \bfu^{\mathrm{v}} _{|\couche} } \in K^{\infty } _{\dregu} (\mathbb{R},H^{m} (-1,0))$.\\

\underline{Step 3:} Let us show that $\bfu^{\mathrm{v}} _{|\Omega } \in H^{1} (\Omega )$. First, step 2 implies $\sbd{\chi  \bfu^{\mathrm{v}} _{|\Gamma } } \in K^{\infty } _{\dregu} (\sbd \Gamma )$ so $\chi  \bfu^{\mathrm{v}} _{|\Gamma } \in H^{1/2} (\Gamma )$. Since $\bfu \in H^{1} _{\mathrm{loc}}$ and $\bfu$ has a compact support by assumption~\ref{hyp: u a support compact}, we get $\bfu^{\mathrm{v}} _{|\Gamma } \in H^{1/2} (\Gamma )$. So the following system has a solution in $H^{1} (\Omega )$ by Lemma~\ref{pb de Helmholtz variationnel}: \[\left\{\begin{array}{r@{\;=\;}ll}
	\mu _{0} \Delta  u^{\mathrm{v}} +\omega ^2 \rho _{0} u^{\mathrm{v}} & f_{1} &\text{ in }\Omega  \\[0.5mm]
	u^{\mathrm{v}} _{|y=0^{+}} & \bfu^{\mathrm{v}} _{|y=0^{-}} &\text{ on }\Gamma  \\[1mm]
	u^{\mathrm{v}} & 0 &\text{ on }\bordgauche
\end{array}\right.\]
But $\bfu^{\mathrm{v}}$ also satisfies these equations, so $u^{\mathrm{v}} +\khif \bfu^{\mathcal{A}} _{|\Omega }$ satisfies the same problem as $\bfu_{|\Omega }$, and they both belong to $H^{1} (\Omega )+\khiz \mathcal{A}(\Omega )$. Then by uniqueness in Theorem~\ref{cadre fonctionnel des champs lointains}, $u^{\mathrm{v}} +\khif \bfu^{\mathcal{A}} _{|\Omega } =\bfu_{|\Omega }$, i.e. $u^{\mathrm{v}} = \bfu^{\mathrm{v}} _{|\Omega }$. Thus, $\bfu^{\mathrm{v}} _{|\Omega } \in H^{1} (\Omega )$.\\

\underline{Step 4:} Let $v:\Omega \rightarrow \mathbb{C}$ defined by $v := \bfu^{\mathrm{v}} _{|\Omega } - v_{\mathrm{lift}}$ with $v_{\mathrm{lift}} (r,\theta ) := \big(1- \frac{\theta }{\Theta } \big) \bfu^{\mathrm{v}} _{|\Gamma } (r)$ in polar coordinates. Let $K^{0} _{\beta } (\sbd \Omega ):=K^{0} _{\beta } (\mathbb{R},L^2 (0,\Theta ))$ for any $\beta \in \mathbb{R}$. We will show by induction on $n$ that: \[\forall n\in \mathbb{N}^*,\; \forall r_{1} <r_{2} <r_{f} ,\qquad  \partial _{t} ^{n} \sbd{\chi \, v}\in K^{0} _{-1} (\sbd \Omega ) \quad  \text{ and } \quad  \nabla  \partial _{t} ^{n} \sbd{\chi \, v}\in L^2 (\sbd \Omega ).\]
For the initial case, we have on the one hand $\chi \bfu^{\mathrm{v}} _{|\Omega } \in H^{1} (\Omega )$ so $\sbd{\chi \bfu^{\mathrm{v}} _{|\Omega } }\in K^{0} _{-1} (\sbd \Omega )$ and $\nabla  \sbd{\chi \bfu^{\mathrm{v}} _{|\Omega } }\in L^2 (\sbd \Omega )$, and on the other $\sbd{\chi  \bfu^{\mathrm{v}} _{|\Gamma } }\in K^{1} _{\dregu} (\sbd \Gamma )$ so $\sbd{\chi  v_{\mathrm{lift}} }\in K^{0} _{-1} (\sbd \Omega )$ and $\nabla  \sbd{\chi  v_{\mathrm{lift}} }\in L^2 (\sbd \Omega )$. Thus the initial case is proven and only the inductive step remains to prove. We assume the property at rank $n$ and we will show it at rank $n+1$, using the method of finite differences. Step 1 implies that: \[\left\{\begin{array}{r@{\;=\;}ll}
	(\mu _{0} \Delta +\omega ^2 \rho _{0} )(\chi  v) & f_{2} := \chi  f_{1} -\chi \,\mu _{0} \Delta  v_{\mathrm{lift}} - 2\mu _{0} \nabla  \chi  \cdot \nabla  \bfu^{\mathrm{v}} -\mu _{0} \Delta  \chi  \cdot \bfu^{\mathrm{v}} &\text{ in }\Omega  \\[0.5mm]
	\chi  v &0 &\text{ on }\bordgauche \cup \Gamma 
\end{array}\right.\]
In addition, we have $\sbd{f_{2} } \in K^{n} _{-2} (\mathbb{R},L^2 (0,\Theta ))$ since:
\begin{itemize}
	\item by step 1, $\sbd{\chi  f_{1}} \in K^{\infty } _{\dregu} (\mathbb{R},L^2 (0,\Theta ))$, so $\sbd{\chi  f_{1}} \in K^{\infty } _{-2} (\mathbb{R},L^2 (0,\Theta ))$,
	\item by step 2, $\sbd{\chi  \bfu^{\mathrm{v}} _{|\Gamma } }\in K^{\infty } _{\dregu} (\sbd \Gamma )$, so $\sbd{\chi  v_{\mathrm{lift}} } \in K^{\infty } _{\dregu} (\mathbb{R},H^2 (0,\Theta ))$, and applying it at $(r_{1} ',r_{2} '):= (r_{2} ,\frac{r_{2} +r_{f} }{2} )$ gives $\sbd{\chi  \Delta  v_{\mathrm{lift}} } \in K^{\infty } _{-2} (\mathbb{R},L^2 (0,\Theta ))$,
	\item and by induction hypothesis applied to $(r_{1} ',r_{2} '):= (r_{2} ,\frac{r_{2} +r_{f} }{2} )$, $\partial _{t} ^{k} \nabla  v$ and $\partial _{t} ^{k} v$ are $L^2 $ on $\{r_{1} <r<r_{2} \}$ for any $k\in [\![0,n]\!]$.
\end{itemize}
Let $s:=\partial _{t} ^{n} (\sbd{\chi  v})$. By changing variables $(x,y) \rightsquigarrow (t,\theta )$ in the previous system and applying $\partial _{t} ^{n}$, we get: \[\left\{\begin{array}{r@{\;=\;}ll}
	e^{-2t} \mu _{0} \Delta s+\omega ^2 \rho _{0} s & f_{3} := e^{-2t} \partial _{t} ^{n} (e^{2t} \sbd{f_{2} }) - \sum\limits _{k=0} ^{n-1} \omega ^2 \rho _{0} 2^{n-k} \partial _{t} ^{k} (\sbd{\chi  v}) &\text{ in }\sbd \Omega  \\[-1mm]
	s & 0 &\text{ on }\sbd\bordgauche \cup \sbd \Gamma 
\end{array}\right.\] with $f_{3} \in K^{0} _{-2} (\sbd \Omega )$ by induction hypothesis. And the variational formulation of this is: \[\forall \varphi \in K^{0} _{-1} (\sbd \Omega ), \ \nabla  \varphi \in L^2 (\sbd \Omega ) \text{ and } \varphi _{|\sbd\bordgauche \cup \sbd \Gamma }=0 \ \Rightarrow  \quad  \int _{\sbd \Omega } (\mu _{0} \nabla  s \cdot \nabla  \varphi +\omega ^2 \rho _{0} e^{2t} s \varphi ) = \int _{\sbd \Omega } e^{2t} f_{3} \varphi .\]
Let us denote $D_{\eta } \varphi (t,\theta ) := \frac{\varphi (t+\eta ,\theta )-\varphi (t,\theta )}{\eta }$ for any $\eta \in \mathbb{R}^*$ and any function $\varphi $. Taking $\varphi  := D_{-\eta } D_{\eta } \bar{s}$ and discretely integrating by parts $D_{-\eta }$ gives:
\[\int _{\sbd \Omega } (\mu _{0} |D_{\eta } \nabla  s|^2  -\omega ^2 \rho _{0} e^{2t} |s|^2 ) = \int _{\sbd \Omega } e^{2t} f_{3} \cdot  D_{-\eta } D_{\eta } \bar{s}\]
Then by coercivity (since $\operatorname{\operatorname{Im}}(\omega ) \neq 0$), we have for any $\delta >0$:
\[\|D_{\eta } \nabla  s\|_{L^2 (\sbd \Omega )}^2  + \|D_{\eta } s\|_{K^{0} _{-1} (\sbd \Omega )}^2  \lesssim  \|f_{3} \|_{K^{0} _{-2} (\sbd \Omega )} \|D_{-\eta } D_{\eta } s\|_{K^{0} _{-1} (\sbd \Omega )} 
\lesssim  {\textstyle \frac{1}{\delta } } \|f_{3} \|_{K^{0} _{-2} (\sbd \Omega )}^2  + \delta  \|D_{\eta } \partial _{t} s\|_{L^2 (\sbd \Omega )}^2 \]
Taking $\delta $ small enough and moving $\|D_{\eta } \partial _{t} s\|_{L^2 (\sbd \Omega )} ^2 $ from the right-hand side to the left-hand one, we get:
\[\|\partial _{t} \nabla  s\|_{L^2 (\sbd \Omega )}^2  + \|\partial _{t} s\|_{K^{0} _{-1} (\sbd \Omega )}^2  
\leqslant  \limsup _{\eta \rightarrow 0} \|D_{\eta } \nabla  s\|_{L^2 (\sbd \Omega )}^2  + \|D_{\eta } s\|_{K^{0} _{-1} (\sbd \Omega )}^2 
\lesssim  \|f_{3} \|_{K^{0} _{-2} (\sbd \Omega )}^2  <\infty .\]
This completes the induction.\\

\underline{Step 5:} We have shown that for any $(n,i)\in \mathbb{N}\times \{0,1\}$ and $r_{1} <r_{2} <r_{f}$, $\partial _{t} ^{n} \partial _{\theta } ^{i} \sbd{\chi  v} \in  K^{0} _{-1} (\sbd \Omega )$. Since $\sbd{\chi  v_{\mathrm{lift}} } \in  K^{\infty } _{-1} (\mathbb{R},H^{1} (0,\Theta ))$, we also have $\partial _{t} ^{n} \partial _{\theta } ^{i} \sbd{\chi  \bfu^{\mathrm{v}} _{|\Omega } } \in  K^{0} _{-1} (\sbd \Omega )$ for these $(n,i)$.\\
To treat higher-order $\theta $-derivatives, we begin with the equality $(\mu _{0} \Delta +\omega ^2 \rho _{0} )(\chi  \bfu^{\mathrm{v}} ) = \chi  f_{1} - 2\mu _{0} \nabla  \chi  \cdot \nabla  \bfu^{\mathrm{v}} -\mu _{0} \Delta  \chi  \cdot \bfu^{\mathrm{v}}$ in $\Omega $, which follows from step 1. Applying the change of variables $(x,y) \rightsquigarrow (t,\theta )$, we get:
\begin{equation}{\label{eq: regularite de u variationnel, etape 5}}
(e^{-2t} \mu _{0} \Delta +\omega ^2 \rho _{0} )(\chi  \bfu^{\mathrm{v}} ) = \sbd{\chi  f_{1} }- 2\mu _{0} \sbd{\nabla  \chi  \cdot \nabla  \bfu^{\mathrm{v}} } -\mu _{0} \sbd{\Delta  \chi  \cdot \bfu^{\mathrm{v}} }.\end{equation}
Now $\sbd{\chi  f_{1}} \in K^{\infty } _{\dregu} (\mathbb{R},H^{m} (0,\Theta ))$ for any $m\in \mathbb{N}$ by step 1. So deriving \eqref{eq: regularite de u variationnel, etape 5} w.r.t. $t$ and $\theta $ enough times gives by induction: $\forall (n,i)\in \mathbb{N}^2 , \forall r_{1} <r_{2} <r_{f} ,\; \partial _{t} ^{n} \partial _{\theta } ^{i} \sbd{\chi  \bfu^{\mathrm{v}} _{|\Omega } } \in  K^{0} _{-1} (\sbd \Omega )$. Combining it with step 2, we get $\sbd{\chi  \bfu^{\mathrm{v}} } \in K^{\infty } _{-1} (\mathbb{R},\mathcal{H}^{m} )$. Finally, we apply it to $\chi :=\khii$ to conclude.
\end{demo}

\begin{textesansboite}
Finally, Propositions~\ref{etape de base du DA des champs lointains} and \ref{regularite de u variationnel} give all the ingredients to prove Theorem~\ref{DA des champs lointains concatenes}. The proof is very similar to the one of Theorem~\ref{DA des champs de coin} on page~\pageref{preuve: DA des champs de coin 2}, so we do not go into details again.
\end{textesansboite}

\addcontentsline{toc}{section}{References}
\setcounter{biburllcpenalty}{7000}
\setcounter{biburlucpenalty}{8000}
\printbibliography

\end{document}